\documentclass[10pt]{article}

\usepackage{amsmath}
\usepackage{amsthm}
\usepackage{amssymb}
\usepackage{amscd}
\usepackage{amsfonts}  
\usepackage[all]{xy} \SelectTips{cm}{} \CompileMatrices
\usepackage{makeidx}     
\usepackage{graphicx}    
\usepackage{multicol}    

\DeclareMathAlphabet\EuR{U}{eur}{m}{n}
\SetMathAlphabet\EuR{bold}{U}{eur}{b}{n}

\makeindex             

\newcounter{commentcounter}
\newcommand{\comment}[1]                      
{\stepcounter{commentcounter}
{\bf Comment \arabic{commentcounter}}: {\ttfamily #1} }

\newcommand{\squarematrix}[4]                
{                                            
\left( \begin{array}{cc} #1 & #2 \\ #3 &
#4
\end{array} \right)
}

\newcommand{\calh}{{\cal H}}

\newcommand{\calk}{{\cal K}}

\newcommand{\bbC}{{\mathbb C}}

\newcommand{\bbQ}{{\mathbb Q}}
\newcommand{\bbR}{{\mathbb R}}

\newcommand{\bbZ}{{\mathbb Z}}

\newcommand{\bfI}{{\mathbf I}}

\newcommand{\curs}{\EuR}

\newcommand{\GROUPS}{\curs{GROUPS}}

\newcommand{\FGINJ}{\curs{FGINJ}}
\newcommand{\GRIFI}{\curs{GRIFI}}
\newcommand{\MODULES}{\curs{MODULES}}
\newcommand{\Or}{\curs{Or}}
\newcommand{\RINGS}{\curs{RINGS}}
\newcommand{\SETS}{\curs{SETS}}

\newcommand{\SPHB}{\curs{SPHB}}
\newcommand{\Sub}{\curs{Sub}}

\newcommand{\calfin}{{\mathcal F}\!{\mathcal I}\!{\mathcal N}}

\newcommand{\ch}{\operatorname{ch}}
\newcommand{\character}{\operatorname{char}}

\newcommand{\cov}{\operatorname{cov}}
\newcommand{\conhom}{\operatorname{conhom}}
\newcommand{\consub}{\operatorname{ccs}}

\newcommand{\cyclic}{\operatorname{cyclic}}
\newcommand{\degree}{\operatorname{deg}}

\newcommand{\edge}{\operatorname{edge}}

\newcommand{\Gen}{\operatorname{Gen}}
\newcommand{\ho}{\operatorname{ho}}

\newcommand{\id}{\operatorname{id}}

\newcommand{\ind}{\operatorname{ind}}

\newcommand{\inn}{\operatorname{inn}}
\newcommand{\inv}{\operatorname{inv}}
\newcommand{\inverselim}{\operatorname{invlim}}

\newcommand{\map}{\operatorname{map}}
\newcommand{\mor}{\operatorname{mor}}

\newcommand{\pr}{\operatorname{pr}}
\newcommand{\pt}{\{\bullet\}}
\newcommand{\Pic}{\operatorname{Pic}}

\newcommand{\res}{\operatorname{res}}
\newcommand{\mrf}{\operatorname{mrf}}

\newcommand{\Sw}{\operatorname{Sw}}

\newcommand{\higherlim}[3]{{\inverselim}_{#1}^{#2}#3}
\newcommand{\highercolim}[3]{{\operatorname{colim}}_{#1}^{#2}#3}
\newcommand{\invlim}[2]{\higherlim{#1}{}{#2}}
\newcommand{\colim}[2]{\highercolim{#1}{}{#2}}

\newcommand{\finsetA}{\overline{A}}
\newcommand{\finsetP}{\overline{P}}
\newcommand{\finsetcharacter}{\overline{\character}}

\newcommand{\widehatA}{\widehat{A}}
\newcommand{\widehatcharacter}{\widehat{\character}}

\newcommand{\hoA}{A_{\ho}}

\newcommand{\invA}{A_{\inv}}
\newcommand{\invP}{P_{\inv}}
\newcommand{\invcharacter}{\character_{\inv}}

\newcommand{\covA}{A_{\cov}}
\newcommand{\covP}{P_{\cov}}

\newcommand{\underlineA}{\underline{A}}
\newcommand{\underlineP}{\underline{P}}
\newcommand{\underlinecharacter}{\underline{\character}}

\newcommand{\arrowchar}{\setbox1=\hbox{\rm char}
        \setbox2=\hbox to \wd1{\leftarrowfill} \ht2=-8pt \dp2=-1pt
        \vtop{\baselineskip=-50pt\box1\box2}
        }

\newtheorem{theorem}{Theorem}[section]
\newtheorem{lemma}[theorem]{Lemma}

\newtheorem{definition}[theorem]{Definition}

\newtheorem{example}[theorem]{Example}
\newtheorem{remark}[theorem]{Remark}

\newtheorem{conjecture}[theorem]{Conjecture}

{\catcode`@=11\global\let\c@equation=\c@theorem}

\newcommand{\OrGF}[2]{\Or_{#2}(#1)}               

\newcommand{\SubGF}[2]{\Sub_{#2}(#1)}               

\newcommand{\comsquare}[8]                   
{\begin{CD}
#1 @>#2>> #3\\
@V{#4}VV @VV{#5}V\\
#6 @>>#7> #8
\end{CD}
}

\begin{document}


\hyphenation{di-men-sio-nal}
\hyphenation{equi-va-riant}
\hyphenation{ho-mo-mor-phism}
\hyphenation{geo-met-ric}
\hyphenation{Hir-ze-bruch}
\hyphenation{iso-mor-phism}
\hyphenation{ma-ni-fold}
\hyphenation{re-pre-sen-ta-tion}
\hyphenation{re-pre-sen-ta-tions}
\hyphenation{Rie-man-nian}
\hyphenation{Rie-mann}
\hyphenation{to-po-lo-gi-cal}


\typeout{----------------------------  burn.tex  ----------------------------}

\title{The Burnside Ring and Equivariant Stable Cohomotopy for Infinite Groups}
\author{Wolfgang L\"uck\thanks{\noindent email:
lueck@math.uni-muenster.de\protect\\
www: ~http://www.math.uni-muenster.de/u/lueck/\protect\\
FAX: 49 251 8338370\protect} \\
Fachbereich Mathematik\\ Universit\"at M\"unster\\
Einsteinstr.~62\\ 48149 M\"unster\\Germany}
\maketitle


\typeout{-----------------------  Abstract  ------------------------}

\begin{abstract}
After we have given a survey on the Burnside ring of a finite
group, we discuss and analyze various extensions of this notion to
infinite (discrete) groups. The first three are the
finite-$G$-set-version, the inverse-limit-version and the
covariant Burnside group. The most sophisticated one is the fourth
definition as the zero-th equivariant stable cohomotopy of the
classifying space for proper actions. In order to make sense of
this definition we define equivariant stable cohomotopy groups of finite
proper equivariant $CW$-complexes in terms of maps between the
sphere bundles associated to equivariant vector bundles. We show
that this yields an equivariant cohomology theory with a
multiplicative structure. We formulate a version of the Segal
Conjecture for infinite groups. All this is analogous and related
to the question what are the possible extensions of the notion of the
representation ring of a finite group to an infinite group. Here
possible candidates are projective class groups, Swan groups and
the equivariant topological $K$-theory of the classifying space
for proper actions.

\end{abstract}


\typeout{--------------------------------   Section 0: Introduction ------------------------------------}

\setcounter{section}{-1}
\section{Introduction}
\label{sec:Introduction}

The basic notions of the Burnside ring and of equivariant stable
cohomotopy have been defined and investigated in detail for finite
groups. The purpose of this article is to discuss how these can be
generalized to infinite (discrete) groups. The guideline will be
the related notion of the representation ring which allows several
generalizations to infinite groups, each of which reflects one
aspect of the original notion for finite groups. Analogously we
will present several possible generalizations of the Burnside ring
for finite groups to infinite (discrete) groups. There seems to be
no general answer to the question which generalization is the
right one. The answer depends on the choice of the background
problem such as universal additive properties, induction theory,
equivariant stable homotopy theory, representation theory,
completion theorems and so on. For finite groups the
representation ring and the Burnside ring are related to all these
topics simultaneously and for infinite groups the notion seems to
split up into different ones which fall together for finite groups
but not in general.

The following table summarizes in the first column the possible generalizations to infinite groups
of the representation ring $R_F(G)$  with coefficients  in a field $F$ of characteristic zero.
In the second column we list the analogous generalizations for the Burnside ring.
In the third column we give key words  for their main property, relevance or application.
Explanations will follow in the main body of the text.\\[4mm]
\renewcommand{\arraystretch}{1.5}
\begin{tabular}{||p{38mm}|p{36mm}|p{33mm}||}
\hline\hline
$R_F(G)$
&
$A(G)$
&
key words
\\
\hline \hline
$K_0(FG)$
&
$\underlineA(G)$
&
universal additive invariant, equivariant Euler characteristic
\\
\hline
$\Sw^f(G;F)$
&
$\finsetA(G)$
&
induction theory, Green functors
\\
\hline
\mbox{$R_ {\cov,F}(G):=$} $\colim{H \in \SubGF{G}{\calfin}} R_F(H)$
&
\mbox{$\covA(G):=$}   $\colim{H \in \SubGF{G}{\calfin}} A(H)$
&
collecting all values for finite subgroups with respect to induction
\\
\hline
\mbox{$R_{\inv,F}(G):=$} $\invlim{H \in \SubGF{G}{\calfin}} R_F(H)$
&
\mbox{$\invA(G):=$} $\invlim{H \in \SubGF{G}{\calfin}} A(H)$
&
collecting all values for finite subgroups with respect to restriction
\\
\hline
$K^0_G(\underline{E}G)$
&
$\hoA(G) := \pi^0_G(\underline{E}G)$
&
completion theorems, equivariant vector bundles
\\
\hline
$K_0^G(\underline{E}G)$
&
$\pi_0^G(\underline{E}G)$
& representation theo\-ry, Baum-Connes Conjecture, equivariant homotopy theory
\\
\hline \hline
\end{tabular}
\renewcommand{\arraystretch}{1.0}
\\[4mm]
The various rings are linked by the following diagram of ring homomorphisms
$$
\xymatrix{
\hoA(G) = \pi^0_G(\underline{E}G) \ar[r]^-{\edge^G}
\ar[dd]^-{\psi^0_G(\underline{E}G)}
&
\invA(G) \ar[d]^-{\invP^G}
&
\finsetA(G) \ar[l]^-{T^G} \ar[d]^-{\finsetP^G}
\\
& R_{\inv,\bbQ}(G) \ar[d]^-{c}
& \Sw^f(G;\bbQ) \ar[l]^-{S^{G,\bbQ}} \ar[d]^-{c}
\\
K_G^0(\underline{E}G)  \ar[r]^-{\edge^G}
& R_{\inv,\bbC}(G)
& \Sw^f(G;\bbC) \ar[l]^-{S^{G,\bbC}}
}
$$
where $c$ denotes the obvious change of coefficients homomorphisms
and the other maps will be  explained later.

We  will also define various pairings which are summarized in the following
diagram which reflects their compatibilities.
\begin{eqnarray}
  \xymatrix@!C=11em{
  \finsetA(G) \times \underlineA(G) \ar[r]^(.66){\mu_A^G} \ar[d]_-{\finsetP^G \times \underlineP^G}
  \ar[ddr]^(.32){T^G \times \id}|\hole
  & {\underlineA(G)} \ar[d]^-{\underlineP^G} \ar[ddr]^(.62){Q_A^G}
  \\
  \Sw^f(G;\bbQ) \times K_0(\bbQ G)  \ar[r]^(.66){\mu_K^G}  \ar@<-1ex>[ddr]_(.15){S^{G,\bbQ}}
  & K_0(\bbQ G) \ar[ddr]^(.35){Q_K^G}|\hole
  \\
  &  \invA(G) \times \underlineA(G) \ar[r]^(.66){\nu_A^G} \ar[d]^-{\invP^G \times (\covP^G \circ (V^G)^{-1})}
  & \bbZ \ar[d]^-{\id}
  \\
  &  R_{\inv,\bbQ}(G) \times  R_{\cov,\bbQ}(G) \ar[r]^-{\nu_R^G} \ar@<-1ex>[uul]_(.6){W^{G,\bbQ}}
   &
  \bbZ
 }
&&
\label{diagram_about_pairings}
\end{eqnarray}

In Section~\ref{sec:Review_of_the_Burnside_Ring_for_Finite_Groups}
we give a brief survey about the Burnside ring $A(G)$ of a finite group $G$
in order to motivate the generalizations. In
Sections~\ref{sec:The_Finite_G-Set_Version_of_the_Burnside_Ring},
\ref{sec:The_Inverse_Limit_Version_of_the_Burnside_Ring} and~
\ref{sec:The_Covariant_Burnside_Group} we treat the
finite-$G$-set-version of the Burnside Ring $\finsetA(G)$, the
inverse-limit-version of the Burnside ring $\invA(G)$ and the
covariant Burnside group $\underlineA(G)$. These definitions are
rather straightforward. The most sophisticated version of the Burnside ring
for infinite groups is the zero-th 
equivariant stable cohomotopy $\pi_G^0(\underline{E}G)$ of the
classifying space $\underline{E}G$ for proper $G$-actions. It will
be constructed in
Section~\ref{sec:Equivariant_Stable_Cohomotopy_in_Terms_of_Real_Vector_Bundles}
after we have explained the notion of an equivariant cohomology
theory with multiplicative structure in
Section~\ref{sec:Equivariant_Cohomology_Theories}. One of the main result of this paper is

{\bf Theorem~\ref{the:Equivariant_Stable_Cohomotopy_in_terms_of_equivariant_vector_bundles}}
\emph{
Equivariant stable cohomotopy $\pi^*_?$ defines an equivariant cohomology theory
with multiplicative structure for finite proper equivariant $CW$-complexes.
For every finite subgroup $H$ of the group $G$  the abelian groups $\pi^n_G(G/H)$ and
$\pi^n_H$ are isomorphic for every $n \in \bbZ$ and the rings $\pi^0_G(G/H)$ and $\pi^0_H = A(H)$ are isomorphic.
}

An important test in the
future will be whether the version of the Segal Conjecture for
infinite groups discussed in
Section~\ref{sec:The_Segal_Conjecture_for_Infinite_Groups} is
true.

The papers is organized as follows:

\begin{tabular}{lll}
\ref{sec:Review_of_the_Burnside_Ring_for_Finite_Groups}. & &
Review of the Burnside Ring for Finite Groups
\\
\ref{sec:The_Finite_G-Set_Version_of_the_Burnside_Ring}. & &
The Finite-$G$-Set-Version of the Burnside Ring
\\
\ref{sec:The_Inverse_Limit_Version_of_the_Burnside_Ring}. & &
The Inverse-Limit-Version of the Burnside Ring
\\
\ref{sec:The_Covariant_Burnside_Group}. & &
The Covariant Burnside Group
\\
\ref{sec:Equivariant_Cohomology_Theories}. & &
Equivariant Cohomology Theories
\\
\ref{sec:Equivariant_Stable_Cohomotopy_in_Terms_of_Real_Vector_Bundles}. & &
Equivariant Stable Cohomotopy in Terms of Real Vector Bundles
\\
\ref{sec:The_Homotopy_Theoretic_Burnside_ring}. & &
The Homotopy Theoretic Burnside Ring
\\
\ref{sec:The_Segal_Conjecture_for_Infinite_Groups}. & &
The Segal Conjecture for Infinite Groups
\\
& &
References
\end{tabular}


\typeout{------------Section 1: Review of the Burnside Ring for Finite Groups --------------------------}

\section{Review of the Burnside Ring for Finite Groups}
\label{sec:Review_of_the_Burnside_Ring_for_Finite_Groups}

In this section we give a brief review of the definition, properties and applications of the Burnside ring
for finite groups in order to motivate our definitions for infinite groups.

\begin{definition}{\bf (Burnside ring of a finite group).}
\label{def:A(G)_for_finite_G}
The isomorphism classes of finite $G$-sets form a commutative associative semi-ring with unit
under disjoint union and cartesian product. The \emph{Burnside ring} $A(G)$ is the
Grothendieck ring associated to this semi-ring.
\end{definition}

As abelian group the Burnside ring $A(G)$ is the free abelian group with the set $\{G/H \mid (H) \in
\consub(G)\}$ as basis, where $\consub(G)$ denotes the set of
conjugacy classes of subgroups of $G$. The zero element is represented by the empty set,
the unit is represented by $G/G$. The interesting feature of the
Burnside ring is its multiplicative structure.

Given a group homomorphism $f \colon G_0 \to G_1$ of finite
groups, restriction with $f$ defines a ring homomorphism $f^*
\colon A(G_1) \to A(G_0)$´. Thus $A(G)$ becomes a contravariant
functor from the category of finite groups to the category of
commutative rings. Induction defines a homomorphism of abelian
groups $f_* \colon A(G_0) \to A(G_1), \; [S] \mapsto [G_1 \times_f
S]$, which is not compatible with the multiplication. Thus $A(G)$
becomes a becomes a covariant functor from the category of finite
groups to the category of abelian groups.


\subsection{The Character Map and the Burnside Ring Congruences}
\label{subsec:The_Character_Map_and_the_Burnside_Ring_Congruences}

Let $G$ be a finite group.
Let $\consub(G)$ be \emph{the set of conjugacy classes $(H)$ of subgroups} $H \subseteq G$.
Define the character map
\begin{eqnarray}
\character^G \colon A(G) & \to & \prod_{(H) \in \consub(G)} \bbZ
\label{character_map_for_A(G)_for_finite_G}
\end{eqnarray}
by sending the class of a finite $G$-set $S$ to the numbers $\{|S^H|
\mid (H) \in \consub(G)\}$.  This is an injective ring homomorphism
whose image can be described by the so called
\emph{Burnside ring congruences} which we explain next.

In the sequel we denote for a subgroup $H \subseteq G$ by $N_GH$ its \emph{normalizer}
$\{g \in G \mid g^{-1}Hg = H\}$, by $C_GH = \{g \in G \mid gh=hg \text{ for } h  \in H\}$
its \emph{centralizer}, by $W_GH$ its \emph{Weyl group} $N_GH/H$ and by $[G:H]$ its \emph{index}.
Let $p_H \colon N_GH \to W_GH$ be the canonical projection. Denote for a cyclic group $C$ by $\Gen(C)$
the set of its generators.
We conclude from \cite[Proposition 1.3.5]{Dieck(1979)}
\begin{theorem}[Burnside ring congruences for finite groups]
\label{the:Burnside_ring_congruences_for_finite_groups}
An element $\{x(H)\} \in \prod_{(H) \in \consub(G)} \bbZ$ lies in the
image of the injective character map $\character^G$ defined in
\eqref{character_map_for_A(G)_for_finite_G} if and only if we have for every $(H) \in
\consub(G)$
$$\sum_{\substack{(C) \in \consub(W_GH)\\C \text{ cyclic}}}
|\Gen(C)| \cdot [W_GH : N_{W_GH}C] \cdot x(p_H^{-1}(C)) ~ \equiv ~
0 \quad \mod |W_GH|.$$
\end{theorem}

\begin{example}[$A(\bbZ/p$] \label{exa:A(Z/p)}
\em Let $p$ be a prime and let $G$ be the cyclic group $\bbZ/p$ of
order $p$. Then $A(G)$ is the free abelian group generated by
$[G]$ and $[G/G]$. The multiplication is determined by the fact
that $[G/G]$ is the unit and $[G] \cdot [G] = p \cdot [G]$.
There is exactly one non-trivial Burnside ring congruence, namely the
one for $H = \{1\}$ which in the notation of
Theorem~\ref{the:Burnside_ring_congruences_for_finite_groups} is
$$x(1) ~ \equiv ~ x(G) \mod p.$$
\em
\end{example}


\subsection{The Equivariant Euler Characteristic}
\label{subsec:The_Equivariant_Euler_Characteristic}

Next we recall the notion of a $G$-$CW$-complex.

\begin{definition}[$G$-$CW$-complex]
\label{def:G-CW-complex}
Let $G$ be a group. A \emph{$G$-$CW$-complex}
$X$ is a $G$-space together with a $G$-invariant filtration
$$\emptyset = X_{-1} \subseteq X_0 \subset X_1 \subseteq \ldots \subseteq
X_n \subseteq \ldots \subseteq\bigcup_{n \ge 0} X_n = X$$
such that $X$ carries the colimit topology
with respect to this filtration
(i.e., a set $C \subseteq X$ is closed if and only if $C \cap X_n$ is closed
in $X_n$ for all $n \ge 0$) and $X_n$ is obtained
from $X_{n-1}$ for each $n \ge 0$ by attaching
equivariant $n$-dimensional cells, i.e., there exists
a $G$-pushout
$$\comsquare
{\coprod_{i \in I_n} G/H_i \times S^{n-1}}
{\coprod_{i \in I_n} q_i^n}
{X_{n-1}}
{}
{}
{\coprod_{i \in I_n} G/H_i \times D^{n}}
{\coprod_{i \in I_n} Q_i^n}
{X_n}
$$
\end{definition}

A $G$-$CW$-complex $X$ is called \emph{finite} if it is built by
finitely many equivariant cells $G/H \times D^n$ and is called
\emph{cocompact} if $G\backslash X$ is compact. The conditions
finite and cocompact are equivalent for a $G$-$CW$-complex.
Provided that $G$ is finite, $X$ is compact if and only if $X$ is
cocompact A $G$-map $f \colon X \to Y$ of $G$-$CW$-complexes is
called \emph{cellular} if $f(X_n) \subseteq Y_n$ holds for all
$n$.

\begin{definition}[Equivariant Euler Characteristic]
\label{def:Equivariant_Euler_Characteristic}
Let $G$ be a finite group and $X$ be a finite $G$-$CW$-complex. Define its
\emph{equivariant Euler characteristic}
$$\chi^G(X) ~ \in A(G)$$
by
$$\chi^G(X) := \sum_{n = 0}^{\infty} (-1)^n \cdot \sum_{i \in I_n} [G/H_i]$$
after choices of the $G$-pushouts as in Definition~\ref{def:G-CW-complex}.
\end{definition}

This definition is independent of the choice of the $G$-pushouts by the next result.
The elementary proofs of the next two results are left to the reader. We denote by
$X^H$ and $X^{>H}$ respectively the subspace of $X$ consisting of elements $x \in X$ whose isotropy
group $G_x$ satisfies $H \subseteq G_x$ and $H\subsetneq G_x$ respectively.

\begin{lemma} \label{lem:basic_properties_of_chi^G(X)}
Let $G$ be a finite group.

\begin{enumerate}

\item \label{basic_properties_of_chi^G(X):homological_description}
Let $X$ be a finite $G$-$CW$-complex. Then
$$\chi^G(X) ~ = ~ \sum_{(H) \in \consub(G)} \chi\left(W_HG\backslash(X^H,X^{>H})\right) \cdot [G/H],$$
where $\chi$ denotes the classical (non-equivariant) Euler characteristic;

\item \label{basic_properties_of_chi^G(X):homotopy_invariance}
If $X$ and $Y$ are $G$-homotopy equivalent finite $G$-$CW$-complexes, then
$$\chi^G(X) = \chi^G(Y);$$

\item \label{basic_properties_of_chi^G(X):Additivity}
If
$$\comsquare{X_0}{i_0}{X_1}{i_2}{}{X_2}{}{X}$$
is a $G$-pushout of finite $G$-$CW$-complexes such that $i_1$ is an inclusion of finite
$G$-$CW$-complexes and $i_2$ is cellular, then
$$\chi^G(X) = \chi^G(X_1)  + \chi^G(X_1) - \chi^G(X_0);$$

\item \label{basic_properties_of_chi^G(X):product_formula}
If $X$ and $Y$ are finite $G$-$CW$-complexes, then $X \times Y$
with the diagonal $G$-action is a finite $G$-$CW$-complex
and
$$\chi^G(X \times Y) ~ = ~ \chi^G(X) \cdot \chi^G(Y);$$

\item \label{basic_properties_of_chi^G(X):character}
The image of $\chi^G(X)$ under the character
map $\character^G$ of \eqref{character_map_for_A(G)_for_finite_G}
is given by the collection of classical (non-equivariant)
Euler characteristics $\{\chi(X^H) \mid (H) \in \consub(G)\}$.

\end{enumerate}
\end{lemma}

An \emph{equivariant additive invariant} for finite $G$-$CW$-complexes is a pair $(A,a)$
consisting of an abelian group and an assignment $a$ which associates to every finite $G$-$CW$-complex
$X$ an element $a(X) \in A$ such that $a(\emptyset) = 0$, $G$-homotopy invariance and  additivity hold, i.e.,
the obvious versions of assertions~\ref{basic_properties_of_chi^G(X):homotopy_invariance}
and \ref{basic_properties_of_chi^G(X):Additivity} appearing in
Lemma~\ref{lem:basic_properties_of_chi^G(X)} are true.  An equivariant additive invariant $(U,u)$
is called \emph{universal} if for every  equivariant additive invariant $(A,a)$ there is precisely one
homomorphism of abelian groups $\phi \colon U \to A$ such that $\phi(u(X)) = a(X)$ holds for
every finite $G$-$CW$-complex $X$. Obviously $(U,u)$ is (up to unique isomorphism) unique if it exists.

\begin{theorem}[The universal equivariant additive invariant]
Let $G$ be a finite group.
The pair $(A(G),\chi^G)$ is the universal equivariant additive invariant for finite $G$-$CW$-complexes.
\end{theorem}


\subsection{The Equivariant Lefschetz Class}
\label{subsec:The_Equivariant_Lefschetz_Class}

The notion of an equivariant Euler characteristic can be
extended to the notion of an equivariant Lefschetz class
as follows.

\begin{definition} \label{def:equivariant_Lefschetz_class_Lambda^G(f)}
Let $G$ be a finite group and $X$ be a finite $G$-$CW$-complex. We
define the \emph{equivariant Lefschetz class} of a cellular
$G$-selfmap $f \colon X \to X$
\begin{eqnarray*}
\Lambda^G(f) & \in & A(G)
\end{eqnarray*}
by
$$\Lambda^G(f) ~ = ~ \sum_{(H) \in \consub(G)} \Lambda\left(W_GH\backslash (f^H,f^{>H})\right) \cdot [G/H],$$
where $\Lambda(W_GH\backslash (f^H,f^{>H})) \in \bbZ$ is the
classical  Lefschetz number of the endomorphism $W_GH\backslash
(f^H,f^{>H})$ of the pair of finite $CW$-complexes $W_GH\backslash
(X^H,X^{>H})$ induced by $f$.
\end{definition}

Obviously $\Lambda^G(\id \colon X \to X)$ agrees with $\chi^G(X)$.
The elementary proof of the next result is left to the reader.

\begin{lemma} \label{lem:basic_properties_of_Lambda^G(f)}
Let $G$ be a finite group.

\begin{enumerate}

\item \label{basic_properties_of_Lambda^G(f):homotopy invariance}
If $f$ and $g$ are $G$-homotopic $G$-selfmaps of a finite $GW$-$CW$-complex $X$, then
$$\Lambda^G(f) = \Lambda^G(g);$$

\item \label{basic_properties_of_Lambda^G(X):Additivity}
Let
$$\comsquare{X_0}{i_0}{X_1}{i_2}{}{X_2}{}{X}$$
is a $G$-pushout of finite $G$-$CW$-complexes such that $i_1$ is an inclusion of finite
$G$-$CW$-complexes and $i_2$ is cellular. Let $f_i \colon X_i \to X_i$
for $i = 0,1,2$ and $f \colon X \to X$ be the $G$-selfmaps
compatible with this $G$-pushout. Then
$$\Lambda^G(f) = \Lambda^G(f_1)  + \Lambda^G(f_2) - \Lambda^G(f_0);$$

\item \label{basic_properties_of_Lambda^G(f):product_formula}
Let $X$ and $Y$ be finite $G$-$CW$-complexes and $f \colon X \to X$ and $g \colon Y \to Y$ be $G$-selfmaps.
Then
$$\Lambda^G(f \times g) ~ = ~ \chi^G(X) \cdot \Lambda^G(g) + \chi^G(Y) \cdot \Lambda^G(f);$$

\item \label{basic_properties_of_Lambda^G(f):trace_property}
Let $f \colon X \to Y$ and $g \colon Y \to X$ be $G$-maps of finite $G$-$CW$-complexes. Then

$$\Lambda^G(f \circ  g) ~ = ~ \Lambda^G(g \circ  f);$$

\item \label{basic_properties_of_Lambda^G(X):character}
The image of $\Lambda^G(f)$ under the character map
$\character^G$ of~\eqref{character_map_for_A(G)_for_finite_G}
is given by the collection of classical (non-equivariant) Lefschetz numbers $\{\Lambda(f^H) \mid (H) \in \consub(G)\}$.

\end{enumerate}
\end{lemma}

One can also give a universal
property characterizing the equivariant Lefschetz class
(see~\cite{Laitinen-Lueck(1989)}).

The equivariant Lefschetz class has also the following homotopy theoretic meaning.

\begin{definition}
\label{def:homotopy_representation}
A \emph{$G$-homotopy representation} $X$ is a finite-dimensional $G$-$CW$-complex
such that for each subgroup $H \subseteq G$ the fixed point set $X^H$ is homotopy equivalent
to a sphere $S^{n(H)}$ for $n(H)$ the dimension of the $CW$-complex $X^H$.
\end{definition}

An example is the unit sphere $SV$ in an orthogonal representation $V$ of $G$. Denote by $[X,X]^G$ the set of $G$-homotopy classes
of $G$-maps $X \to X$. The proof of the next theorem can be found in \cite[Theorem~3.4 on page~139]{Lueck(1988)}
and is a consequence of the equivariant Hopf Theorem
(see for instance \cite[page 213]{Dieck(1979)}, \cite[II.4]{Dieck(1987)}, \cite{Laitinen(1986)}).

\begin{theorem} \label{the:[x,x]_cong_A(G)}
Let $X$ be a $G$-homotopy representation of the finite group $G$.
Suppose that

\begin{enumerate}

\item Every subgroup $H \subseteq G$ occurs as isotropy group of $X$;

\item $\dim(X^G) \ge 1$;

\item The group $G$ is nilpotent or for  every subgroup
$H \subseteq G$ we have $\dim(X^{>H}) + 2 \le \dim(X^H)$.

\end{enumerate}

Then the following map is an bijection of monoids, where the
monoid structure on the source comes from the composition and the
one on the target from the multiplication

$$\degree^G \colon [X,X]^G \xrightarrow{\cong} A(G),
\quad [f] ~ \mapsto (\Lambda^G(f) -1) \cdot (\chi^G(X) -1).$$
\end{theorem}

We mention that the image of the element $\degree^G(f)$ for a self-$G$-map of a $G$-homotopy representation
under the character map $\character^G$ of \eqref{character_map_for_A(G)_for_finite_G}
is given by the collection of (non-equivariant) degrees $\{\degree(f^H) \mid (H) \in \consub(G)\}$.


\subsection{The Burnside Ring and Equivariant Stable Cohomotopy}
\label{subsec:The_Burnside_Ring_and_Equivariant_Stable_Cohomotopy}

Let $X$ and $Y$ be two finite pointed $G$-$CW$-complexes. Pointed
means that we have specified an element in its $0$-skeleton which
is fixed under the $G$-action. If $V$ is a real
$G$-representation, let $S^V$ be its one-point compactification.
We will use the point at infinity as base point for $S^V$. If $V$
is an orthogonal representation, i.e., comes with an  $G$-invariant
scalar product, then $S^V$ is $G$-homeomorphic to the unit sphere
$S(V \oplus \bbR)$. Given two pointed $G$-$CW$-complexes $X$ and
$Y$ with base points $x$ and $y$, define their
\emph{one-point-union} $X \vee Y$ to be the pointed
$G$-$CW$-complex $X \times \{y\} \cup \{x\} \times Y \subseteq X
\times Y$ and their \emph{smash product} $X \wedge Y$ to be the
pointed $G$-$CW$-complex  $X \times Y/X \vee Y$.

We briefly introduce equivariant stable homotopy groups following the approach due to tom Dieck~\cite[II.6]{Dieck(1987)}.

If $V$ and $W$ are two complex $G$-representations, we write
$V \le W$ if there exists a complex $G$-representation $U$ and a linear $G$-isomorphism
$\phi \colon U \oplus V \to W$. If $\phi \colon U \oplus V \to W$ is a linear $G$-isomorphism,
define a map
$$b_{V,W} \colon [S^V \wedge X,S^V \wedge Y]^G \to [S^W \wedge X,S^W \wedge Y]^G$$
by the composition
\begin{multline*}
[S^V \wedge X,S^V \wedge Y]^G \xrightarrow{u_1}  [S^U \wedge S^V \wedge X,S^U \wedge S^V \wedge Y]^G
\\
\xrightarrow{u_2}  [S^{U \oplus V} \wedge X, S^{U \oplus V} \wedge Y]^G
\xrightarrow{u_3} [S^W \wedge X,S^W \wedge Y]^G,
\end{multline*}
where the map $u_1$ is given by $[f] \mapsto [\id_{S^U} \wedge f]$, the map $u_2$ comes from the
obvious $G$-homeomorphism $S^{U \oplus V} \xrightarrow{\cong} S^U \wedge S^V$ induced
by the inclusion $V \oplus W \to S^V \wedge S^W$ and the map
$u_3$ from the $G$-homeomorphism $S^{\phi} \colon S^{U \oplus V} \xrightarrow{\cong} S^W$.
Any two linear $G$-isomorphisms $\phi_0,\phi_1 \colon V_1 \to V_2$ between to complex
$G$-representations are isotopic as linear $G$-isomorphisms. (This is not true
for real $G$-representations.) This implies that the map
$b_{V,W}$ is indeed independent of the choice of $U$ and $\phi$. One easily checks
that $b_{V_2,V_1} \circ b_{V_0,V_1} = b_{V_0,V_2}$ holds for complex $G$-representations
$V_0$, $V_1$ and $V_2$ satisfying $V_0 \le V_1$ and $V_1 \le V_2$.

Let $I$ be the set of complex $G$-representations with underlying complex vector space $\bbC^n$ for some $n$.
(Notice that the collections of all complex $G$-representations does not form a set.)
Define on the disjoint union
$$\coprod_{V \in I} [S^V\wedge X,S^V \wedge Y]^G$$
an equivalence relation by calling $f \in [S^V \wedge X,S^V \wedge Y]^G$ and
$g \in [S^W\wedge X,S^W \wedge Y]^G$ equivalent if there exists a
representation $U \in I$ with $V \le U$ and $W \le U$ such that
$b_{V,U}([f]) = b_{W,U}([g])$ holds. Let $\omega_0^G(X,Y)$ for two pointed $G$-$CW$-complexes
$X$ and $Y$ be the set of equivalence classes.

If $V$ is any complex $G$-representation (not necessarily in I)
and $f \colon S^V \wedge X \to S^V \wedge Y$ is any
$G$-map, there exists an element $W \in I$
with $V \le W$ and we get an element in $\omega_n^G(X,Y)$  by $b_{V,W}([f])$.
This element is independent of the choice of $W$ and also denoted by
$[f] \in \omega_n^G(X,Y)$.

One can define the structure of an abelian group on the set $\omega_0^G(X,Y)$ as follows.
Consider elements $x,y \in \omega_0^G(X,Y)$. We can choose an element of the shape $\bbC \oplus U$
in $I$ for $\bbC$ equipped with the trivial $G$-action and $G$-maps
$f,g \colon S^{\bbC \oplus U} \wedge X \to S^{\bbC \oplus U} \wedge X$ representing
$x$ and $y$. Now using the standard pinching map $\nabla \colon S^{\bbC} \to S^{\bbC} \wedge S^{\bbC}$
one defines $x + y$ as the class of the $G$-map
\begin{multline*}
S^{\bbC \oplus U} \wedge X \xrightarrow{\cong} S^{\bbC}
\wedge S^U \wedge X \xrightarrow{\nabla \wedge \id_{S^U} \wedge \id_X}
(S^{\bbC} \vee S^{\bbC}) \wedge S^U \wedge X
\\
\xrightarrow{\cong}
(S^{\bbC} \wedge S^U \wedge X) \vee (S^{\bbC} \wedge S^U \wedge X)
\xrightarrow{\cong}
(S^{\bbC \oplus U} \wedge X) \vee (S^{\bbC\oplus U} \wedge X)
\xrightarrow{f \vee g}  S^{\bbC \oplus U} \wedge X .
\end{multline*}
The inverse of $x$ is defined by the class of
$$
S^{\bbC \oplus U} \wedge X \xrightarrow{\cong} S^{\bbC} \wedge S^U \wedge X \xrightarrow{d \wedge f}
S^{\bbC} \wedge S^U \wedge X
\xrightarrow{\cong} S^{\bbC \oplus U} \wedge X
$$
where $d \colon S^{\bbC} \to S^{\bbC}$ is any pointed map of degree $-1$.
This is indeed independent of the choices of $U$, $f$ and $g$.

We define the abelian groups
$$
\begin{array}{lclll}
\omega_n^G(X,Y)   & = & \omega_0^G(S^n \wedge X,Y) & &   n \ge 0;
\\
\omega_n^G(X,Y)   & = & \omega_0^G(X,S^{-n}Y)      & &   n \le 0;
\\
\omega_G^n(X,Y)   & = & \omega_{-n}^G(X,Y)         & &   n \in \bbZ;
\end{array}
$$
Obviously $\omega_n^G(X,Y)$ is functorial, namely contravariant in $X$ and covariant in $Y$.

Let $X$ and $Y$ be (unpointed) $G$-$CW$-complexes.
Let $X_+$ and $Y_+$ be the pointed $G$-$CW$-complexes
obtained from $X$ and $Y$ by adjoining a disjoint base point.
Denote by $\pt$ the one-point-space. Define abelian groups
$$
\begin{array}{lclll}
\pi_n^G(Y)        & = & \omega_n^G(\pt_+,Y_+)          & &   n \in \bbZ;
\\
\pi^n_G(X)        & = & \omega^n_G(X_+,\pt_+)          & &   n \in \bbZ;
\\
\pi^G_n           & = & \pi^G_n(\pt)                   & &   n \in \bbZ;
\\
\pi_G^n           & = & \pi^G_{-n}                     & &   n \in \bbZ.
\end{array}
$$
The abelian group $\pi^G_0 = \pi_G^0$ becomes a ring by the composition of
maps. The abelian groups $\pi_n^G(Y)$ define covariant functors in $Y$ and are called the
\emph{equivariant stable homotopy groups} of $Y$.
The abelian groups $\pi_G^n(X)$  define contravariant functors in $X$ and are called the
\emph{equivariant stable cohomotopy groups} of $X$.

We emphasize that our input
in $\pi_G^n$ and $\pi^G_n$ are unpointed $G$-$CW$-complexes. This is later consistent
with our constructions for infinite groups, where
all $G$-$CW$-spaces must be proper and therefore have empty $G$-fixed point sets and cannot have base points.

Theorem~\ref{the:[x,x]_cong_A(G)} implies the following result due to Segal~\cite{Segal(1971)}.

\begin{theorem} \label{the:pi^0_G_cong_A(G)}
The isomorphism $\deg^G$ appearing in Theorem~\ref{the:[x,x]_cong_A(G)}
induces an isomorphism of rings
$$\deg^G \colon \pi_G^0 \xrightarrow{\cong} A(G).$$
\end{theorem}

For a more sophisticated and detailed construction of and more information about the equivariant stable homotopy category
we refer for instance to
\cite{Greenlees-May(1995b)}, \cite{Lewis-May-Steinberger(1986)}.


\subsection{The Segal Conjecture for Finite Groups}
\label{subsec:The_Segal_Conjecture_for_Finite_Groups}

The equivariant stable cohomotopy groups $\pi^n_G(X)$ are modules over
the ring $\pi^0_G = A(G)$, the module structure is given by
composition of maps. The \emph{augmentation homomorphism}
$\epsilon^G \colon A(G) \to \bbZ$ is the ring homomorphism 
sending the class of a finite
set $S$ to $|S|$ which is just the component belonging to the trivial
subgroup of the character map defined in
\eqref{character_map_for_A(G)_for_finite_G}. The
\emph{augmentation ideal} $\bfI_G \subseteq A(G)$ is the kernel of
the augmentation homomorphism $\epsilon^G$.

For an (unpointed) $CW$-complex $X$
we denote by $\pi^n_s(X)$ the (non-equivariant) stable cohomotopy group of
$X_+$. This is in the previous notation for equivariant stable cohomotopy the same
as $\pi^n_{\{1\}}(X)$ for $\{1\}$ the trivial group. If $X$ is a finite $G$-$CW$-complex,
we can consider $\pi^n_s(EG \times_G X)$. Since $\pi^n_G(X)$ is a $A(G)$-module, we can also consider its
$\bfI_G$-adic completion denoted by $\pi^n_G(X)\widehat{_{\bfI_G}}$.
The following result is due to Carlsson~\cite{Carlsson(1984)}.

\begin{theorem}[Segal Conjecture for finite groups]
\label{the:Segal_Conjecture_for_finite_groups}
The \emph{Segal Conjecture for finite groups $G$} is true, i.e., for every finite group $G$ and
finite $G$-$CW$-complex $X$ there is an isomorphism
$$\pi^n_G(X)\widehat{_ {\bfI_G}} ~ \xrightarrow{\cong} ~ \pi^n_s(EG \times_G X).$$
\end{theorem}

In particular we get in the case  $X = \pt$ and $n = 0$ an isomorphism
\begin{eqnarray}
A(G)\widehat{_{\bfI_G}} & \xrightarrow{\cong} & \pi^0_s(BG).
\label{Segal_conjecture_for_finite_G_and_X=pt}
\end{eqnarray}
Thus the Burnside ring is linked via its $\bfI_G$-adic completion
to the stable cohomotopy of the classifying space $BG$
of a finite group $G$.

\begin{example}[Segal Conjecture for $\bbZ/p$] \label{exa:Segal_Conjecture_for_Z/p} \em
Let $G$ be the cyclic group $\bbZ/p$ of order $p$. We have computed $A(G)$ in
Example~\ref{exa:A(Z/p)}. If we put $x = [G] -  p \cdot [G/G]$, then the augmentation ideal
is generated by $x$. Since
$$x^2 ~ = ~ ([G] - p)^2 ~ = ~ [G]^2 - 2p \cdot [G] + p^2 = (-p) \cdot x,$$
we get $x^n = (-p)^{n-1} x$ and hence $\bfI_G^n = p^{n-1} \cdot \bfI_G$ for $n \in \bbZ$, $n \ge 1$.
This implies
$$A(G)\widehat{_{\bfI_G}} ~ = ~ \invlim{n \to \infty}{\bbZ \oplus \bfI_G/\bfI_G^n} ~ = ~ \bbZ \times \bbZ\widehat{_p},$$
where $\bbZ\widehat{_p}$ denotes the ring of $p$-adic integers.
\em
\end{example}


\subsection{The Burnside Ring as a Green Functor}
\label{subsec:The_Burnside_Ring_as_a_Green_Functor}

Let $R$ be an associative commutative ring with unit.
Let $\FGINJ$ be the category of finite
groups with injective group homomorphisms as morphisms.
Let $M\colon \FGINJ \to R\text{-}\MODULES$ be a bifunctor, i.e., a pair $(M_*,M^*)$
consisting of a covariant functor $M_*$ and a contravariant
functor $M^*$ from $\FGINJ$ to $R\text{-}\MODULES$ which agree on objects.
We will often denote for an injective  group homomorphism $f\colon H \to G$
the map $M_*(f)\colon M(H) \to M(G)$ by $\ind_f$ and
the map $M^*(f)\colon M(G) \to M(H)$ by $\res_f$ and write
$\ind_H^G = \ind_f$ and $\res_G^H = \res_f$ if $f$ is an inclusion of
groups. We call such a bifunctor $M$ a \emph{Mackey functor}
with values in $R$-modules if
\begin{enumerate}
\item
For an inner automorphism $c(g)\colon G \to G$ we have
$M_*(c(g)) = \id\colon  M(G) \to M(G)$;

\item
For an isomorphism of groups $f\colon G \xrightarrow{\cong} H$
the composites $\res_f \circ \ind_f$ and $\ind_f \circ \res_f$
are the identity;

\item Double coset formula\\[1mm]
We have for two subgroups $H,K \subset G$
$$\res_G^K \circ \ind_H^G ~ = ~\sum_{KgH \in K\backslash G/H}
\ind_{c(g)\colon H\cap g^{-1}Kg \to K}
 \circ \res_{H}^{ H\cap g^{-1}Kg},$$
where $c(g)$ is conjugation with $g$, i.e., $c(g)(h) = ghg^{-1}$.
\end{enumerate}

Let $\phi \colon R \to S$ be a homomorphism of associative commutative rings with unit.
Let $M$ be a Mackey functor with values in $R$-modules and let
$N$ and $P$  be Mackey functors with values in $S$-modules.
A \emph{pairing} with respect to $\phi$ is a family of  maps
$$m(H)\colon  M(H) \times N(H) \to P(H), \quad (x,y) \mapsto m(H)(x,y) =: x\cdot y,$$
where $H$ runs through the finite groups and we require the
following properties for all injective group homomorphisms
$f\colon  H \to K$ of finite groups:
$$\begin{array}{ll}
(x_1 + x_2) \cdot y = x_1 \cdot y + x_2 \cdot y
& \text{for }
x_1,x_2 \in M(H),y \in N(H);
\\
x \cdot (y_1 + y_2) = x\cdot y_1 + x \cdot y_2
& \text{for }
x \in M(H),y_1,y_2 \in N(H);
\\
(rx) \cdot y = \phi(r)(x\cdot y)
& \text{for }
r \in R, x \in M(H), y \in N(H);
\\
x \cdot sy = s (x\cdot y)
& \text{for }
s \in S, x \in M(H), y \in N(H);
\\
\res_f(x\cdot y) = \res_f(x)\cdot \res_f(y)
& \text{for }
x \in M(K), y \in N(K);
\\
\ind_f(x) \cdot y = \ind_f\left(x \cdot \res_f(y)\right)
& \text{for }
x \in M(H), y \in N(K);
\\
x \cdot \ind_f(y) = \ind_f\left(\res_f(x) \cdot y\right)
& \text{for }
x \in M(K), y \in N(H).
\end{array}$$

A \emph{Green functor} with values in $R$-modules is a Mackey functor
$U$ together with a pairing with respect to
$\id\colon  R \to R$ and elements $1_H \in U(H)$ for
each finite group $H$ such that for each finite group $H$ the pairing
$U(H) \times U(H) \to U(H)$ induces the structure of an $R$-algebra
on $U(H)$ with unit $1_H$ and for any morphism $f\colon H \to K$ in
$\FGINJ$ the map $U^*(f)\colon U(K) \to U(H)$
is a homomorphism of $R$-algebras with unit.
Let $U$ be a Green functor with values
in $R$-modules and $M$ be a Mackey functor
with values in $S$-modules. A (left) $U$-module
structure on $M$ with respect to the ring
homomorphism $\phi \colon  R \to S$ is a pairing
such that any of the maps $U(H) \times M(H) \to M(H)$ induces the structure
of a (left) module over the $R$-algebra $U(H)$ on the $R$-module $\phi^*M(H)$
which is obtained from the $S$-module $M(H)$ by $rx := \phi(r)x$ for
$r \in R$ and $x \in M(H)$.

\begin{theorem} \label{the:Burnside_ring_and_Green_functor}
\begin{enumerate}
\item \label{the:Burnside_ring_and_Green_functor:A_is_a_Green_functor}
The Burnside ring defines a Green functor with values in $\bbZ$-modules;

\item \label{the:Burnside_ring_and_Green_functor:A(G)-module structure}
If $M$ is a Mackey functor with values in $R$-modules, then $M$ is in a canonical way
a module over the Green functor given by the Burnside ring with respect to the canonical
ring homomorphism $\phi \colon \bbZ \to R$.
\end{enumerate}
\end{theorem}
\begin{proof}
\ref{the:Burnside_ring_and_Green_functor:A_is_a_Green_functor} Let
$f \colon H \to G$ be an injective homomorphism of groups. Define
$\ind_f \colon A(H) \to A(G)$ by sending the class of a finite
$H$-set $S$ to the class of the finite $G$-set $G \times_f S$.
Define $\res_f \colon A(G) \to A(H)$ by considering a finite
$G$-set as an $H$-set by restriction with $f$. One easily verifies
that the axioms of a Green functor with values in $\bbZ$-modules
are satisfied.
\\[1mm]
\ref{the:Burnside_ring_and_Green_functor:A(G)-module structure}
We have to specify for any finite group $G$ a pairing $m(G) \colon A(G) \times M(G) \to M(G)$.
This is done by the formula
$$m(G)\left(\sum_{i} n_i \cdot [G/H_i],x\right) ~ := ~
\sum_{i} n_i \cdot \ind_{H_i}^G \circ \res_G^{H_i}(x).$$
 One easily verifies that the axioms of a module over the Green functor
given by the Burnside ring are satisfied.
\end{proof}

Theorem~\ref{the:Burnside_ring_and_Green_functor} is the main
reason why the Burnside ring plays an important role in induction
theory. Induction theory addresses the question whether one can
compute the values of a Mackey functor on a finite group by its
values on a certain class of subgroups such as the family of
cyclic or hyperelementary groups. Typical examples of such Mackey
functors are the representation ring $R_F(G)$ or algebraic $K$ and
$L$-groups $K_n(RG)$ and $L_n(RG)$ of groups rings. The
applications require among other things a good understanding of
the prime ideals of the Burnside ring. For more information about
induction theory for finite groups we refer to the fundamental
papers by Dress~\cite {Dress(1973)},~\cite{Dress(1975)} and for
instance to~\cite[Chapter~6]{Dieck(1979)}.
Induction theory for infinite groups is developed in
\cite{Bartels-Lueck(2004b)}.

As an illustration we give an example how the Green-functor
mechanism works.

\begin{example}[Artin's Theorem] \label{exa:Artins_Theorem}
\em
Let $R_{\bbQ}(G)$ be the rational representation ring of the finite group $G$.
For any finite cyclic group $C$ one can construct an element
$$\theta_C \in R_{\bbQ}(C)$$
which is uniquely determined by the property that its character function
sends a generator of $C$ to $|C|$ and every other element of $C$ to zero.

Let $G$ be a finite group.
Let $\bbQ$ be the trivial $1$-dimensional rational $G$-representation.
It is not hard to check by a calculation with characters
that
\begin{eqnarray}
|G| \cdot[\bbQ] ~ = ~ \sum_{\substack{C \subset G\\C \cyclic}} \ind_C^H\theta_C.
\label{|G|cdot[Q]=sum_{Ccyclci}ind_C^Gtheta_C}
\end{eqnarray}

Assigning to a finite group $G$ the rational representation ring
$R_{\bbQ}(G)$ inherits the structure of a Green functor with
values in $\bbZ$-modules by induction and restriction. Suppose
that $M$ is a Mackey functor with values in $\bbZ$-modules which
is a module over the Green functor $R_{\bbQ}$. Then for every
finite group $G$  the cokernel of the map
$$\bigoplus_{\substack{C \subset G\\C \cyclic}} \ind_C^G \colon
\bigoplus_{\substack{C \subset G\\C \cyclic}} M(C) ~ \to ~ M(G)$$
is annihilated by multiplication with $|G|$. This follows from the
following calculation for $x \in M(G)$ based
on~\eqref{|G|cdot[Q]=sum_{Ccyclci}ind_C^Gtheta_C} and the axioms
of a Green functor and a module over it
\begin{multline*}
|G| \cdot x ~ = ~ (|G| \cdot [\bbQ]) \cdot x ~ = ~
\sum_{\substack{C \subset G\\C \cyclic}} \ind_C^H(\theta_C) \cdot x
~ = ~ \sum_{\substack{C \subset G\\C \cyclic}}  \ind_C^H(\theta_C \cdot \res_H^C x).
\end{multline*}
Examples for $M$ are algebraic $K$- and $L$-groups
$K_n(RG)$ and $L_n(RG)$ for any ring $R$ with $\bbQ \subseteq R$.
We may also take $M$ to be $R_F$ for any field $F$ of characteristic zero and then the statement
above is Artin's Theorem (see \cite[Theorem~26 on page 97]{Serre(1977)}.
\em
\end{example}


\subsection{The Burnside Ring and Rational Representations}
\label{subsec:The_Burnside_Ring_and_Rational_Representations}

Let $R_{\bbQ}(G)$ be the representation ring of finite-dimensional
rational $G$-re\-pre\-sen\-ta\-tion. Given a finite $G$-set $S$, let $\bbQ[S]$
be the rational $G$-representation given by the $\bbQ$-vector space
with the set $S$ as basis. The next result is due to Segal~\cite{Segal(1972)}.

\begin{theorem}{\bf (The Burnside ring and the rational representation ring for finite groups).}
\label{the:The_Burnside_ring_and_the_rational_representation_ring_for_finite_groups}
Let $G$ be a finite group. We obtain a ring homomorphism
$$P^G \colon A(G) \to R_{\bbQ}(G), \quad [S] \mapsto [\bbQ[S]].$$
It is rationally surjective. If $G$ is a $p$-group for some prime $p$,
it is surjective. It is bijective if and only  if $G$ is cyclic.
\end{theorem}


\subsection{The Burnside Ring and Homotopy Representations}
\label{subsec:The_Burnside_Ring_and_Homotopy_Representations}

We have introduced the notion of a $G$-homotopy representation in
Definition~\ref{def:homotopy_representation}. The join of two
$G$-homotopy representations is again a $G$-homotopy
representation. We call two $G$-homotopy representations $X$ and
$Y$ stably $G$-homotopy equivalent if for some $G$-homotopy
representation $Z$ the joins $X \ast Z$ and $Y \ast Z$ are
$G$-homotopy equivalent. The stable $G$-homotopy classes of
$G$-homotopy representations together with the join define an
abelian semi-group. The \emph{$G$-homotopy representation group}
$V(G)$ is the associated Grothendieck group. It may be viewed as
the homotopy version of the representation ring. Taking the unit
sphere yields a group homomorphism $R_{\bbR}(G) \to V(G)$.

The  \emph{dimension function} of a $G$-homotopy representation $X$
$$\dim(X) \in  \prod_{(H)} \bbZ$$
associates to the conjugacy class $(H)$ of a subgroup $H  \subseteq G$ the dimension of
$X^H$. The question which elements in $\prod_{(H)} \bbZ$ occur as $\dim(X)$ is studied for instance in
\cite{Dieck(1982)}, \cite{Dieck(1986)}, \cite[III.5]{Dieck(1987)} and \cite{Dieck-Petrie(1982)}.
Define $V(G,\dim)$ by the exact sequence
$$0 \to V(G,\dim) \to V(G)  \xrightarrow{\dim}  \prod_{(H)} \bbZ.$$
Let $\Pic(A(G))$ be the Picard group of the Burnside ring, i.e., 
the abelian group of projective $A(G)$-module of rank one
with respect to the tensor product.
The next result is taken from~\cite[6.5]{Dieck-Petrie(1982)}.

\begin{theorem}[$V(G,\dim)$ and the Picard group of $A(G)$]
There is an isomorphism
$$V(G,\dim)  \xrightarrow{\cong} \Pic(A(G)).$$
\end{theorem}

Further references about the Burnside ring of finite groups are
\cite{Bouc(2000)},  \cite{Dress(1969)},\cite{Dress-Siebeneicher-Yoshida(1991)},  \cite{Dress-Siebeneicher-Yoshida(1992)},
 \cite{Laitinen(1979)}, \cite {Laitinen(1981)}, \cite{Morimoto(2002)}, \cite{Dieck(1975)},
\cite{Yoshida(1990a)}.


\typeout{------------Section 2: The Finite-$G$-Set-Version of the Burnside Ring -------------}

\section{The Finite-$G$-Set-Version of the Burnside Ring}
\label{sec:The_Finite_G-Set_Version_of_the_Burnside_Ring}

From now on $G$ can be any (discrete) group and needs not to be finite
anymore.  Next we give a first definition of the Burnside ring for infinite
groups.

\begin{definition}{\bf (The finite-$G$-set-version of the Burnside ring).}
\label{def:The_finite-$G$-set-version_of_the_Burnside_ring}
The isomorphism classes of finite $G$-sets form a commutative associative semi-ring with unit
under the disjoint union and the cartesian product. The  \emph{finite-$G$-set-version of the Burnside ring} $\finsetA(G)$
is the  Grothendieck ring associated to this semi-ring.
\end{definition}

To avoid any confusion, we emphasize that finite $G$-set means a
finite set with a $G$-action.  This definition is word by word the same as
given for a finite group in Definition~\ref{def:A(G)_for_finite_G}.

Given a group homomorphism $f \colon G_0 \to G_1$ of groups, restriction with $f$ defines a ring homomorphism
$f^* \colon A(G_1) \to A(G_0)$´. Thus $A(G)$ becomes a contravariant functor from the category
of groups to the category of commutative rings. Provided that the image of $f$ has finite index,
induction defines a homomorphism of abelian groups
$f_* \colon A(G_0) \to A(G_1), \; [S] \mapsto [G_1 \times_f S]$, which is not compatible
with the multiplication.


\subsection{Character Theory and Burnside Congruences for the Finite-$G$-Set-Version}
\label{subsec:Character_Theory_and_Burnside_Congruences_for_the_Finite-G-Set-Version}

The definition of the character map~\eqref{character_map_for_A(G)_for_finite_G}
makes also sense for infinite groups and we denote it by

\begin{eqnarray}
\finsetcharacter^G \colon \finsetA(G) & \to & \prod_{(H) \in \consub(G)} \bbZ,
\quad [S] ~ \mapsto [|(S^H|)_{(H)}.
\label{character_map_for_finsetA(G)}
\end{eqnarray}

Given a group homomorphism $f \colon G_0 \to G_1$, define a ring homomorphism
\begin{eqnarray}
f^* \colon \prod_{(H_1) \in \consub(G_1)} & \to & \prod_{(H_0) \in \consub(G_0)} \bbZ
\label{f^* for prod_{consub(?)}}
\end{eqnarray}
by sending $\{x(H_1)\}$ to $\{x(f(H_1))\}$. One easily checks
\begin{lemma} \label{lem:char_is_natural}
  The following diagram of commutative rings with unit commutes for
  every group homomorphism $f \colon G_0 \to G_1$
  $$
  \comsquare{\overline{A}(G_1)}{\overline{\character}^{G_1}}{\prod_{(H_1)
      \in \consub(G_1)} \bbZ} {f^*}{f^*}
  {\overline{A}(G_0)}{\overline{\character}^{G_0}}{\prod_{(H_0) \in
      \consub(G_0)} \bbZ}
  $$
\end{lemma}

\begin{theorem}[Burnside ring congruences for $\finsetA(G)$]
\label{the:Burnside_ring_congruences_for_finsetA(G)}
The character map $\finsetcharacter^G$ is an injective ring homomorphism.

Already the composition
$$\finsetA(G) ~ \xrightarrow{\finsetcharacter^G} ~ \prod_{(H) \in \consub(G)} \bbZ
~ \xrightarrow{\pr} ~ \prod_{\substack{(H) \in \consub(G)\\ [G:H] < \infty}} \bbZ$$
for $\pr $ the obvious projection is injective.

An element $x = \{x(H)\} \in \prod_{(H) \in \consub(G)} \bbZ$ lies in
the image of the character map $\finsetcharacter^G$ defined in
\eqref{character_map_for_finsetA(G)} if and only if it satisfies
the following two conditions:

\begin{enumerate}

\item \label{the:Burnside_ring_congruences_for_finsetA(G):conditionK_x} There
  exists a normal subgroup $K_x \subseteq G$ of finite index such that
  $x(H) = x (H \cdot K_x)$ holds for all $H \subseteq G$, where
  $H \cdot K_x$ is the subgroup $\{hk \mid h \in H, k \in K_x\}$;

\item \label{the:Burnside_ring_congruences_for_finsetA(G):condition_congruences} We
  have for every $(H) \in \consub(G)$ with $[G :H ] < \infty$:
  $$\sum_{\substack{(C) \in \consub(W_GH)\\C \text{ cyclic}}}
  |\Gen(C)| \cdot [W_GH : N_{W_GH}C] \cdot x(p_H^{-1}(C)) ~ \equiv ~ 0
  \quad \mod |W_GH|,$$
  where $p_H \colon N_GH \to W_GH$ is the
  obvious projection.
\end{enumerate}
\end{theorem}
\begin{proof}
  Obviously $\finsetcharacter^G$ is a ring homomorphism.

  Suppose that
  $x \in \finsetA(G)$ lies in the kernel of $\finsetcharacter^G$. For any finite $G$-set the intersection of all its isotropy groups is
  a normal subgroup of finite index in $G$.  Hence can find
  an epimorphism $p_x \colon G \to Q_x$ onto a finite group $Q_x$ and $\overline{x}\ \in A(Q_x)$ such that
  $x$ lies in the image of $p_x^* \colon \overline{A}(Q_x) \to \overline{A}(G)$. Since the map
  $$p_x^* \colon  \prod_{(H) \in \consub(Q_x)} \bbZ ~ \to ~  \prod_{(K) \in \consub(G)} \bbZ
  $$
  is obviously injective and the character map $\finsetcharacter^{Q_x}$ is injective by
  Theorem~\ref{the:Burnside_ring_congruences_for_finite_groups}, we conclude $x = 0$.
  Hence $\finsetcharacter^G$ is injective.

  Suppose that $y$ lies in the image of $\overline{\character}^G$.
  Choose $x \in \finsetA(G)$ with $\overline{\character}^G(x) = y$.
  As explained above we can find
  an epimorphism $p_x \colon G \to Q_x$ onto a finite group $Q_x$ and $\overline{x}\ \in A(Q_x)$ such that
  $p_x^* \colon \overline{A}(Q_x) \to \overline{A}(G)$ maps $\overline{x}$ to $x$.
  Then Condition~\ref{the:Burnside_ring_congruences_for_finsetA(G):conditionK_x} is
  satisfied by Lemma~\ref{lem:char_is_natural} if we take $K_x$ to be the kernel of $p_x$.
  Condition~\ref{the:Burnside_ring_congruences_for_finsetA(G):condition_congruences}
  holds for $x$ since the proof of
  Theorem~\ref{the:Burnside_ring_congruences_for_finite_groups} carries
  though word by word to the case, where $G$ is possibly infinite but $H \subseteq
  G$ is required to have finite index in $G$ and hence $W_GH$ is
  finite.

  We conclude that $\finsetcharacter^G(x) = 0$ if and only if $\pr \circ \finsetcharacter^G(x) = 0$ holds.
  Hence $\pr \circ \finsetcharacter^G$ is injective.

  Now suppose that $x = \{y(H)\} \in \prod_{(H) \in \consub(G)} \bbZ$
  satisfies Condition~\ref{the:Burnside_ring_congruences_for_finsetA(G):conditionK_x} and
  Condition~\ref{the:Burnside_ring_congruences_for_finsetA(G):condition_congruences}.
  Let $Q_x= G/K_x$ and let $p_x\colon G \to Q_x$ be the
  projection. In the sequel we abbreviate $Q = Q_x$ and $p = p_x$.
  Then Condition~\ref{the:Burnside_ring_congruences_for_finsetA(G):conditionK_x}
  ensures that $x$ lies in the image of the injective
  map
  $$p^* \colon \prod_{(H) \in \consub(Q)} \bbZ \to \prod_{(K) \in
    \consub(G)} \bbZ.$$
  Let $y \in \prod_{(H) \in \consub(Q)} \bbZ$ be
  such a preimage.  Because of Lemma~\ref{lem:char_is_natural} it suffices to prove that
  $y$ lies in the image of the character map
  $$\character^{Q} \colon A(Q) \to \prod_{(H) \in
    \consub(Q)} \bbZ.$$
  By Theorem~\ref{the:Burnside_ring_congruences_for_finite_groups}
  this is true if and only if for every subgroup
  $K \subseteq Q$ the congruence
  $$\sum_{\substack{(C) \in \consub(W_QK)\\C \text{ cyclic}}}
  |\Gen(C)| \cdot [W_QK : N_{W_QK}C] \cdot
  y\left(p_K^Q)^{-1}(C)\right) ~ \equiv ~ 0 \quad \mod |W_QK|$$
  holds, where $p^Q_K \colon N_QK \to W_GK$ is the projection. Fix a
  subgroup $K \subseteq Q$. Put $H = p^{-1}(K) \subseteq G$. The
  epimorphism $p \colon G \to Q$ induces an isomorphism $\overline{p}
  \colon W_GH \xrightarrow{\cong} W_QK$.
  Condition~\ref{the:Burnside_ring_congruences_for_finsetA(G):condition_congruences}
  applied to $x$ and $H$ yields
  $$\sum_{\substack{(C) \in \consub(W_GH)\\C \text{ cyclic}}}
  |\Gen(C)| \cdot [W_GH : N_{W_GH}C] \cdot x(p_H^{-1}(C)) ~ \equiv ~ 0
  \quad \mod |W_GH|.$$
  For any cyclic subgroup $C \subseteq W_GH$ we
  obtain a cyclic subgroup $\overline{p}(C) \subseteq W_QK$ and we
  have
\begin{eqnarray*}
|\Gen(C)| & = & |\Gen(\overline{p}(C))|;
\\
{[W_GH : N_{W_GH}C]} & = &  {[W_QK : N_{W_QK}\overline{p}(C)]};
\\
x\left(p_H^{-1}(C)\right) & = & y\left((p_K^Q)^{-1}(\overline{p}(C))\right).
\end{eqnarray*}
Now the desired congruence for y follows. This finishes the proof of
Theorem~\ref{the:Burnside_ring_congruences_for_finsetA(G)}.
\end{proof}

\begin{example}[$\finsetA$ of the integers] \label{exa:finsetA(Z)} \em
  Consider the infinite cyclic group $\bbZ$. Any subgroup of
  finite index is of the form $n\bbZ$ for some $n \in \bbZ$, $n \ge
  1$. As an abelian group $\finsetA(Z)$ is generated by the
  classes $[\bbZ/n\bbZ]$ for $n \in \bbZ$, $n \ge 1$.  The
  condition~\ref{the:Burnside_ring_congruences_for_finsetA(G):condition_congruences}
  appearing in Theorem~\ref{the:Burnside_ring_congruences_for_finsetA(G)}
  reduces to the condition that for every subgroup
  $n\bbZ$ for $n \in \bbZ$, $n \ge 1$ the congruence
  $$\sum_{m \in \bbZ, m\ge 1, m|n} \phi\left(\frac{n}{m}\right) \cdot x(m)
  ~ \equiv ~ 0 \quad \mod n$$
  holds, where $\phi$ is the
  Euler function, whose value $\phi(k)$ is $|\Gen(\bbZ/k\bbZ)|$.  The
  condition~\ref{the:Burnside_ring_congruences_for_finsetA(G):conditionK_x}
  reduces to the condition that there exists $n_x \in \bbZ$, $n_x \ge
  1$ such that for all $m \in \bbZ$, $m \ge 1$ we have $x(m\bbZ) =
  x(gcd(m,n_x)\bbZ)$, where $gcd(m,n_x)$ is the greatest common
  divisor of $m$ and $n_x$.  \em
\end{example}

\begin{remark}[The completion $\widehatA(G)$ of $\overline{A}(G)$]
\label{rem:The_completion_widehat(A)_of_overlineA(G)}
\em
We call a $G$-set \emph{almost finite} if each isotropy group has finite index
and for every positive integer $n$ the number of orbits
$G/H$ in $S$ with $[G:H] \le n$ is finite. A $G$-set $S$ is almost finite if and only if for every subgroup $H \subseteq G$ of finite index
the $H$-fixed point set $S^H$ is finite and $S$ is the union $\bigcup_{\substack{(H) \in \consub(G)\\ [G:H] < \infty}} S^H$.
Of course every finite $G$-set $S$ is almost finite.
The disjoint union and the cartesian product with the diagonal $G$-action
of two almost finite $G$-sets is again almost finite.
Define $\widehatA(G)$ as the Grothendieck ring of the semi-ring of almost finite $G$-sets
under the disjoint union and the cartesian product. There is an obvious inclusion of rings
$\finsetA(G) \to \widehatA(G)$. We can define as before a character map
\begin{eqnarray}
\widehatcharacter^G \colon \widehatA(G) & \to & \prod_{\substack{(H) \in \consub(G)\\ [G:H] < \infty}} \bbZ,
\quad [S] ~ \mapsto (|(S^H|)_{(H)}.
\label{character_map_for_widehatA(G)}
\end{eqnarray}
We leave it to the reader to check that $\widehatcharacter^G$ is injective,
and that an element $x$ in $\prod_{\substack{(H) \in \consub(G)\\ [G:H] < \infty}} \bbZ$ lies in its image
if and only if $x$ satisfies
condition~\ref{the:Burnside_ring_congruences_for_finsetA(G):condition_congruences}
appearing in Theorem~\ref{the:Burnside_ring_congruences_for_finsetA(G)}.

Dress and Siebeneicher~\cite{Dress-Siebeneicher(1988)} analyze
$\widehat{A}(G)$ for profinite groups $G$ and put it into relation
with the Witt vector construction. They also explain that
$\widehat{A}(G)$ is a completion of $\finsetA(G)$. The ring
$\widehat{A}(\bbZ)$ is studied and put in relation to the necklace
algebra, $\lambda$-rings and the universal ring of Witt vectors in
\cite{Dress-Siebeneicher(1989)}. \em
\end{remark}


\subsection{The Finite-$G$-Set-Version and the Equivariant Euler Characteristic and the Equivariant Lefschetz Class}
\label{subsec:The_Finite_G_Set_Version_and_Equivariant_Euler_Characteristic_and_the_Equivariant_Lefschetz_Class}

The results of
Sections~\ref{subsec:The_Equivariant_Euler_Characteristic}
and~\ref{subsec:The_Equivariant_Lefschetz_Class} carry over to
$\finsetA(G)$ if one considers only finite $G$-$CW$-complexes $X$
whose isotropy group all have finite index in $G$. But this is not
really new since for any such $G$-$CW$-complex $X$ there is a
subgroup $H \subseteq G$, namely the intersection of all isotropy
groups,  such that $H$ is normal, has finite index in $G$ and acts
trivially on $X$. Thus  $X$ is a finite $Q$-$CW$-complex for the
finite group $Q = G/H$ and all these invariant are obtained from
the one over $Q$ by the applying the obvious ring homomorphism
$A(Q) = \finsetA(Q) \to \finsetA(G)$ to the invariants already
defined over the finite group $Q$.


\subsection{The Finite-$G$-Set-Version as a Green Functor}
\label{subsec:The_Finite_G_Set_Version_as_a_Green_Functor}

The notions and results of
Subsection~\ref{subsec:The_Burnside_Ring_as_a_Green_Functor} carry
over to the finite-$G$-set-version $\finsetA(G)$ for an infinite
group $G$, we replace the category $\FGINJ$ by the category
$\GRIFI$ whose objects are groups and whose morphisms are
injective group homomorphisms whose image has finite index in the
target. However, for infinite groups this does not seem to be the
right approach to induction theory. The approach presented in
Bartels-L\"uck~\cite{Bartels-Lueck(2004b)} is more useful. It is based on
classifying spaces for families and aims at reducing the family of
subgroups, for instance from all finite subgroups to all
hyperelementary finite subgroups or from all virtually cyclic
subgroups to the family of subgroups which admit an epimorphism to
a hyperelementary group and whose kernel is trivial or infinite
cyclic.


\subsection{The Finite-$G$-Set-Version and the Swan Ring}
\label{subsec:The_Finite_G_Set_Version_and_the_Swan_Ring}

Let $R$ be a commutative ring. Let $\Sw^f(G;R)$ be the abelian group which is generated
by the $RG$-isomorphisms classes of $RG$-modules which are finitely generated free over $R$
with the relations $[M_0] - [M_1]-[M_2] = 0$ for any short exact $RG$-sequence
$0 \to M_0 \to M_1 \to M_2 \to 0$ of such $RG$-modules.
It becomes a commutative ring, the so called \emph{Swan ring} $\Sw^f(G;R)$, by the tensor product
$\otimes_R$. If $G$ is finite and $F$ is a field, then $\Sw^f(G;F)$ is the same as the representation ring
$R_F(G)$ of (finite-dimensional) $G$-representations over $K$.

Let $G_0(RG)$ be the abelian group which is generated
by the $RG$-isomorphism classes of finitely generated $RG$-modules
with the relations $[M_0] - [M_1]-[M_2] = 0$ for any short exact $RG$-sequence
$0 \to M_0 \to M_1 \to M_2 \to 0$ of such $RG$-modules. There is an obvious map
$$\phi \colon \Sw^f(G;R) \to G_0(RG)$$ of abelian groups. It is an isomorphism
if $G$ is finite and $R$ is a principle ideal domain. This follows from
\cite[Theorem~38.42 on page 22]{Curtis-Reiner(1987)}.

We obtain a ring homomorphism
\begin{eqnarray}
\finsetP^G \colon \finsetA(G) ~ \to ~ \Sw^f(G;\bbQ), \quad [S] \mapsto [R[S]],
\label{finsetA(G)_to_Sw(G;R)}
\end{eqnarray}
where $R[S]$ is the finitely generated free $R$-module with the finite set $S$ as
basis and becomes a $RG$-module by the $G$-action on $S$.

Theorem~\ref{the:The_Burnside_ring_and_the_rational_representation_ring_for_finite_groups}
does not carry over to $\finsetA(G)$ for infinite groups. For instance, the determinant
induces a surjective homomorphism
$$\det \colon \Sw(\bbZ;\bbQ) \to \bbQ^*, \quad [V] ~ \mapsto~ \det(l_t \colon V \to V),$$
where $l_t$ is left multiplication with a fixed generator $t \in \bbZ$. Given a finite $\bbZ$-set $S$, the map
$l_t \colon \bbQ [S] \to \bbQ [S]$ satisfies $(l_t)^n = \id$ for some $n \ge 1$ and
hence $\finsetP^{\bbZ}(\bbQ [S]) = \pm -1$. Hence the image of the composition
$\det \circ \finsetP_{\bbZ}$ is contained in $\{\pm 1\}$. Therefore the  map
$\finsetP^{\bbZ}$ of~\eqref{finsetA(G)_to_Sw(G;R)} is not rationally surjective.


\subsection{Maximal Residually Finite Quotients}
\label{subsec:Maximal_Residually_Finite_Quotients}

  Let $G$ be a group. Denote by $G_0$ the intersection of all normal
  subgroups of finite index.  This is a normal subgroup. Let
  $p \colon  G \to G/G_0$ be the projection. Recall that $G$ is called \emph{residually finite}
  if for every element $g \in G$ with $g \not= 1$ there exists a homomorphism onto a finite group
  which sends $g$ to an element different from $1$. If $G$ is countable, then $G$ is residually finite if and only if
  $G_0$ is trivial. The projection $p \colon  G \to G_{\mrf} := G/G_0$
  is the projection onto the maximal residually finite quotient
  of $G$, i.e., $G_{\mrf}$ is residually finite and every epimorphism $f \colon G \to Q$
  onto a residually finite group $Q$ factorizes
  through $p$ into a composition $G \xrightarrow{p} G_{\mrf} \xrightarrow{\overline{f}} Q$.
  If $G$ is a finitely generated subgroup of $GL_n(F)$ for some field $F$, then
  $G$ is residually finite (see \cite{Mal'cev(1940)},
  \cite[Theorem 4.2]{Wehrfritz(1973)}). Hence for every finitely generated group $G$ each
  $G$-representation $V$ with coefficients in a field $F$ is obtained by restriction
  with $p \colon  G \to G_{\mrf}$ from a $G_{\mrf}$-representation. In particular
  every $G$-representation with coefficient in a field $F$ is trivial if $G$ is finitely generated and
  $G_{\mrf}$ is trivial.

  One easily checks that
  $$p^* \colon \finsetA(G_{\mrf}) \xrightarrow{\cong} \finsetA(G)$$
  is an isomorphism. In particular we have $\finsetA(G) = \bbZ$ if $G_{\mrf}$ is trivial.
  If $G$ is finitely generated, then
  $$p^* \colon  \Sw^f(G_{\mrf};F) \xrightarrow{\cong} \Sw^f(G;F)$$
  is an isomorphism. In particular we have $ \Sw^f(G;F) = \bbZ$ if $G$ is finitely generated and
  $G_{\mrf}$ is trivial.

\begin{example}[$\finsetA(\bbZ/p^{\infty})$ and $\Sw^f(\bbZ/p^{\infty};\bbQ)$]
\label{exa:finsetA(Z/p^{infty} and Sw^f(Z/p^{infty};Q)} \em
  Let $\bbZ/p^{\infty}$ be the \emph{Pr\"ufer group}, i.e., the colimit
  of the directed system of injections of abelian groups $\bbZ/p \to
  \bbZ/p^2 \to \bbZ/p^3 \to \cdots$.  It can  be identified with both
  $\bbQ/\bbZ_{(p)}$ and $\bbZ[1/p]/\bbZ$.  We want to show that the following diagram
  is commutative and consists of isomorphisms
  $$
  \begin{CD}\finsetA(\{1\}) = \bbZ @>p^* > \cong > \finsetA(\bbZ/p^{\infty})
  \\
  @V \finsetP^{\{1\}} V \cong V @V \finsetP^{\bbZ/p^{\infty}} V \cong V
  \\
  \Sw^f(\{1\};\bbQ) = \bbZ @> p^* > \cong > \Sw^f(\bbZ/p^{\infty};\bbQ)
  \end{CD}
  $$
  where $p \colon \bbZ/p^{\infty} \to \{1\}$ is the projection.
  Obviously the diagram commutes and the left vertical arrow is bijective.
  Hence it remains to show that the horizontal arrows are bijective.

  Let $f \colon \bbZ/p^{\infty} \to Q$ be any
  epimorphism onto a finite group.  Since $\bbZ/p^{\infty}$ is
  abelian, $Q$ is a finite abelian group. Since any element in
  $\bbZ/p^n$ has $p$-power order, we conclude from the definition of
  $\bbZ/p^{\infty}$ as a colimit that $Q$ is a finite abelian
  $p$-group. Since $\bbQ$ is $p$-divisible, the quotient $Q$ must be
  $p$-divisible. Therefore $Q$ must be trivial. Hence $(\bbZ/p^{\infty})_{\mrf}$ is trivial
  and the upper horizontal arrow is bijective.

  In order to show that the lower horizontal arrow is bijective, it suffices to show
  that every (finite-dimensional) rational $\bbZ/p^{\infty}$-representation  $V$ is trivial.
  It is enough to show that for every subgroup $\bbZ/p^m$ its restriction $\res_{i_m} V$
  for the inclusion $i_m \colon \bbZ/p^m \to \bbZ/p^{\infty} $ is trivial.
  For this purpose choose a positive integer $n$ such that $\dim_{\bbQ}(V) < (p-1) \cdot p^n$.
  Consider the rational $\bbZ/p^{m+n}$-representation $\res_{i_{m+n}} V$.
  Let $p_k^{m+n} \colon \bbZ/p^{m+n} \to \bbZ/p^k$ be the canonical projection.
  Let $\bbQ(p^k)$ be the rational $\bbZ/p^k$-representation given by adjoining a primitive $p^k$-th
  root of unity to $\bbQ$. Then the dimension of $\bbQ(p^k)$ is $(p-1)\cdot p^{k-1}$. A complete system of representatives
  for the isomorphism classes of irreducible rational
  $\bbZ/p^{m+n}$-representations is $\{\res_{p^{m+n}_k} \bbQ(p^k) \mid k = 0,1,2, \ldots m+n\}$.
  Since $\dim_{\bbQ}(V) < (p-1)\cdot p^n$, there exists a rational $\bbZ/p^n$-representation with
  $W$ with $\res_{i_{m+n}} V \cong \res_{p^{m+n}_n} W$. Hence we get an isomorphism of rational $\bbZ/p^m$-representations
  $$\res_{i_{m}} V \cong \res_{i_{m,m+n}} \circ \res_{p^{m+n}_n} W$$
  where $i_{m,m+n} \colon \bbZ/p^m \to \bbZ/p^{m+n}$ is the inclusion. Since the composition
  $p^{m+n}_n \circ i_{m,m+n}$ is trivial, the  rational $\bbZ/p^m$-representation
  $\res_{i_{m}} V$ is trivial.
  \em
  \end{example}

  It is not true that $$\Sw^f(\bbZ/p^{\infty};\bbC) \to  \Sw^f(\{1\};\bbC) = \bbZ$$
  is bijective because $\Sw^f(\bbZ/p^{\infty};\bbC)$ has as abelian group infinite rank
  (see Example~\ref{R_{inv,C}(Z/p^{infty})}).


\typeout{------------ Section 3: The inverse-limit-Version of the Burnside Ring-----------}

\section{The Inverse-Limit-Version of the Burnside Ring}
\label{sec:The_Inverse_Limit_Version_of_the_Burnside_Ring}

In this section we present the  inverse-limit-definition of the Burnside ring for infinite groups.

The \emph{orbit category} $\OrGF{G}{}$ has as objects homogeneous spaces
$G/H$ and as morphisms $G$-maps.  Let $\SubGF{G}{}$ be the category whose
objects are subgroups $H$ of $G$. For two subgroups $H$ and $K$ of $G$
denote by $\conhom_G(H,K)$ the set of group homomorphisms $f\colon H
\to K$, for which there exists an element $g \in G$ with $gHg^{-1}
\subset K$ such that $f$ is given by conjugation with $g$, i.e.,  $f =
c(g)\colon H \to K, ~ h \mapsto ghg^{-1}$.  Notice that $c(g) = c(g')$
holds for two elements $g,g' \in G$ with $gHg^{-1} \subset K$ and
$g'H(g')^{-1} \subset K$ if and only if $g^{-1}g'$ lies in the
centralizer $C_GH = \{g \in G \mid gh=hg \mbox{ for all } h \in H\}$
of $H$ in $G$. The group of inner automorphisms of $K$ acts on
$\conhom_G(H,K)$ from the left by composition. Define the set of
morphisms
$$\mor_{\SubGF{G}{}}(H,K) ~ := ~ \inn(K)\backslash \conhom_G(H,K).$$

There is a natural projection $\pr\colon \OrGF{G}{} \to \SubGF{G}{}$ which
sends a homogeneous space $G/H$ to $H$.  Given a $G$-map $f\colon G/H
\to G/K$, we can choose an element $g \in G$ with $gHg^{-1} \subset K$
and $f(g'H) = g'g^{-1}K$. Then $\pr(f)$ is represented by $c(g)\colon
H \to K$. Notice that $\mor_{\SubGF{G}{}}(H,K)$ can be identified with the
quotient $\mor_{\OrGF{G}{}}(G/H,G/K)/C_GH$, where $g \in C_GH$ acts on
$\mor_{\OrGF{G}{}}(G/H,G/K)$ by composition with $R_{g^{-1}}\colon G/H \to
G/H, \hspace{3mm} g'H \mapsto g'g^{-1}H$.  We mention as illustration
that for abelian $G$ the set of morphisms $\mor_{\SubGF{G}{}}(H,K)$ is empty if $H$ is not a
subgroup of $K$, and consists of precisely one element given by the
inclusion $H \to K$ if $H$ is a subgroup in $K$.

Denote by $\OrGF{G}{\calfin} \subset \OrGF{G}{}$ and $\SubGF{G}{\calfin} \subset
\SubGF{G}{}$ the full subcategories, whose objects $G/H$ and $H$ are given
by finite subgroups $H \subset G$.

\begin{definition}{\bf (The inverse-limit-version of the Burnside ring).}
\label{def:The_inverse-limit-version_of_the_Burnside_ring}
The \emph{inverse-limit-version of the Burnside ring} $\invA(G)$ is defined to be the
commutative ring with unit given by the inverse limit of the
contravariant functor
$$A(?) \colon \SubGF{G}{\calfin} \to \RINGS, \quad H ~ \mapsto ~ A(H).$$
\end{definition}

Since inner automorphisms induce the identity on $A(H)$, the
contravariant functor appearing in the definition above is
well-defined.

Consider a group homomorphism $f \colon G_0 \to G_1$. We obtain a
covariant functor $\SubGF{f}{\calfin} \colon \SubGF{G_0}{\calfin} \to
\SubGF{G_1}{\calfin}$ sending an object $H$ to $f(H)$.  A morphism $u
\colon H \to K$ given by $c(g) \colon H \to K$ for some $g \in G$ with
$gHg^{-1} \subseteq K$ is sent to the morphism given by $c(f(g))
\colon f(H) \to f(K)$. There is an obvious transformation from the
composite of the functor $A(?) \colon
\SubGF {G_1}{\calfin} \to \RINGS$ with $\SubGF{f}{\calfin}$ to the
functor $\invA(?) \colon \SubGF{G_0}{\calfin} \to \RINGS$.  It
is given for an object $H \in \SubGF{G_0}{\calfin}$ by the ring
homomorphism $A(f(H)) \to A(H)$ induced by the group homomorphism
$f|_H \colon H \to f(H)$. Thus we obtain a ring homomorphism $\invA(f)
\colon \invA(G_1) \to \invA(G_0)$.  So $\invA$ becomes a contravariant functor
$\GROUPS\to \RINGS$.

Definition~\ref{def:The_inverse-limit-version_of_the_Burnside_ring}
reduces to the one for finite groups presented in
Subsection~\ref{sec:Review_of_the_Burnside_Ring_for_Finite_Groups}
since for a finite group $G$ the object $G \in \SubGF{G}{\calfin}$ is a
terminal object.

There is an obvious ring homomorphism, natural in $G$,
\begin{eqnarray}
T^G \colon \finsetA(G) & \to & \invA(G)
\label{T^G:finsetA(G)_to_invA(G)}
\end{eqnarray}
which is induced from the various ring homomorphisms
$\finsetA(i_H) \colon \finsetA(G) \to \finsetA(H) = A(H)$
for the inclusions $i_H \colon H \to G$ for each finite subgroup $H
\subseteq G$.  The following examples show that it is neither
injective nor surjective in general.


\subsection{Some Computations of the Inverse-Limit-Version}
\label{subsec:Some_Computations_of_the_Inverse-Limit-Version}

\begin{example}[$\invA(G)$ for torsionfree $G$]
\label{exa:invA(G)_for_torsionfree_G} \em Suppose that $G$ is
torsionfree. Then $\SubGF{G}{\calfin}$ is the trivial category
with precisely one object and one morphisms. This implies that the
projection $\pr \colon G \to \{1\}$ induces a ring isomorphism
$$
\invA(\pr) \colon \invA(\{1\}) = \bbZ ~ \xrightarrow{\cong} ~ \invA(G).$$
In particular we conclude from Example~\ref{exa:finsetA(Z)} that the
canonical ring homomorphism
$$T^{\bbZ} \colon \finsetA(\bbZ) \to \invA(\bbZ)$$
of \eqref{T^G:finsetA(G)_to_invA(G)} is not injective.  \em
\end{example}

\begin{example}[Groups with appropriate maximal finite subgroups] \label{exa:conditions_M_and_NM} \em
  Let $G$ be a discrete group which is not torsionfree.  Consider the following assertions
  concerning $G$:
\begin{itemize}

\item[(M)] Every non-trivial finite subgroup of $G$ is contained in a
  unique maximal finite subgroup;

\item[ (NM)] If $M \subseteq G$ is maximal finite, then $N_GM = M$.

\end{itemize}

The conditions (M) and (NM) imply the following: Let $H$ be a
non-trivial finite subgroup of $G$. Then there is a unique maximal
finite subgroup $M_H$ with $H \subseteq M_H$ and the set of
morphisms in $\SubGF{G}{\calfin}$ from $H$ to $M_H$ consists of
precisely one element which is represented by the inclusion $H \to
M_H$. Let $\{M_i \mid i \in I\}$ be a complete set of
representatives of the conjugacy classes of maximal finite
subgroups of $G$. Denote by $j_i \colon M_i \to G$, $k_i \colon
\{1\} \to M_i$ and $k \colon \{1\} \to G$ the inclusions. Then we
obtain a short exact sequence
\begin{multline*}
0 \to \invA(G) \xrightarrow{\invA(j_{\{1\}}) \times \prod_{i \in I} \invA(j_i)}
\invA(\{1\}) \times \prod_{i \in I} \invA(M_i)
\\
\xrightarrow{\Delta - \prod_{i
    \in I} \invA(k_i)} \prod_{i \in I} \invA(\{1\} \to 0,
\end{multline*}
where $\Delta
\colon \invA(\{1\} \to \prod_{i \in I} \invA(\{1\}$ is the diagonal embedding.
If we define $\widetilde{\invA}(G)$ as the kernel of $\invA(G) \to \invA(\{1\})$,
this gives an isomorphism
$$\widetilde{\invA}(G) \xrightarrow{\prod_{i \in I} \widetilde{\invA}(j_i)}
\prod_{i \in I} \widetilde{\invA}(M_i).$$

Here are some examples of groups $Q$ which satisfy conditions (M) and
(NM):
\begin{itemize}

\item Extensions $1 \to \bbZ^n \to G \to F \to 1$ for finite $F$ such
  that the conjugation
  action of $F$ on $\bbZ^n$ is free outside $0 \in \bbZ^n$. \\[1mm]
  The conditions (M) and (NM) are satisfied by \cite[Lemma
  6.3]{Lueck-Stamm(2000)}.

\item Fuchsian groups $F$ \\[1mm]
  See for instance \cite[Lemma 4.5]{Lueck-Stamm(2000)}).  In
  \cite{Lueck-Stamm(2000)} the larger class of cocompact planar groups
  (sometimes also called cocompact NEC-groups) is treated.

\item Finitely generated one-relator groups $G$\\[1mm]
  Let $G = \langle (q_i)_{i \in I} \mid r \rangle$ be a presentation
  with one relation.  Let $F$ be the free group with basis $\{q_i \mid
  i \in I\}$. Then $r$ is an element in $F$. There exists an element
  $s \in F$ and an integer $m \ge 1$ such that $r = s^m$, the cyclic
  subgroup $C$ generated by the class $\overline{s} \in Q$ represented
  by $s$ has order $m$, any finite subgroup of $G$ is subconjugated to
  $C$ and for any $q \in Q$ the implication $q^{-1}Cq \cap C \not= 1
  \Rightarrow q \in C$ holds.  These claims follow from
  \cite[Propositions 5.17, 5.18 and 5.19 in II.5 on pages 107 and
  108]{Lyndon-Schupp(1977)}.  Hence $Q$ satisfies (M) and (NM) and the
  inclusion $i \colon C \to G$ induces an isomorphism
  $$\invA(i) \colon \invA(G) \xrightarrow{\cong} \invA(C).$$
\end{itemize}
\em
\end{example}

\begin{example}[Olshanskii's group] \em
  There is for any prime number $p > 10^{75}$ an infinite finitely
  generated group $G$ all of whose proper subgroups are finite of
  order $p$ \cite{Olshanskii(1982c)}. Obviously $G$ contains no
  subgroup of finite index.  Hence the inclusion $i \colon \{1\} \to
  G$ induces an isomorphism
  $$\finsetA(i) \colon \finsetA(G) \to A(\{1\}) = \bbZ$$
  (see Subsection~\ref{subsec:Maximal_Residually_Finite_Quotients}). Let $H$ be a finite
  non-trivial subgroup of $G$. Then $H$ is isomorphic to $\bbZ/p$ and
  agrees with its normalizer.  So the conditions appearing in Example~\ref{exa:conditions_M_and_NM}
  are satisfied. Hence we obtain
  an isomorphism
  $$\widetilde{\invA}(G) \xrightarrow{\cong} \prod_{\substack{(H) \in \consub_f(G)\\H \not= \{1\}}}
  \widetilde{\invA}(\bbZ/p),$$
  where $\consub_f(G)$ is the set of conjugacy classes of finite subgroups. This implies
  that the natural map
  $$T^{G} \colon \finsetA(G) \xrightarrow{\cong} \invA(G)$$
  of
  \eqref{T^G:finsetA(G)_to_invA(G)} is not surjective.  \em
\end{example}

\begin{example}[Extensions of $\bbZ^n$ with $\bbZ/p$ as quotient] \label{exa:extensions_by_Z/p} \em
  Suppose that $G$ can be written as an extension $1 \to A \to G \to
  \bbZ/p \to 1$ for some fixed prime number $p$ and for $A = \bbZ^n$
  for some integer $n \ge 0$ and that $G$ is not torsionfree.  The
  conjugation action of $G$ on the normal subgroup $A$ yields the
  structure of a $\bbZ[\bbZ/p]$-module on $A$. Every non-trivial
  element $g \in G$ of finite order $G$ has order $p$ and satisfies
  $$N_G\langle g \rangle = C_G\langle g \rangle = A^{\bbZ/p} \times
  \langle g \rangle.$$
  In particular the conjugation action of
  $N_G\langle g \rangle$ on $\langle g \rangle$ is trivial.  There is
  a bijection
  $$\mu \colon H^1(\bbZ/p;A) ~ \xrightarrow{\cong} ~ \consub_f(G),$$
  where $H^1(\bbZ/p;A)$ is the first cohomology of $\bbZ/p$ with
  coefficients in the $\bbZ[\bbZ/p]$-module $A$. If we fix an element
  $g \in G$ of order $p$ and a generator $s \in \bbZ/p$, the bijection
  $\mu$ sends $[u] \in H^1(\bbZ/p;A)$ to $(\langle ug \rangle)$ of the
  cyclic group $\langle ug\rangle$ of order $p$ if $[u] \in
  H^1(\bbZ/p;A)$ is represented by the element $u$ in the kernel of
  the second differential $A \to A, ~ a \mapsto \sum_{i=0}^{p-1} s^i
  \cdot a$.  Hence we obtain an exact sequence
  $$0 \to \invA(G) \to \invA(\{1\}) \times \prod_{H^1(\bbZ/p;A)} \invA(\bbZ/p) \to
  \prod_{H^1(\bbZ/p;A)} \invA(\{1\}) \to 0$$
  This gives
  an isomorphism
  $$\widetilde{\invA}(G) \xrightarrow{\cong} \prod_{H^1(\bbZ/p;A)}
  \widetilde{A}(\bbZ/p).$$
  \em
\end{example}


\subsection{Character Theory and Burnside Congruences for the Inverse-Limit-Version}
\label{subsec:Character_Theory_and_Burnside_Congruences_for_the_Inverse-Limit-Version}

Next we define a character map for infinite groups $G$ and determine
its image generalizing Theorem~\ref{the:Burnside_ring_congruences_for_finite_groups}.

Let $\consub_f(G)$ be the set of conjugacy classes $(H)$ of finite
subgroups $H \subseteq G$.  Given a group homomorphism $f \colon G_0
\to G_1$, let $\consub_f(G_0) \to \consub_f(G_1)$ be the map sending
the $G_0$-conjugacy class of a finite subgroup $H\subseteq G_0$ to the
$G_1$-conjugacy class of $f(H) \subseteq G_1$.  We obtain a covariant
functor
$$\consub \colon \SubGF{G}{\calfin} \to \SETS, \quad H ~ \mapsto
\consub(H).$$
For each finite subgroup $H \subseteq G$ the inclusion
$H \to G$ induces a map $\consub(H) \to \consub_f(G)$ sending $(K)$ to
$(K)$. These fit together to a bijection of sets
\begin{eqnarray}
\phi^G \colon \colim{H \in \SubGF{G}{\calfin}}{\consub(H)}  & \xrightarrow{\cong}& \consub_f(G).
\label{phi^g colimt consub(H) to consub_f(G)}
\end{eqnarray}
One easily checks that $\phi^G$ is well-defined and surjective.
Next we show injectivity.
Consider two elements $x_0$ and $x_1$ in the source of $\phi^G$ with
$\phi^G(x_0) = \phi^G(x_1)$. For $i = 0,1$ we can choose an object
$H_i \in \SubGF{G}{\calfin}$ and an element $(K_i) \in \consub_f(H_i)$
such that the structure map of the colimit for the object $(H_i)$
sends $(K_i)$ to $x_i$. Then $\phi^G(x_0) = \phi^G(x_1)$ means that
the subgroups $K_0$ and $K_1$ of $G$ are conjugated in $G$. Hence we
can find $g \in G$ with $gK_0g^{-1} = K_1$. The morphism $K_0 \to H_0$
induced by the inclusion yields a map $\consub(K_0) \to \consub(H_0)$
sending $(K_0)$ to $(K_0)$.  The morphism $K_0 \to H_1$ induced by the
conjugation homomorphism $c(g) \colon K_0 \to H_1$ yields a map
$\consub(K_0) \to \consub(H_1)$ sending $(K_0)$ to $(K_1)$. This
implies $x_0 = x_1$. 

By the universal property of the colimit we
obtain an isomorphism of abelian groups
\begin{multline}
  \psi^G \colon \map\left(\colim{H \in \SubGF{G}{\calfin}}{\consub(H)},\bbZ\right)
  \\
  \xrightarrow{\cong} \invlim{(H) \in \SubGF {G}{\calfin}}{\map(\consub(H)};\bbZ).
\label{psi^G(map(colim,Z)) to invlim map(-,Z)}
\end{multline}

Define the character map
\begin{eqnarray}
 \invcharacter^G \colon \invA(G) & \to & \prod_{(H) \in \consub_f(G)} \bbZ
\label{character^G_invA(G)_to_prod_{(H),|H|<infty}Z}
\end{eqnarray}
to be the map for which the composition with the isomorphism
\begin{multline*}
  \prod_{(H) \in \consub_f(G)} \bbZ ~ = ~ \map(\consub_f(G),\bbZ)
  \xrightarrow{\map(\phi^G,\id)}
  \map\left(\colim{H \in \SubGF{G}{\calfin}}{\consub(H)},\bbZ\right)
  \\
  \xrightarrow{\psi^G} \invlim{H \in \SubGF{G}{\calfin}}{\map(\consub(H),\bbZ)}
\end{multline*}
is the map
\begin{multline*}
  \invA(G) ~ = ~ \invlim{H \in \SubGF{G}{\calfin}}{\map(\consub(H),\bbZ)}
  \xrightarrow{\invlim{H \in \SubGF{G}{\calfin}}{\character^H}}
  \\
  \invlim{H \in \SubGF{G}{\calfin}}{\map(\consub(H);\bbZ)},
\end{multline*}
where $\character^H \colon A(H) \to \map(\consub(H),\bbZ)$ is the map
defined in \eqref{character_map_for_A(G)_for_finite_G}.

\begin{theorem}[Burnside ring congruences for $\invA(G)$]
\label{the:Burnside_ring_congruences_for_invA(G)}
Let $x$ be an element in $\prod_{(H) \in \consub_f(G)} \bbZ$. Then:
\begin{enumerate}

\item \label{the:Burnside_ring_congruences_for_invA(G):injectivity}

  The character map
  $$\invcharacter^G\colon \invA(G) \to \prod_{(H) \in \consub_f(G)} \bbZ$$
  of
  \eqref{character^G_invA(G)_to_prod_{(H),|H|<infty}Z} is injective;

\item \label{the:Burnside_ring_congruences_for_invA(G):congruences} The
  element $x$ lies in the image of the character map
  $$\invcharacter^G\colon \invA(G) \to \prod_{(H) \in \consub_f(G)} \bbZ$$
  of
  \eqref{character^G_invA(G)_to_prod_{(H),|H|<infty}Z} if and only if
  for every finite subgroup $K \subseteq G$ the following condition
  $C(K)$ is satisfied: The image of $x$ under the map induced by the
  inclusion $i_K \colon K \to G$
  $$(i_K)^* \prod_{(H) \in \consub_f(G)} \bbZ ~ =
  ~\map(\consub(H),\bbZ) ~ \to ~ \prod_{(L) \in \consub_f(K)} \bbZ ~ =
  ~ \map(\consub(K),\bbZ)$$
  satisfies the Burnside ring congruences
  for the finite group $K$ appearing in
  Theorem~\ref{the:Burnside_ring_congruences_for_finite_groups};

\item \label{the:Burnside_ring_congruences_for_invA(G):redundance} If
  $K_0 \subseteq K_1 \subseteq G$ are two subgroups, then condition
  $C(K_1)$ implies condition $C(K_0)$.
\end{enumerate}
\end{theorem}
\begin{proof}
  This follows from Theorem~\ref{the:Burnside_ring_congruences_for_finite_groups} and the fact that the inverse limit is left exact.
\end{proof}

\begin{example}{\bf (Finitely many conjugacy classes of finite subgroups).}
\label{exa:rank_of_invA(G)_for_|ccs_f(G)|_finite} \em
Suppose that $G$ has only finitely many conjugacy classes of finite subgroups.
Then we conclude from Theorem~\ref{the:Burnside_ring_congruences_for_invA(G)} that the cokernel of
the injective character map $\invcharacter^G\colon \invA(G) \to \prod_{(H) \in \consub_f(G)} \bbZ$
is finite. Hence $\invA(G)$ is a finitely generated free abelian group of rank $|\consub_f(G)|$.
\em
\end{example}

\begin{example}[$\invA(\bbZ/p^{\infty})$] \label{exa:invA(Z/p^{infty}} \em
  We have introduced in Example~\ref{exa:finsetA(Z/p^{infty} and Sw^f(Z/p^{infty};Q)} the
  Pr\"ufer group $\bbZ/p^{\infty}$ as $\colim{n \to \infty}{\bbZ/p^n}$.
  Each $\bbZ/p^n$ represents a finite subgroup and each
  finite subgroup arises in this way. Hence $\consub(\bbZ/p^{\infty})$
  is on one-to-one-correspondence with $\bbZ^{\ge 0} = \{n \in \bbZ
  \mid n \ge 0\}$. Thus $x \in \prod_{(H) \in
    \consub(\bbZ/p^{\infty})} \bbZ$ can be written as a sequence
  $\{x(n)\} = \{x(n) \mid n \in \bbZ^{\ge 0}\}$, where $x(n)$
  corresponds to the value of $x$ at $\bbZ/p^n$.

  Consider the finite subgroup $\bbZ/p^m$. Its subgroups are given by
  $\bbZ/p^k$ for $k = 0,1,2 \ldots m$.  Then condition $C(\bbZ/p^m)$
  reduces to the set of congruences for each $k = 0,1,2, \ldots, m-1$
  $$\sum_{C \subset (\bbZ/p^m)/(\bbZ/p^k)} \Gen(C) \cdot
  x(p_k^{-1}(C)) ~ \equiv ~ 0 \quad \mod p^{m-k},$$
  where $p_k \colon
  \bbZ/p^m \to (\bbZ/p^m)/(\bbZ/p^k)$ is the projection.  More
  explicitly, the condition $C(\bbZ/p^m)$ reduces to the set of
  congruences for each $k = 0,1,2, \ldots, m-1$
  $$x(k) + \sum_{i = 1}^{m-k} p^{i-1} \cdot (p-1) \cdot x(k+i) ~
  \equiv ~ 0 \quad \mod p^{m-k},$$
  which can be rewritten as
  $$
  \sum_{i = 0}^{m-k-1} p^i \cdot (x(k+i) - x(k+i+1)) ~ \equiv ~ 0
  \quad \mod p^{m-k}.$$
  One can see that $C(\bbZ/p^{m_1})$ implies
  $C(\bbZ/p^{m_0})$ for $m_0 \le m_1$ as predicted by
  Theorem~\ref{the:Burnside_ring_congruences_for_invA(G)}
  \ref{the:Burnside_ring_congruences_for_invA(G):redundance}.

  Suppose that $x$ satisfies $C(\bbZ/p^m)$ for $m = 0,1,2, \ldots $.
  We want to show inductively for $l = 0,1,2 \ldots $ that $x(j)
  \equiv x(j+1) \mod p^l$ holds for $j = 0,1,2, \ldots$.  The
  induction begin $l = 0$ is trivial, the induction step from $l-1$ to
  $l \ge 1$ done as follows.  The $m$-th equation appearing in
  condition $C(l+m)$ yields
  $$\sum_{i = 0}^{l-1} p^i \cdot (x(m+i) - x(m+i+1)) ~ \equiv ~ 0
  \quad \mod p^{l}.$$
  Since by induction hypothesis $x(k+i) - x(k+i+1)
  \equiv 0 \mod p^{l-1}$ holds, this reduces to
  $$x(m) - x(m+1) ~ \equiv ~ 0 \quad \mod p^{l}.$$
  This finishes the
  induction step.

  Since $x(j) \equiv x(j+1) \mod p^l$ holds for $l = 0,1,2, \ldots$,
  we conclude $x(j) = x(j+1)$ for $j = 0,1,2, \ldots $. On the other
  hand, if $x(j) = x(j+1)$ holds for $j = 0,1,2, \ldots $, then $x$
  obviously satisfies the conditions $C(\bbZ/p^m)$ for $m = 0,1,2,
  \ldots$. Theorem~\ref{the:Burnside_ring_congruences_for_invA(G)}
  \ref{the:Burnside_ring_congruences_for_invA(G):injectivity} and
  \ref{the:Burnside_ring_congruences_for_invA(G):congruences} shows
  that the character map
  $$\invcharacter^{\bbZ/p^{\infty}} \colon \invA(\bbZ/p^{\infty}) ~ \to ~
  \prod_{(H) \in \consub_f(\bbZ/p^{\infty})} \bbZ$$
  is injective and
  its image consists of the copy of the integers given by the constant
  series. This implies that the projection $\pr \colon \bbZ/p^{\infty}
  \to \{1\}$ induces a ring isomorphism
  $$
  \invA(\pr) \colon \invA(\{1\}) = \bbZ ~ \xrightarrow{\cong} ~
  \invA(\bbZ/p^{\infty}).$$
  In particular we conclude from
  Example~\ref{exa:finsetA(Z/p^{infty} and Sw^f(Z/p^{infty};Q)} that the canonical ring
  homomorphism
  $$T^{\bbZ/p^{\infty}} \colon \finsetA(\bbZ/p^{\infty})
  \xrightarrow{\cong} \invA(\bbZ/p^{\infty})$$
  of \eqref{T^G:finsetA(G)_to_invA(G)} is bijective.  \em
\end{example}


\subsection{The-Inverse-Limit-Version of the Burnside Ring and Rational Representations}
\label{subsec:The_Inverse_Limit_Version_of_the_Burnside_Ring_and_Rational_Representations}

Analogously to $\invA(G)$ one defines
$R_{\inv,F}(G)$ for a field $F$ to be the
commutative ring with unit given by the inverse limit of the
contravariant functor
$$R_{\inv,F}(?) \colon \SubGF{G}{\calfin} \to \RINGS, \quad H ~ \mapsto ~ R_F(H).$$
This functor has been studied for $F = \bbC$ for instance
in \cite{Adem(1992a)}, \cite{Adem(1993b)}.
The system of maps $P^H \colon A(H) \to R_{\bbQ}(H)$ for the finite subgroups $H \subseteq G$ appearing
in Theorem~\ref{the:The_Burnside_ring_and_the_rational_representation_ring_for_finite_groups}
defines a ring homomorphism
\begin{eqnarray}
\invP^G \colon \invA(G) & \to & R_{\inv,\bbQ}(G).
\label{inv^P^G}
\end{eqnarray}
The system of the restriction maps for every finite subgroup $H \subseteq G$ induces a
homomorphism
\begin{eqnarray}
S^{G,F} \colon \Sw^f(G;F) & \to & R_{\inv,F}(G).
\label{S^{G,F}}
\end{eqnarray}
Although each of the maps $P^H$ for the finite subgroups $H \subseteq G$ are rationally surjective by
Theorem~\ref{the:The_Burnside_ring_and_the_rational_representation_ring_for_finite_groups},
the map $\invP$ need not to be rationally surjective in general, since
inverse limits do not respects surjectivity or rationally surjectivity in general.

\begin{example}[$R_{\inv;\bbQ}(\bbZ/^{\infty})$]
\label{R_{inv,Q}(Z/p^{infty})} \em
Since every finite subgroup of $\bbZ/p^{\infty}$ is cyclic, we conclude from
Theorem~\ref{the:The_Burnside_ring_and_the_rational_representation_ring_for_finite_groups}
that the map
$$\invP^{\bbZ/p^{\infty}} \colon \invA(\bbZ/p^{\infty}) \xrightarrow{\cong} R_{\inv,\bbQ}(\bbZ/p^{\infty})$$
is bijective. We have already seen in Example~\ref{exa:invA(Z/p^{infty}}
that $p^* \colon \invA(\{1\}) \to \invA(\bbZ/p^{\infty})$ is bijective.
We conclude from Example~\ref{exa:finsetA(Z/p^{infty} and Sw^f(Z/p^{infty};Q)}
that the following diagram  is commutative and consists of isomorphisms
$$
\begin{CD}
\finsetA(\bbZ/p^{\infty}) @> T^{\bbZ/p^{\infty}} > \cong > \invA(\bbZ/p^{\infty})
\\
@V \finsetP^{\bbZ/p^{\infty}} V \cong V @V\invP^{\bbZ/p^{\infty}} V \cong V
\\
\Sw^f(\bbZ/p^{\infty};\bbQ)  @> S^{\bbZ/p^{\infty};\bbQ}_{\inv} > \cong >
R_{\inv,\bbQ}(\bbZ/p^{\infty})
\end{CD}
$$
and is isomorphic by the maps induced by the projection $p \colon \bbZ/p^{\infty} \to \{1\}$
to the following commutative diagram whose corners are all isomorphic to $\bbZ$ and whose arrows are all the identity
under this identification.
$$
\begin{CD}
\finsetA(\{1\}) @> T^{\{1\}} > \cong > \invA(\{1\})
\\
@V \finsetP^{\{1\}} V \cong V @V\invP^{\{1\}} V \cong V
\\
\Sw^f(\{1\};\bbQ)  @> S^{\{1\};\bbQ}_{\inv} > \cong >
R_{\inv,\bbQ}(\{1\})
\end{CD}
$$
\em
\end{example}

\begin{example}[$\Sw^f(\bbZ/p^{\infty};\bbC)$ and $R_{\inv;\bbC}(\bbZ/^{\infty})$]
\label{R_{inv,C}(Z/p^{infty})} \em
On the other hand let us consider $\bbC$ as coefficients.
Consider the canonical map
$$S^{\bbZ/p^{\infty},\bbC} \colon \Sw^f(\bbZ/p^{\infty};\bbC) \to R_{\inv,\bbC}(\bbZ/p^{\infty})$$
which is induced by the restriction maps for all inclusions $H \to G$ of finite subgroups.
If $\phi_H \colon R_{\inv,\bbC}(\bbZ/p^{\infty}) \to R_{\bbC}(H)$ is the structure map of the inverse limit defining
$R_{\inv,\bbC}(\bbZ/p^{\infty})$ for the finite subgroup $H \subseteq \bbZ/p^{\infty}$, then the composition
$$\Sw^f(\bbZ/p^{\infty};\bbC) \xrightarrow{S^{\bbZ/p^{\infty},\bbC} } R_{\inv,\bbC}(\bbZ/p^{\infty}) \xrightarrow{\psi_H} R_{\bbC}(H)$$
is the map given by restriction with the inclusion of the finite subgroup $H$ in $\bbZ/p^{\infty}$.
We claim that this composition is surjective.  Choose $n$ with $H = \bbZ/p^n$.
We have to find for every $1$-dimensional complex $\bbZ/p^n$-representation
$V$ a $1$-dimensional complex $\bbZ/p^{\infty}$-representation $W$ such that $V$ is the restriction
of $V$. If $V$ is given by the homomorphism $\bbZ/p^n \to S^1, \overline{k} \mapsto \exp(2\pi i k/p^n)$,
then the desired $W$ is given by the homomorphism
$$\bbZ/p^{\infty} = \bbZ[1/p]/\bbZ \to S^1, \quad k \mapsto \exp(2\pi i k).$$
This implies that both $\Sw^f(\bbZ/p^{\infty};\bbC)$ and $R_{\inv,\bbC}(\bbZ/p^{\infty})$
have infinite rank as abelian groups.
\em
\end{example}


\typeout{---------------------- Section 4: The Covariant Burnside Group --------------------------------}

\section{The Covariant Burnside Group}
\label{sec:The_Covariant_Burnside_Group}

Next we give a third version for infinite groups which however will only be an abelian group,
not necessarily a ring.

\begin{definition}[Covariant Burnside group] \label{def:covariant_Burnside_group}
  Define the \emph{covariant Burnside group} $\underlineA(G)$ of a
  group $G$ to be the Grothendieck group which is associated to the
  abelian monoid under disjoint union of $G$-isomorphism classes of
  proper cofinite $G$-sets $S$, i.e., $G$-sets $S$ for which the
  isotropy group of each element in $S$ and the quotient $G\backslash
  S$ are finite.
\end{definition}

The cartesian product of two proper cofinite $G$-sets with the
diagonal action is proper but not cofinite unless $G$ is finite. So
for infinite group $G$ we do not get a ring structure on the
Burnside group $\underlineA(G)$.  If $G$ is finite the underlying
abelian group of the Burnside ring $A(G)$ is just $\underlineA(G)$.
Given a group homomorphism $f \colon G_0 \to G_1$, induction yields a
homomorphism of abelian group $\underlineA(G_0) \to
\underlineA(G_1)$ sending $[S]$ to $[G_1 \times_f S]$. Thus
$\underlineA$ becomes a covariant functor from $\GROUPS$ to $\bbZ-\MODULES$.

In the sequel we denote by $R[S]$ for a commutative ring $R$ and a set
$S$ the free $R$-module with the set $S$ as $R$-basis. We
obtain an isomorphism of abelian groups
\begin{eqnarray}
\beta^G \colon \bbZ[\consub_f(G)] & \xrightarrow{\cong} & \underlineA(G),
\quad (H) ~ \mapsto ~ [G/H].
\label{iso_alpha^G:Z[consub_f(G)]_to_covA(G)}
\end{eqnarray}

The elementary proof of the following lemma is left to the reader.

\begin{lemma} \label{lem:orbits_and_Weyl_groups}
  Let $H$ and $K$ be subgroups of $G$.  Then

\begin{enumerate}

\item \label{lem:orbits_and_Weyl_groups:G/H^K} $G/H^K ~ = ~ \{gH
  \mid g^{-1}Kg \subset H\};$

\item \label{lem:orbits_and_Weyl_groups:WKbackslashG/H^K} The map
  $$\phi\colon G/H^K \rightarrow \consub(H), \hspace*{10mm} gH \mapsto
  g^{-1}Kg$$
  induces an injection
  $$W_GK\backslash(G/H^K) \rightarrow \consub(H);$$

\item \label{lem:orbits_and_Weyl_groups:WK isotropy} The
  $W_GK$-isotropy group of $gH \in G/H^K$ is $(gHg^{-1}\cap N_GK)/K$;

\item \label{lem:orbits_and_Weyl_groups:G/H^K and its WK-orbits} If
  $H$ is finite, then $G/H^K$ is a finite union of $W_GK$-orbits of
  the shape $W_GK/L$ for finite subgroups $L \subset W_GK$. \qed

\end{enumerate}
\end{lemma}

The next definition makes sense because of Lemma~\ref{lem:orbits_and_Weyl_groups} above.

\begin{definition}[$L^2$-character map] \label{def:L^2-character_map_for_the_Burnside_group}
  Define for a finite subgroup $K \subset G$ the \emph{$L^2$-character
    map at $(K)$ }
  $$\underlinecharacter_K^{G} \colon \underlineA(G) \rightarrow
  \bbQ, \quad [S] \mapsto \sum_{i=1}^r |L_i|^{-1}$$
  if $W_GK/L_1$,
  $W_GK/L_2$,\ldots , $W_GK/L_r$ are the $W_GK$-orbits of $S^K$.
  Define the \emph{global $L^2$-character map} by
  $$\underlinecharacter^{G} \colon   \underlineA(G) \to \bbQ[\consub_f(G)], \quad [S]
  ~ \mapsto \sum_{(K) \in \consub_f(G)}
  \underlinecharacter^G_K([S]) \cdot (K).$$
\end{definition}

Notice that one gets from Lemma~\ref{lem:orbits_and_Weyl_groups} the
following explicit formula for the value of
$\underlinecharacter^{G}_K(G/H)$. Namely, define
\begin{eqnarray*}
{\cal L}_K(H) & := & \{(L) \in \consub(H) \mid
L \mbox{ conjugate to } K \text{ in }G\}.
\end{eqnarray*}
For $(L) \in {\cal L}_K(H)$ choose $L \in (L)$ and $g \in G$ with
$g^{-1}Kg = L$. Then
\begin{eqnarray*}
g(H \cap N_GL)g^{-1} & = & gHg^{-1} \cap N_GK;
\\
|(gHg^{-1} \cap N_GK)/K|^{-1} & = & \frac{|K|}{|H \cap N_GL|}.
\end{eqnarray*}
This implies
\begin{eqnarray}
\label{value of chi^{Gamma}_K(Gamma/H))}
\underlinecharacter_K^{G}(G/H) & = & \sum_{(L) \in {\cal L}_K(H)}
\frac{|K|}{|H \cap N_GL|}.
\label{formula for character at (K)}
\end{eqnarray}

\begin{remark}[Burnside integrality relations] \label{rem:Burnside_integrality_relations} \em
  Let $T \subseteq \consub_f(G)$ be a finite subset closed under
  taking subgroups, i.e., if $(H) \in T$, then $(K) \in T$ for every
  subgroup $K \subseteq H$.  Since a finite subgroup contains only
  finitely many subgroups, one can write $\consub_f(G)$ as the union
  of such subsets $T$.  The union of two such subsets is again such a
  subset. So $R[\consub_f(G)]$ is the colimit of the finitely
  generated free $R$-modules $R[T]$ if $T$ runs to the finite subsets
  of $\consub_f(G)$ closed under taking subgroups.

  Fix a subset $T$ of $\consub_f(G)$ closed under taking subgroups.
  One easily checks using Lemma~\ref{lem:orbits_and_Weyl_groups} that
  the composition
  $$\bbZ[\consub_f(G)] \xrightarrow{\beta^G} \underlineA(G)
  \xrightarrow{\underlinecharacter^G} \bbQ[\consub_f(G)]$$
  maps
  $\bbZ[T]$ to $\bbQ[T]$. We numerate the elements in $T$ by $(H_1)$,
  $(H_2)$, $\ldots$, $(H_r)$ such that $H_i$ is subconjugated to
  $(H_j)$ only if $i \le j$ holds. Then the composition
  $$\bbQ[\consub_f(G)] \xrightarrow{\beta^G \otimes_{\bbZ} \bbQ}
  \underlineA(G) \otimes_{\bbZ} \bbQ
  \xrightarrow{\underlinecharacter^G} \bbQ[\consub_f(G)]$$
  induces
  a $\bbQ$-homomorphism
  $$\bbQ[T] \xrightarrow{A_T} \bbQ[T]$$
  given with respect to the
  basis $\{(H_i) \mid i = 1,2, \ldots r\}$ by a matrix $A$ which is
  triangular and all whose diagonal entries are equal to $1$.
  The explicit values of the entries in $A_T$ are given by \eqref{formula for character at (K)}.
  The matrix $A_T$ is invertible and one can actually write down an explicit formula for
  its inverse matrix $B_T$ in terms of M\"obius inversion \cite[Chaper
  IV]{Aigner(1979)}.  The matrix $B_T$ yields an isomorphism
  $$B_T \colon\bbQ[T] \xrightarrow{\cong} \bbQ[T].$$

  Given an element $x \in \bbQ[\consub_f(G)]$, we can find a finite
  subset $T \subseteq \consub_f(G)$ closed under taking subgroups such
  that $x$ lies already in $\bbQ[T]$. Then $x$ lies in the image of
  the injective $L^2$-character map
  $$\underlinecharacter^{G} \colon \underlineA(G) \to \bbQ[\consub_f(G)]$$
  of
  Definition~\ref{def:L^2-character_map_for_the_Burnside_group} if
  and only if
  $$B_T \colon\bbQ[T] \xrightarrow{\cong} \bbQ[T]$$
  maps $x$ to an
  element in $\bbZ[T]$. This means that the following rational numbers
  $$\sum_{j=1}^r B_T(i,j) \cdot x(j)$$
  for $i = 1,2 \ldots, r$ are
  integers, where $B_T(i,j)$ and $x(j)$ are the components of $B_T$
  and $x$ belonging to $(i,j)$ and $j$. We call the condition that these
  rational numbers are integral numbers the
  \emph{Burnside integrality relations}.

  Now suppose that $G$ is finite. Then the global $L^2$-character of
  Definition~\ref{def:L^2-character_map_for_the_Burnside_group} is
  related to the classical character map~\eqref{character_map_for_A(G)_for_finite_G}
  by the factors $|W\!K|^{-1}$, i.e., we have for each
  subgroup $K$ of $G$ and any finite $G$-set $S$
\begin{eqnarray}
\ch^{G}_K(S) & = & |W\!K|^{-1} \cdot |S^K|.
\label{relation between L^2-character map and classical character map}
\end{eqnarray}
One easily checks that for finite $G$  under the identification
(\ref{relation between L^2-character map and classical character map}) the Burnside
integrality relations can be reformulated as a set of
congruences, which consists of one congruence modulo $|W_GH|$ for every
subgroup $H \subseteq G$ (compare
Subsection~\ref{subsec:The_Character_Map_and_the_Burnside_Ring_Congruences}).
\em
\end{remark}


\subsection{Relation to $L^2$-Euler characteristic and Universal Property of the Covariant Burnside Group}
\label{subsec:Relation_to_$L^2$-Euler_characteristic_and_universal_property_of_the_Covariant_Burnside_Group}

  The Burnside group $\underlineA(G)$ can be characterized as the
  universal additive invariant for finite proper $G$-$CW$-complexes
  and the universal equivariant Euler characteristic of a finite
  proper $G$-$CW$-complex is mapped to the $L^2$-Euler characteristics
  of the $W_GH$-$CW$-complexes $X^H$ by the character map at $(H)$ for
  every finite subgroup $H \subseteq G$.  In particular it is
  interesting to investigate the universal equivariant Euler
  characteristic of the classifying space for proper $G$-actions
  $\underline{E}G$ provided that there is a finite $G$-$CW$-model for
  $\underline{E}G$.  All this is explained in \cite[Section
  6.6.2]{Lueck(2002)}.

  The relation of the universal equivariant Euler characteristic to the equivariant Euler class
  which is by definition  the class of the Euler operator on a cocompact proper
  smooth  G-manifold with $G$-invariant Riemannian metric in equivariant $K$-homology
  defined by Kasparov is analyzed in \cite{Lueck-Rosenberg(2003b)}. 
  Equivariant Lefschetz classes for $G$-maps of finite proper $G$-$CW$-complexes are studied
  in \cite{Lueck-Rosenberg(2003a)}.


\subsection{The Covariant Burnside Group and the Colimit-Version of the Burnside Ring Agree}
\label{subsec:The_Covariant_Burnside_Group_and_the_Colimit_Version_of_the_Burndside_Ring_Agree}

Instead of the inverse-limit-version one may also consider the colimit-version
$$\covA(G) ~:= \colim{H \in \SubGF{G}{\calfin}} A(H)$$
where we consider $A$ as a covariant functor from
$\SubGF{G}{\calfin}$ to the category of $\bbZ$-modules by
induction.

\begin{theorem}[$\covA(G)$ and $\underlineA(G)$ agree]
\label{the:covA(G)_and_underlineA(G)_agree}
There obvious map induced by the various inclusions
of a finite subgroup $H \subseteq G$
$$V^G \colon \covA(G) \xrightarrow{\cong}  \underlineA(G)$$
is a bijection of abelian groups.
\end{theorem}
\begin{proof}
Recall that $\underlineA(G)$ is the free abelian group with the set
$\consub_f(G)$ of conjugacy classes of finite subgroups as basis.
Now the claim follows from the bijection~\eqref{phi^g colimt consub(H) to consub_f(G)}.
\end{proof}

The analogue for the representation ring is an open conjecture.
Namely if we define for a field $F$ of characteristic zero
$$R_{\cov,F}(G) ~:= \colim{H \in \SubGF{G}{\calfin}} R_F(H)$$
we can consider
\begin{conjecture} \label{con:R_{cov,F}_cong_K_0(FG)}
The obvious map
$$W^{G,F} \colon R_{\cov,F}(G) \to K_0(F[G])$$
is a bijection of abelian groups.
\end{conjecture}
This conjecture follows from the Farrell-Jones Conjecture for algebraic $K$-theory for
$F[G]$ as explained in \cite[Conjecture~3.3]{Lueck-Reich(2004)}. No counterexamples are known at the time of writing.
For a status report about the Farrell-Jones Conjecture we refer for instance to
\cite[Section~5]{Lueck-Reich(2004)}.

Let
\begin{eqnarray}
\covP^G \colon \covA(G) & \to & R_{\cov,F}(G).
\label{covP}
\end{eqnarray}
be the map induced by the maps $P^H \colon A(H) \to R_F(H), [S] \mapsto [F[S]]$ 
for the various finite subgroups $H \subseteq G$.


\subsection{The Covariant Burnside Group and the Projective Class Group}
\label{subsec:The_Covariant_Burnside_Group_and_the_Projective_Class_Group}

Given a finite proper $G$-set, the $\bbQ$-module $\bbQ[S]$ with the set $S$ as basis
becomes a finitely generated projective $\bbQ G$-module by the $G$-action on $S$.
Thus we obtain a homomorphism
\begin{eqnarray}
\underlineP^G \colon \underlineA(G) & \to & K_0(\bbQ G).
\label{underlineP}
\end{eqnarray}

\begin{conjecture}
The map $\underlineP^G \colon \underlineA(G) \to K_0(\bbQ G)$ is rationally surjective.
\end{conjecture}

This conjecture is motivated by the fact that it is implied by
Theorem~\ref{the:The_Burnside_ring_and_the_rational_representation_ring_for_finite_groups} and
Theorem~\ref{the:covA(G)_and_underlineA(G)_agree} together
with Conjecture~\ref{con:R_{cov,F}_cong_K_0(FG)}.


\subsection{The Covariant Burnside Group as Module over the Finite-Set-Version}
\label{subsec:The_Covariant_Burnside_Group_as_Module_over_the_Finite-Set-Version}

If $S$ is a finite $G$-set and $T$ is a cofinite proper $G$-set, then their cartesian product
with the diagonal $G$-action is a cofinite proper $G$-set. Thus we obtain a pairing
\begin{eqnarray}
\mu_A^G \colon \finsetA(G) \times \underlineA(G) \to \underlineA(G), \quad ([S],[T]) ~ \mapsto ~ [S \times T]
\label{pairingfinset(A)_times_covA(G)_to_covA(G)}
\end{eqnarray}
Analogously one defines a pairing
\begin{eqnarray}
\mu_K^G \colon \Sw^f(G;\bbQ) \times  K_0(\bbQ G) \to K_0(\bbQ G) , \quad ([M],[P]) ~ \mapsto ~ [M \otimes_{\bbQ} P]
\label{pairing_Sw(G;Q)_times_K_0(QG)_to__K_0(QG)}
\end{eqnarray}
which turns $K_0(\bbQ G)$ into a $\Sw^f(G;\bbQ)$-module. These two pairings are compatible in the obvious sense
(see \eqref{diagram_about_pairings}).


\subsection{A Pairing between the Inverse-Limit-Version and the Covariant Burnside Group}
\label{subsec:A_Pairing_between_the_Inverse-Limit-Version_and_the_Covariant_Burnside_Group}

Given a finite group $H$, we obtain a homomorphism of abelian groups
$$\nu^H \colon A(H) \to \hom_{\bbZ}(A(H),\bbZ), \quad [S] \mapsto \nu^H(S),$$
where $\nu^H(S) \colon A(H) \to \bbZ$ maps $[T]$ to $|G\backslash(S \times T)|$
for the diagonal $G$-operation on $S \times T$. A group homomorphism
$f \colon H \to K$ induces a homomorphism of abelian groups $\res_{\phi} \colon A(K) \to A(H)$ by restriction
and a homomorphism of abelian groups $\ind_{\phi} \colon A(H) \to A(K)$ by induction. The latter induces
a homomorphism of abelian groups $(\ind_{\phi})^* \colon A(K) \to A(H)$. One easily checks
that the collection of the homomorphisms $\nu^H$ for the subgroups $H \subseteq G$ induces
a natural transformation of the contravariant functors from $\SubGF{G}{\calfin}$ to $\bbZ$-modules
given by $A(?)$ and $\hom_{\bbZ}(A(?),\bbZ)$. Passing to the inverse limit, the canonical isomorphism of
abelian groups
$$\hom_{\bbZ}\left(\colim{H \in \SubGF{G}{\calfin}} A(H), \bbZ\right)
~ \xrightarrow{\cong} ~
\invlim{H \in\SubGF{G}{\calfin}} \hom_{\bbZ}(A(H),\bbZ)
$$
and the isomorphism appearing in Theorem~\ref{the:covA(G)_and_underlineA(G)_agree} yield
a homomorphism of abelian groups
$$\nu_A^G \colon \invA(G) \to \hom_{\bbZ}(\underline{A}(G);\bbZ)$$
which we can also write a bilinear pairing
\begin{eqnarray}
\nu_A^G \colon \invA(G) \times  \underlineA(G) &  \to & \bbZ.
\label{pairing_nu_A invA(G)_times_underlineA(G)_to_Z}
\end{eqnarray}

For a field $F$ of characteristic zero, there is an analogous pairing
\begin{eqnarray}
\nu_R^G \colon R_{\inv,F}(G) \times R_{\cov,F}(G) & \to & \bbZ
\label{pairing_nu_K}
\end{eqnarray}
which comes from the various homomorphisms of abelian groups for each finite subgroup $H \subseteq G$
$$R_F(H) \to \hom_{\bbZ}(R_F(H);\bbZ), [V] ~ \mapsto ~ \left([W] \mapsto \dim_F(F \otimes_{FG} (V \otimes_F W))\right).$$
The pairings $\nu_A^G$ and $\nu_R^G$ are compatible in the obvious sense
(see \eqref{diagram_about_pairings}).

The homomorphism of abelian groups $\overline{\nu_R^G} \colon
R_{\inv,F}(G) \to \hom_{\bbZ}(R_{\cov,F}(G),\bbZ)$ associated to
the pairing $\nu_R^G$ is injective. Its cokernel is finite if $G$
has only finitely many conjugacy classes of finite subgroups. It
is rationally surjective if there is an upper bound on the orders
of finite subgroups of $G$.

Define the homomorphism $Q^G_A \colon \underlineA(G) \to \bbZ$ and
$Q^G_K \colon K_0(\bbQ G) \to \bbZ$  respectively by sending $[S]$ to $|G\backslash S|$ and 
$[P]$ to $\dim_{\bbQ}(\bbQ \otimes_{\bbQ G} P)$ respectively.
Then the pairings $\mu_A^G$, $\nu_A^G$, $\mu_K^G$ and $\nu_R^G$ are compatible in the obvious sense
(see \eqref{diagram_about_pairings}).


\subsection{Some Computations of the Covariant Burnside Group}
\label{subsec:Some_Computations_of_the_Covariant_Burnside_Group}

\begin{example}[$\underlineA(G)$ for torsionfree $G$]
\label{exa:covA(G)_for_torsionfree_G} \em Suppose that $G$ is
torsionfree. Then the inclusion $i\colon \{1\} \to G$ induces a
$\bbZ$-isomorphism
$$
\underlineA(i) \colon \underlineA(\{1\}) = \bbZ ~
\xrightarrow{\cong} ~ \underlineA(G).$$
\em
\end{example}

\begin{example}[Extensions of $\bbZ^n$ with $\bbZ/p$ as quotient]
\label{exa:extensions_by_Z/p_for_covA} \em
Suppose that $G$ can be written as an extension $1 \to A \to G \to
\bbZ/p \to 1$ for some fixed prime number $p$ and for $A = \bbZ^n$ for
some integer $n \ge 0$ and that $G$ is not torsionfree.  We use the
notation of Example~\ref{exa:extensions_by_Z/p} in the sequel. We
obtain an exact sequence
$$0 \to \bigoplus_{H^1(\bbZ/p;A)} \underlineA(\{1\}) \to
\underlineA(\{1\}) \oplus \bigoplus_{H^1(\bbZ/p;A)}
\underlineA(\bbZ/p) \to \underlineA(G) \to 0$$
If we define
$\widetilde{\underlineA}(G)$ as the kernel of $\underlineA(G) \to
\underlineA(\{1\})$, we obtain an isomorphism
$$\bigoplus_{H^1(\bbZ/p;A)} \widetilde{A}(\bbZ/p)
\xrightarrow{\cong} \widetilde{\underlineA}(G).$$
Let $H_0$ be the
trivial subgroup and $H_1$, $H_2$, $\ldots $, $H_r$ be a complete set
of representatives of the conjugacy classes of finite subgroups.
Then $r = |H^1(\bbZ/p;A)|$ and $\underlineA(G)$ is the free abelian group of rank $r+1$
with $\{[G/H_0], [G/H_1], \ldots [G/H_r]\}$ as basis. Each
$H_i$ is isomorphic to $\bbZ/p$.  We compute using \eqref{formula for
  character at (K)}
$$\begin{array}{lcllll} \ch^{G}_{H_0}(G/H_0) & = & 1; & & &
  \\
  \ch^{G}_{H_0}(G/H_j) & = & \frac{1}{p} & & & j = 1,2, \ldots , r;
  \\
  \ch^{G}_{H_i}(G/H_j) & = & 1 & & & i = j, ~ i,j = 1,2, \ldots , r;
  \\
  \ch^{G}_{H_i}(G/H_j) & = & 0 & & & i \not= j, ~ i,j = 1,2, \ldots ,
  r.
\end{array}$$
The Burnside integrality conditions
(see Remark~\ref{rem:Burnside_integrality_relations}) become in this case for $x =
(x(i)) \in \bigoplus_{i = 0}^r \bbQ$
\begin{eqnarray*}
x(0) - \frac{1}{p} \cdot \sum_{i = 1}^r x(i) & \in & \bbZ;
\\
x(i) & \in &\bbZ \quad i = 1, 2 ,\ldots , r.
\end{eqnarray*}
\em
\end{example}

\begin{example}[Groups with appropriate maximal finite subgroups]
\label{exa:conditions_M and_NM_for_covA(G)} \em
Consider the groups appearing in Example~\ref{exa:conditions_M_and_NM}.
In the notation of Example~\ref{exa:conditions_M_and_NM} we
get an isomorphism of $\bbZ$-modules
$$\bigoplus_{i \in I} \widetilde{A}(M_i)
\xrightarrow{\bigoplus_{i \in I} \widetilde{\underlineA}(j_i)}
\widetilde{\underlineA}(G).$$
\em
\end{example}


\typeout{--------- Section 5: Equivariant Cohomology Theories-----------}

\section{Equivariant Cohomology Theories}
\label{sec:Equivariant_Cohomology_Theories}

In this section we recall the axioms
of a (proper) equivariant cohomology theory of \cite{Lueck(2004i)}.
 They are dual to the ones
of a (proper) equivariant homology theory as described in
\cite[Section 1]{Lueck(2002b)}.

\subsection{Axiomatic Description of a $G$-Cohomology Theory}
\label{subsec:Axiomatic_Description_of_a_G-Cohomology_Theory}

Fix a group $G$ and an commutative ring $R$.
A $G$-$CW$-pair $(X,A)$ is a pair of $G$-$CW$-complexes.
Recall that a $G$-$CW$-complex $X$ is proper if and only if all isotropy groups of $X$ are finite,
and is finite  if $X$ is obtained from $A$ by attaching
finitely many equivariant cells, or, equivalently, if
$G\backslash X$ is compact.
A \emph{$G$-cohomology theory $\calh^*_G$ with values
in $R$-modules} is a collection of
covariant functors $\calh_G^n$ from the category of
$G$-$CW$-pairs to the category of
$R$-modules indexed by $n \in \bbZ$ together with natural transformations
$\delta^n_G(X,A)\colon \calh^n_G(X,A) \to
\calh^{n+1}_G(A):= \calh^{n+1}_G(A,\emptyset)$ for $n \in \bbZ$
such that the following axioms are satisfied:
\begin{itemize}

\item $G$-homotopy invariance \\[1mm]
If $f_0$ and $f_1$ are $G$-homotopic maps $(X,A) \to (Y,B)$
of $G$-$CW$-pairs, then $\calh^n_G(f_0) = \calh^n_G(f_1)$ for $n \in \bbZ$;

\item Long exact sequence of a pair\\[1mm]
Given a pair $(X,A)$ of $G$-$CW$-complexes,
there is a long exact sequence
$$\cdots \xrightarrow{\delta^{n-1}_G}
\calh^{n}_G(X,A) \xrightarrow{\calh^{n}_G(j)}
\calh^n_G(X) \xrightarrow{\calh^{n}_G(i)} \calh^n_G(A)
\xrightarrow{\delta^n_G} \cdots,$$
where $i\colon A \to X$ and $j\colon X \to (X,A)$ are the inclusions;

\item Excision \\[1mm]
Let $(X,A)$ be a $G$-$CW$-pair and let
$f\colon A \to B$ be a cellular $G$-map of
$G$-$CW$-complexes. Equip $(X\cup_f B,B)$ with the induced structure
of a $G$-$CW$-pair. Then the canonical map
$(F,f)\colon (X,A) \to (X\cup_f B,B)$ induces an isomorphism
$$\calh^n_G(F,f)\colon \calh^n_G(X,A) \xrightarrow{\cong}
\calh^n_G(X\cup_f B,B);$$

\item Disjoint union axiom\\[1mm]
Let $\{X_i\mid i \in I\}$ be a family of
$G$-$CW$-complexes. Denote by
$j_i\colon  X_i \to \coprod_{i \in I} X_i$ the canonical inclusion.
Then the map
$$\prod_{i \in I} \calh^{n}_G(j_i)\colon  \calh^n_G\left(\coprod_{i \in I} X_i\right)
\xrightarrow{\cong}  \prod _{i \in I} \calh^n_G(X_i)
$$
is bijective.

\end{itemize}

If $\calh^*_G$ is defined or considered only for proper $G$-$CW$-pairs
$(X,A)$, we call it a
\emph{proper $G$-cohomology theory $\calh^*_G$ with values
in $R$-modules}.

\subsection{Axiomatic Description of an Equivariant Cohomology Theory}
\label{subsec:Axiomatic_Description_of_a_Equivariant_Cohomology_Theory}

Let $\alpha\colon H \to G$ be a group homomorphism.
Given an $H$-space $X$, define the \emph{induction of $X$ with $\alpha$}
to be the $G$-space $\ind_{\alpha} X$ which  is the quotient of
$G \times X$ by the right $H$-action
$(g,x) \cdot h := (g\alpha(h),h^{-1} x)$
for $h \in H$ and $(g,x) \in G \times X$.
If $\alpha\colon H \to G$ is an inclusion, we also write $\ind_H^G$ instead of
$\ind_{\alpha}$.

A \emph{(proper) equivariant cohomology theory
$\calh^*_?$ with values in $R$-modules}
consists of a collection of (proper)
$G$-cohomology theories $\calh^*_G$ with values in $R$-modules for each group $G$
together with the following so called \emph{induction structure}:
given a group homomorphism $\alpha\colon  H \to G$ and  a (proper) $H$-$CW$-pair
$(X,A)$ there are for each $n \in \bbZ$ natural homomorphisms
\begin{eqnarray}
\ind_{\alpha}\colon \calh^n_G(\ind_{\alpha}(X,A))
&\to &
\calh^n_H(X,A)\label{induction structure}
\end{eqnarray}
satisfying

\begin{enumerate}

\item Bijectivity\\[1mm]
If $\ker(\alpha)$ acts freely on $X$, then
$\ind_{\alpha}\colon \calh^n_G(\ind_{\alpha}(X,A))
\to \calh^n_H(X,A)$
is bijective for all $n \in \bbZ$;

\item Compatibility with the boundary homomorphisms\\[1mm]
$\delta^n_H \circ \ind_{\alpha} = \ind_{\alpha} \circ \delta^n_G$;

\item Functoriality\\[1mm]
Let $\beta\colon  G \to K$ be another group
homomorphism. Then we have for $n \in \bbZ$
$$\ind_{\beta \circ \alpha} ~ = ~
\ind_{\alpha} \circ \ind_{\beta} \circ \calh^n_K(f_1) \colon
\calh^n_H(\ind_{\beta\circ\alpha}(X,A)) \to \calh^n_K(X,A),$$
where $f_1\colon \ind_{\beta}\ind_{\alpha}(X,A)
\xrightarrow{\cong} \ind_{\beta\circ \alpha}(X,A),
\hspace{1mm} (k,g,x) \mapsto (k\beta(g),x)$
is the natural $K$-homeo\-mor\-phism;

\item Compatibility with conjugation\\[1mm]
For $n \in \bbZ$, $g \in G$ and a (proper) $G$-$CW$-pair $(X,A)$
the homomorphism $\ind_{c(g)\colon G \to G}\colon
\calh^n_G(\ind_{c(g)\colon G \to G}(X,A)) \to \calh^n_G(X,A)$ agrees with
$\calh^n_G(f_2)$, where $f_2$ is the $G$-homeomorphism
$f_2\colon (X,A) \to \ind_{c(g)\colon  G \to G} (X,A), \hspace{1mm}
x \mapsto(1,g^{-1}x)$ and $c(g)(g') = gg'g^{-1}$.

\end{enumerate}

This induction structure links the
various $G$-cohomology theories for different groups $G$.

 Sometimes we will need the following lemma
 whose elementary proof is analogous to the
 one in \cite[Lemma 1.2]{Lueck(2002b)}.
 \begin{lemma}
 \label{lem:calh_G(G/H)_and_calh_H(*)}
 Consider finite subgroups $H,K \subseteq G$ and an element
 $g \in G$ with $gHg^{-1} \subseteq K$.
 Let $R_{g^{-1}}\colon G/H \to G/K$ be the $G$-map
 sending $g'H$ to $g'g^{-1}K$ and
 $c(g)\colon H \to K$ be the homomorphism sending $h$ to $ghg^{-1}$.
 Let $\pr\colon (\ind_{c(g)\colon H \to K}\pt) \to \pt$ be the projection
 to the one-point space $\pt$.
 Then the following diagram commutes
 $$\begin{CD}
 \calh_G^n(G/K)   @> \calh^n_G(R_{g^{-1}}) >> \calh_G^n(G/H)
 \\
 @ V \ind_K^G V \cong V @V \ind_H^G V\cong V
 \\
 \calh^n_K(\pt) @>  \ind_{c(g)} \circ \calh_K^n(\pr) >>  \calh^n_H(\pt)
 \end{CD}
 $$
 \end{lemma}


\subsection{Multiplicative Structures}
\label{subsec:Multiplicative_Structures}

Let  $\calh^*_G$ be a (proper) $G$-cohomology theory. A \emph{multiplicative structure}
assigns to a  (proper) $G$-$CW$-complex $X$ with
$G$-$CW$-subcomplexes $A,B \subseteq X$ natural $R$-homomorphisms
\begin{eqnarray}
\cup \colon \calh^{m}_G(X,A) \otimes_R \calh^{n}_G(X,B) & \to &
\calh^{m+n}_G(X,A\cup B).
\label{equivariant cup -product}
\end{eqnarray}
This product is required to be compatible with the boundary
homomorphism of  the long exact sequence of a $G$-$CW$-pair,
namely,
for $u \in \calh^m_G(A)$ and $v\in \calh^n{G}(X)$ and $i \colon A \to X$ the inclusion we have
$\delta(u \cup  v) = \delta(u \cup \calh^n(i)(v)$.
Moreover, it is required to be graded commutative, to be associative and to have a unit $1 \in \calh^0_G(X)$ for every (proper)
$G$-$CW$-complex $X$.

Let $\calh^*_?$ be a (proper) equivariant cohomology theory. A
\emph{multiplicative structure} on it assigns a multiplicative
structure to the associated (proper) $G$-coho\-mo\-logy theory $\calh^*_G$ for
every group $G$ such that for each group homomorphism
$\alpha \colon H \to G$ the maps given by the induction structure of \eqref{induction structure}
\begin{eqnarray*}
\ind_{\alpha}\colon \calh^n_G(\ind_{\alpha}(X,A))
&\xrightarrow{\cong} &
\calh^n_H(X,A)
\end{eqnarray*}
are in the obvious way compatible with the multiplicative structures on $\calh^*_G$ and
$\calh^*_H$.

\begin{example}{\bf Equivariant cohomology theories coming from non-equi\-va\-riant ones).}
 \label{exa:equivariant_cohomology_theories}
\em
Let $\calk^*$ be a (non-equivariant) cohomology theory with multiplicative structure, for instance
singular cohomology or topological $K$-theory. We can assign to it an equivariant cohomology theory
with multiplicative structure $\calh^*_?$ in two ways. Namely, for a group $G$ and a pair
of $G$-$CW$-complexes $(X,A)$ we define $\calh^n_G(X,A)$ by
$\calk^n(G\backslash (X,A))$ or by $\calk^n(EG \times_G(X,A))$.
\em
\end{example}


\subsection{Equivariant Topological $K$-Theory}
\label{subsec:Equivariant_Topological_K-Theory}

In \cite{Lueck-Oliver(2001a)} equivariant topological $K$-theory is defined for
finite proper equivariant $CW$-complexes in terms of equivariant vector bundles.
For finite $G$ it reduces to the classical notion which appears for instance in
\cite{Atiyah(1967)}. Its relation to equivariant $KK$-theory is
explained in \cite{Phillips(1988)}.
This definition is extended to (not necessarily finite) proper equivariant $CW$-complexes in
\cite{Lueck-Oliver(2001a)} in terms of equivariant spectra using $\Gamma$-spaces.
This equivariant cohomology theory $K_?^*$ has the property that for any finite subgroup $H$ of a group
$G$ we have
$$
K^n_G(G/H) ~ = ~ K^n_H(\pt) ~ = ~
\left\{
\begin{array}{lll}
R_{\bbC}(H) & & n \text{ even};
\\
\{0\}       & & n \text{ odd}.
\end{array}
\right.
$$


\typeout{--------- Section 5: Equivariant Cohomology Theories-----------}

\section{Equivariant Stable Cohomotopy in Terms of Real Vector Bundles}
\label{sec:Equivariant_Stable_Cohomotopy_in_Terms_of_Real_Vector_Bundles}

In this section we give a construction  of equivariant stable cohomotopy for
finite proper $G$-$CW$-complexes for infinite groups $G$ in terms of real vector
bundles and maps between the associated sphere bundles.  The result will be an equivariant
cohomology theory with multiplicative structure for finite proper equivariant $CW$-complexes.
It generalizes the well-known approach for finite groups in terms of representations.
We will first give the construction, show why it reduces to the classical
construction for a finite group and explain why we need to consider equivariant vector bundles
and not only representations in the case of an infinite group.


\subsection{Preliminaries about Equivariant Vector Bundles}
\label{subsec:Preliminaries_about_Equivariant_Vector_Bundles}

We will need the following notation.  Given a finite-dimensional
(real) vector space $V$, we denote by $S^V$ its \emph{one-point
  compactification}. We will use the point at infinity as the base
point of $S^V$ in the sequel. Given two finite-dimensional vector
spaces $V$ and $W$, the obvious inclusion $V \oplus W \to S^V \wedge
S^W$ induces a natural homeomorphism
\begin{eqnarray}
\phi(V,W) \colon S^{V \oplus W} &\xrightarrow{\cong} & S^V \wedge S^W.
\label{phi(V,W):S^{V oplus W}_to_S^V_wedge_S^W}
\end{eqnarray}
Let
\begin{eqnarray}
\nabla \colon S^{\bbR} & \to  & S^{\bbR} \vee S^{\bbR}.
\label{nabla:S^R_to_S^R_vee_S^R}
\end{eqnarray}
be the pinching map, which sends $x > 0$ to $\ln(x) \in \bbR \subseteq
S^{\bbR}$ in the first summand, $x < 0$ to $- \ln(-x) \in \bbR \subseteq
S^{\bbR}$ in the second summand and $0$ and $\infty$ to the base point
in $S^{\bbR} \vee S^{\bbR}$.  This is under the identification
$S^{\bbR} = S^1$ the standard pinching map $S^1 \to S^1/S^0 \cong S^1
\vee S^1$, at least up to pointed homotopy.

We need some basics about $G$-vector bundles over proper
$G$-$CW$-complexes.
More details can be found for instance in
\cite[Section 1]{Lueck-Oliver(2001a)}.
Recall that a $G$-$CW$-complex is proper
if and only if all its isotropy groups are finite. A \emph{$G$-vector bundle}
$\xi \colon E \to X$ over $X$ is a real vector bundle with a
$G$-action on $E$ such that $\xi$ is $G$-equivariant and for each $g
\in G$ the map $l_g \colon E \to E$ is fiberwise a linear isomorphism.
Such a $G$-vector bundle is automatically trivial in the equivariant
sense that for each $x \in X$ there is a $G$-neighborhood $U$, a
$G$-map $f \colon U \to G/H$ and a $H$-representation $V$ such that
$\xi|_U$ is isomorphic as $G$-vector bundle to the pullback of the
$G$-vector bundle $G \times_HV \to G/H$ by $f$. Denote the fiber
$\xi^{-1}(x)$ over a point $x \in X$ by $E_x$. This is a
representation of the finite isotropy group $G_x$ of $x \in X$.  A
\emph{map of $G$-vector bundles} $(\overline{f},f)$ from $\xi_0 \colon
E_0 \to X_0$ to $\xi_1 \colon E_1 \to X_1$ consists of $G$-maps
$\overline{f} \colon E_0 \to E_1$ and $f \colon X_0 \to X_1$ with
$\xi_1 \circ \overline{f} = f \circ \xi_0$ such that $\overline{f}$ is
fiberwise a (not necessarily injective or surjective) linear map.

Given a $G$-vector bundle $\xi \colon E \to X$, let $S^{\xi} \colon
S^E \to X$ be the locally trivial $G$-bundle whose fiber over $x \in
X$ is $S^{E_x}$. Consider two $G$-vector bundles $\xi \colon E \to X$
and $\mu \colon F \to X$. Let $S^{\xi} \wedge_X S^{\mu} \colon S^E
\wedge_X S^F \to X$ be the $G$-bundle whose fiber of $x \in X$ is
$S^{E_x} \wedge S^{F_x}$, in other word it is obtained from $S^{\xi}
\colon S^{E} \to X$ and $S^{\mu} \colon S^F \to X$ by the fiberwise
smash product.  Define $S^{\xi} \vee_X S^{\mu}$ analogously. 
From \eqref{phi(V,W):S^{V oplus W}_to_S^V_wedge_S^W}
we obtain a natural $G$-bundle isomorphism
\begin{eqnarray}
\phi(\xi,\mu) \colon S^{\xi \oplus \mu} &\xrightarrow{\cong} & S^{\xi} \wedge_X S^{\mu}.
\label{phi(xi,mu):S^{xi oplus mu}_to_S^xi_wedge_S^mu}
\end{eqnarray}

The next basic lemma is proved in \cite[Lemma~3.7]{Lueck-Oliver(2001a)}.

\begin{lemma} \label{lem:extending_vector_bundles_in_LO1}
Let $f \colon X \to Y$ be a $G$-map between finite proper $G$-$CW$-complexes
and $\xi$ a $G$-vector bundle over $X$. Then there is a $G$-vector bundle
$\mu$ over $Y$ such that $\xi$ is a direct summand in $f^*\mu$.
\end{lemma}


\subsection{The Definition of Equivariant Stable Cohomotopy Groups}
\label{subsec:The_Definition_of_Equivariant_Stable_Cohomotopy_Groups}

Fix a proper $G$-$CW$-complex $X$. Let $\SPHB^G(X)$ be the following
category.  Objects are $G$-$CW$-vector bundles $\xi \colon E \to X$
over $X$. A morphism from $\xi \colon E \to X$ to $\mu \colon F \to X$
is a bundle map $u \colon S^{\xi} \to S^{\mu}$ covering the identity
$\id \colon X \to X$ which fiberwise preserve the base points.  (We do
not require that $u$ is fiberwise a homotopy equivalence.)

A \emph{homotopy $h$ between two morphisms} $u_0,u_1$ from $\xi \colon
E \to X$ to $\mu \colon F \to X$ is a bundle map $h \colon S^{\xi}
\times [0,1] \to S^{\mu}$ from the bundle $S^{\xi} \times \id_{[0,1]}
\colon S^E \times [0,1] \to X \times [0,1]$ to $S^{\mu}$ which covers the
projection $X \times [0,1] \to X$ and fiberwise preserve the base points such that its restriction to $X
\times \{i\}$ is $u_i$ for $i = 0,1$.

Let $\underline{\bbR^k}$ be the
trivial vector bundle $X \times \bbR^k \to X$. We consider it as a
$G$-vector bundle using the trivial $G$-action on $\bbR^k$. Fix an integer $n \in \bbZ$.
Given two objects $\xi_i$, two non-negative integers $k_i$ with $k_i + n \ge 0$
and two morphisms
$$u_i \colon S^{\xi_i \oplus \underline{\bbR^{k_i}}} \to
S^{\xi_i \oplus \underline{\bbR^{k_i + n}}}$$
for $i = 0,1$, we call $u_0$ and $u_1$
\emph{equivalent}, if there are objects $\mu_i$ in $\SPHB^G(X)$ for $i = 0,1$ and an
isomorphism of $G$-vector bundles
$v \colon \mu_0 \oplus \xi_0  \xrightarrow{\cong} \mu_1  \oplus \xi_1$ such that the following
diagram in $\SPHB^G(X)$ commutes up to homotopy
$$
\begin{CD}
S^{\underline{\mu_0 \oplus \bbR^{k_1}}} \wedge_X S^{\xi_0  \oplus \underline{\bbR^{k_0}}}
@> \id  \wedge_X u_0 >>
S^{\underline{\mu_0 \oplus \bbR^{k_1}}} \wedge_X S^{\xi_0  \oplus \underline{\bbR^{k_0+n}}}
\\
@V \sigma_1 VV @V \sigma_2 VV
\\
S^{\mu_0 \oplus \xi_0 \oplus \underline{\bbR^{k_0 + k_1}}} & &
S^{\mu_0 \oplus \xi_0 \oplus \underline{\bbR^{k_0 + k_1 + n}}}
\\
@V S^{v \oplus \id}   VV @V S^{v \oplus \id}   VV
\\
S^{\mu_1 \oplus \xi_1 \oplus \underline{\bbR^{k_0 + k_1}}} & &
S^{\mu_1 \oplus \xi_1 \oplus \underline{\bbR^{k_0 + k_1 + n}}}
\\
@V \sigma_3 VV @V \sigma_4 VV
\\
S^{\mu_1  \oplus\underline{\bbR^{k_0}}} \wedge_X S^{\xi_1  \oplus \underline{\bbR^{k_1}}}
@> \id\wedge_X u_1 >>
S^{\mu_1  \oplus\underline{\bbR^{k_0}}} \wedge_X S^{\xi_1  \oplus \underline{\bbR^{k_1+n}}}
\end{CD}
$$
where $\sigma_i$ stands for the obvious isomorphism coming from
\eqref{phi(xi,mu):S^{xi oplus mu}_to_S^xi_wedge_S^mu} and permutation.

We define $\pi^n_G(X)$ to be the set of equivalence classes of such
morphisms $u \colon S^{\xi \oplus \underline{\bbR^k}} \to S^{\xi \oplus \underline{\bbR^{k+n}}}$
under the equivalence relation mentioned above. It becomes an abelian
group as follows.

The zero element is represented by the class of any morphism
$c \colon S^{\xi \oplus \underline{\bbR^k}} \to S^{\xi \oplus \underline{\bbR^{k+n}}}$
which is fiberwise the constant
map onto the base point.

Consider classes $[u_0]$ and $[u_1]$ represented by two morphisms of
the shape $u_i \colon S^{\xi_i \oplus \underline{\bbR^{k_i}}} \to
S^{\xi_i \oplus \underline{\bbR^{k_i + n}}}$ for $i = 0,1$.  Define their sum by the class
represented by the morphism
\begin{multline*}
  S^{\xi_0 \oplus \xi_1 \oplus \underline{\bbR^{k_0 + k_1 +1}}}
  \xrightarrow{\sigma_1}
  S^{\xi_0 \oplus \underline{\bbR^{k_0}}} \wedge_X S^{\xi_1
\oplus \underline{\bbR^{k_1}}} \wedge_X  S^{\underline{\bbR}}
  \\
  \xrightarrow{\id \wedge_X \id \wedge_X  \underline{\nabla}}
  S^{\xi_0 \oplus \underline{\bbR^{k_0}}} \wedge_X S^{\xi_1  \oplus \underline{\bbR^{k_1}}} \wedge_X
  \left(S^{\underline{\bbR}} \vee_X S^{\underline{\bbR}}\right)
  \\
  \xrightarrow{\tau}
  \left(S^{\xi_0 \oplus \underline{\bbR^{k_0}}} \wedge_X S^{\xi_1  \oplus \underline{\bbR^{k_1}}}
  \wedge_X S^{\underline{\bbR}}\right) \vee_X
  \left(S^{\xi_0 \oplus \underline{\bbR^{k_0}}} \wedge_X S^{\xi_1  \oplus \underline{\bbR^{k_1}}}
  \wedge_X S^{\underline{\bbR}}\right)
  \\
  \xrightarrow{\left(u_0 \wedge_X \id \wedge_X \id\right) \vee_X
    \left(\id \wedge_X u_1 \wedge_X\id\right)}
  \\
  \left(S^{\xi_0 \oplus \underline{\bbR^{k_0 +n}}} \wedge_X S^{\xi_1  \oplus \underline{\bbR^{k_1}}}
  \wedge_X S^{\underline{\bbR}}\right) \vee_X
  \left(S^{\xi_0 \oplus \underline{\bbR^{k_0}}} \wedge_X S^{\xi_1  \oplus \underline{\bbR^{k_1 +n }}}
  \wedge_X S^{\underline{\bbR}}\right)
  \\
   \xrightarrow{\sigma_5 \vee_X \sigma_6}
   \left(S^{\xi_0 \oplus \xi_1 \oplus \underline{\bbR^{k_0 + k_1 + 1 + n}}}\right)
   \vee_X
   \left(S^{\xi_0 \oplus \xi_1 \oplus \underline{\bbR^{k_0 + k_1 + 1 +n}}}\right)
   \\
   \xrightarrow{\id  \vee_X \id}
   S^{\xi_0 \oplus \xi_1 \oplus \underline{\bbR^{k_0 + k_1 + 1 +n}}}
\end{multline*}
where  the isomorphisms $\sigma_i$ are given by permutation and the
isomorphisms~\eqref{phi(xi,mu):S^{xi oplus mu}_to_S^xi_wedge_S^mu}
$\tau$ is given by the distributivity law for smash
and wedge-products and $\underline{\nabla}$ is defined
fiberwise by the map $\nabla$ of \eqref{nabla:S^R_to_S^R_vee_S^R}.

Consider a class $[u]$ represented by the morphisms of the shape
$u \colon S^{\xi \oplus \underline{\bbR^k}} \to S^{\xi \oplus \underline{\bbR^{k+n}}}$.
Define its inverse as the class represented by the composition
$$
S^{\xi \oplus  \underline{\bbR^{k+1}}}
\xrightarrow{\sigma_1} S^{\xi  \oplus \underline{\bbR^k}}
\wedge_X S^{\underline{\bbR}}
\xrightarrow{u \wedge_X \underline{-\id}} S^{\xi \oplus \underline{\bbR^{k+n}}} \wedge_X
S^{\underline{\bbR}}
\xrightarrow{\sigma_2} S^{\xi \oplus \underline{\bbR^{k+1+n}}}, $$
where $\underline{-\id}$ is fiberwise the map  $- \id \colon \bbR
\to \bbR$.  The proof that this defines the structure of an abelian
group is essentially the same as the one that the abelian group structure on the
stable homotopy groups of a space is well-defined.

Next consider a pair $(X,A)$ of proper $G$-$CW$-complexes.  In order
to define the abelian group $\pi_G^n(X,A)$ we consider morphisms
$u\colon S^{\xi  \oplus \underline{\bbR^k}} \to S^{\xi \oplus \underline{\bbR^{k+n}}}$
with $k + n \ge 0$ in $\SPHB^G(X)$
such that $u$ is trivial over $A$, i.e.,  for every point $a \in A$ the
map $u_a \colon S^{\xi_a \oplus \bbR^k} \to S^{\xi_a \oplus \bbR^{k+n}}$ is the constant
map onto the base point.  In the definition of the equivalence
relation for such pairs we require that the homotopies of two
morphisms are stationary over $A$.  Then define $\pi_G^n(X,A)$ as the
set of equivalence classes of morphism
$u\colon S^{\xi  \oplus \underline{\bbR^k}} \to S^{\xi \oplus \underline{\bbR^{k+n}}}$
in $\SPHB^G(X)$ with $k + n \ge 0$ which are trivial over $A$
using this equivalence relation. The definition of the abelian group
structure goes through word by word.

Notice that in the definition of $\pi_G^n(X,A)$ we cannot use as in
the classical settings cones or suspensions since these contain
$G$-fixed points and hence are not proper unless $G$ is finite. The
properness is needed to ensure that certain basic facts about bundles
carry over to the equivariant setting.


\subsection{The Proof of the Axioms of an Equivariant Cohomology Theory with
Multiplicative Structure}
\label{subsec:The_Proof_of_the_Axioms_of_an_Equivariant_Cohomology_Theory_with_Multiplicative_Structure}

In this subsection we want to prove

\begin{theorem}[Equivariant stable cohomotopy in terms of equivariant vector bundles]
\label{the:Equivariant_Stable_Cohomotopy_in_terms_of_equivariant_vector_bundles}
Equivariant stable cohomotopy $\pi^*_?$ defines an equivariant cohomology theory
with multiplicative structure for finite proper equivariant $CW$-complexes.
For every finite subgroup $H$ of the group $G$  the abelian groups $\pi^n_G(G/H)$ and
$\pi^n_H$ are isomorphic for every $n \in \bbZ$ and the rings $\pi^0_G(G/H)$ and $\pi^0_H = A(H)$ are isomorphic.
\end{theorem}

Consider a $G$-map $f \colon (X,A) \to (Y,B)$ of pairs of proper
$G$-$CW$-complexes. Using the pullback construction one defines a
homomorphism of abelian groups
$$\pi_G^n(f) \colon \pi_G^n(Y,B) \to \pi_n^G(X,A).$$
Thus $\pi_G^n$
becomes a contravariant functor from the category of proper
$G$-$CW$-pairs to the category of abelian groups.

\begin{lemma} \label{lem:homotopy_invariance_of_pi^n_G}
  Let $f_0,f_1 \colon (X,A) \to (Y,B)$ be two $G$-maps of pairs of
  proper $G$-$CW$-complexes.  If they are $G$-homotopic, then
  $\pi_G^n(f_0) = \pi_G^n(f_1)$.
\end{lemma}
\begin{proof}
  By the naturality of $\pi_G^n$ it suffices to prove that $\pi_G^n(h)
  = \id$ holds for the $G$-map
  $$h \colon (X,A) \times [0,1] \to (X,A) \times [0,1], \quad (x,t)
  ~\mapsto ~ (x,0).$$
  Let the element $[u] \in \pi_G^n((X,A) \times [0,1])$ be given
  by the morphism $u \colon S^{\xi\oplus \underline{\bbR^k}} \to S^{\xi \oplus \underline{\bbR^{k+n}}}$ in
  $\SPHB^G(X\times [0,1])$ with $k + n \ge 0$ which is trivial over $X \times \{0,1\} \cup A \times [0,1]$. By
  \cite[Theorem~1.2]{Lueck-Oliver(2001a)} there is an isomorphism of
  $G$-vector bundles $v \colon \xi \xrightarrow{\cong} h^*\xi$
  covering the identity $\id \colon X \times [0,1] \to X \times [0,1]$
  such that $v$ restricted to $X \times \{0\}$ is the identity $\id
  \colon \xi|_{X \times \{0\}} \to \xi|_{X \times \{0\}}$.  The
  composition of morphisms in $\SPHB^G(X \times [0,1])$
  $$u' \colon S^{h^*\xi  \oplus \underline{\bbR^k}} \xrightarrow{S^{v^{-1} \oplus \id}}
  S^{\xi \oplus \underline{\bbR^k}} \xrightarrow{u}
  S^{\xi \oplus \underline{\bbR^{k+n}}}
  \xrightarrow{S^{v} \oplus \id} S^{h^*\xi \oplus \underline{\bbR^{k+n}}} $$
  has the property that its restriction to $X
  \times \{0\}$ agrees with the restriction of $h^*u$ to $X \times
  \{0\}$. Hence this composite $u'$ is homotopic to the morphism
  $h^*u$ itself. Namely, if we write $h^*\xi = i_0^*\xi \times [0,1]$ for
  $i_0 \colon X \to X \times [0,1], x \mapsto (x,0)$, then the
  homotopy is given at time $s \in [0,1]$ by the morphism
  $$S^{i_0^*\xi \oplus \underline{\bbR^k}} \times [0,1] \to S^{i_0^*\xi \oplus
    \underline{\bbR^{k+n}}} \times [0,1], \quad (z,t) ~ \mapsto ~ (\pr \circ
  u'(z,st),t)$$
  for $\pr \colon S^{i_0^*\xi \oplus \underline{\bbR^{k+n}}}
  \times [0,1] \to S^{i_0^*\xi \oplus \underline{\bbR^{k+n}}}$ the
  projection. Obviously this homotopy is stationary over $A$.  We
  conclude from the equivalence relation appearing in the definition
  of $\pi_G^n$ and the definition of $\pi_G^n(h)$ that
  $\pi_G^n(h)([u]) = [u'] =[u]$ holds.
\end{proof}

Next we define the suspension homomorphism
\begin{eqnarray}
\sigma^n_G(X,A)  \colon \pi_G^n(X,A) & \to & \pi_G^{n+1}((X,A)  \times ([0,1],\{0,1\})).
\label{suspension homomorphism}
\end{eqnarray}

For a $G$-vector bundle $\xi$ over $X$ let $\xi \times [0,1]$ be the
obvious $G$-vector bundle over $X \times [0,1]$, which is the same as
the pullback of $\xi$ for the projection $X \times [0,1] \to X$.
Consider a morphism $u \colon S^{\xi} \to S^{\mu}$ in $\SPHB^G(X)$ which
is trivial over $A$. Let
$$\sigma(u) \colon S^{\xi \times [0,1]} = S^{\xi} \times [0,1] ~\to ~
S^{(\mu \oplus \underline{\bbR}) \times [0,1]} = \left(S^{\mu}
  \wedge_X S^{\underline{\bbR}}\right) \times [0,1]$$
be the morphism in $\SPHB^G(X \times [0,1])$ given by
$$(z,t) \in S^{\xi} \times [0,1] ~ \mapsto ~
\left((u(z)\wedge e(t)\right),t) \in \left(S^{\mu} \wedge_X
  S^{\underline{\bbR}}\right) \times [0,1],$$
where $e \colon [0,1] \to S^{\underline{\bbR}}$ comes from the homeomorphism $(0,1) \to \bbR, t
~ \mapsto \ln(x) - \ln(1-x)$. The morphism $\sigma(u)$ is trivial over
$X \times \{0,1\} \cup A \times [0,1]$. We define the in $(X,A)$
natural homomorphism of abelian groups $\sigma^n_G(X,A)$ by sending
the class of $u$ to the class of $\sigma(u)$.

\begin{lemma} \label{sigma^n_G(X,A)_is_a_bijection}
  The homomorphism $\sigma^n_G(X,A)$ of \eqref{suspension homomorphism}
is bijective for all pairs of
  proper $G$-$CW$-complexes $(X,A)$.
\end{lemma}
\begin{proof}
  We want to construct an inverse
  $$\tau_G^{n+1} \colon \pi_G^{n+1}((X,A) \times ([0,1],\{0,1\})) ~
  \to ~ \pi_G^n(X,A)$$
  of $\sigma^n_G(X,A)$.  Consider two $G$-vector
  bundles $\xi$ and $\mu$ over $X$ and a morphism $v \colon S^{\xi
    \times [0,1]} \to S^{\mu \times [0,1]}$ in $\SPHB^G(X \times [0,1])$
  which is trivial over $A \times \{0,1\}$. Define a morphism in
  $\SPHB^G(X)$ which is trivial over $A$
  $$\tau(v) \colon S^{\xi \oplus \underline{\bbR}} = S^{\xi} \wedge_X
  S^{\underline{\bbR}} \to S^{\mu}$$
  by sending $(z,(e(t)) \in S^{\xi}
  \wedge_X S^{\underline{\bbR}}$ to $\pr \circ v(z,t)$ for $\pr \colon
  S^{\mu \times [0,1]} = S^{\mu} \times [0,1] \to S^{\mu}$ the
  projection, $e \colon [0,1] \to S^{\underline{\bbR}}$ the map defined above and
  $t \in [0,1]$. Next consider an element $[u] \in \pi_G^{n+1}(X,A)$
  represented by a morphism
  $u \colon S^{\xi\oplus \underline{\bbR^{k}}} \to S^{\xi \oplus
    \underline{\bbR^{k+n+1}}}$ in $\SPHB^G(X \times [0,1])$ for $k +n \ge 0$ which is
  trivial over $X \times \{0,1\} \cup A \times [0,1]$. Choose an
  isomorphism of $G$-vector bundles $v \colon \xi_0 \times [0,1]
  \xrightarrow{\cong} \xi$ which covers the identity on $X \times
  [0,1]$ and is the identity over $X \times \{0\}$, where $\xi_0$ is the
  restriction of $\xi$ to $X = X \times \{0\}$ (see
  \cite[Theorem~1.2]{Lueck-Oliver(2001a)}). Let
  $u' \colon S^{(\xi_0 \bigoplus \underline{\bbR^k}) \times [0,1]} \to
  S^{(\xi_0 \bigoplus \underline{\bbR^n})\times
    [0,1]}$ be the composition
  \begin{multline*}
  S^{(\xi_0 \bigoplus \underline{\bbR^k}) \times [0,1]} \xrightarrow{S^{v \oplus \id}}
  S^{(\xi \bigoplus \underline{\bbR^k})} \xrightarrow{u}
  S^{\xi \oplus \underline{\bbR^{k+n+1}}}
  \xrightarrow{S^{v^{-1} \oplus \id}}
  S^{(\xi_0 \bigoplus \underline{\bbR^{k+n+1}})\times [0,1]}.
  \end{multline*}
  Notice that $[u] = [u']$ holds in $\pi_G^{n+1}((X,A) \times
  ([0,1],\{0,1\}))$. Define $\tau_G^{n+1}([u])$ by the class of
  $\tau(u') \colon S^{\xi_0 \bigoplus \underline{\bbR^{k+1}}} \to S^{\xi_0 \bigoplus \underline{\bbR^{k+n+1}}}$.

  Consider $[u] \in \pi^n_G(X,A)$ represented by the morphism $u
  \colon S^{\xi \oplus \underline{\bbR^{k}}} \to S^{\xi \oplus \underline{\bbR^{k+n}}}$ in $\SPHB^G(X)$
  which is trivial over $A$. Then $\tau^{n+1}_G \circ \sigma^n_G([u])$
  is represented by the morphism $\tau \circ \sigma(u)$ which can be
  identified with
  $$S^{\xi \oplus \underline{\bbR^{k+1}}} \xrightarrow{\sigma_1} S^{\xi \oplus \underline{\bbR^{k}}}
  \wedge_X S^{\underline{\bbR}} \xrightarrow{u \wedge_X \id} S^{\xi
    \bigoplus \underline{\bbR^{k + n}}} \wedge_X S^{\underline{\bbR}}
  \xrightarrow{\sigma_2} S^{\xi \oplus \underline{\bbR^{k+1+n}} \oplus
    \underline{\bbR}}$$
  But the latter morphism represents the same
  class as $u$. This shows $\tau^{n+1}_G \circ
  \sigma^n_G([u]) = [u]$ and hence $\tau^{n+1}_G \circ \sigma^n_G =
  \id$. The proof of $\sigma^n_G([u]) \circ \tau^{n+1}_G = \id$ is
  analogous.
\end{proof}

So far we have only assumed that the $G$-$CW$-complex $X$ is proper.
In the sequel we will need additionally that it is finite since this condition
appears in Lemma~\ref{lem:extending_vector_bundles_in_LO1}.

\begin{lemma} \label{lem:excision_for_pi_G^n}
  Let $(X_1,X_0)$ be a pair of finite proper $G$-$CW$-complexes and $f
  \colon X_0 \to X_2$ be a cellular $G$-map of finite proper
  $G$-$CW$-complexes.  Define the pair of finite proper
  $G$-$CW$-complexes $(X,X_2)$ by the cellular $G$-pushout
  $$\comsquare{X_0}{f}{X_2}{}{}{X_1}{F}{X}$$
  Then the
  homomorphism
  $$\pi_G^n(F,f) \colon \pi_G^n(X,X_2) \xrightarrow{\cong}
  \pi_G^n(X_1,X_0)$$
  is bijective for all $n \in \bbZ$.
\end{lemma}
\begin{proof}
  We begin with surjectivity. Consider an element $a \in \pi_G^n(X_1,X_0)$ represented by a
  morphism
  $$u \colon S^{\xi \oplus \underline{\bbR^k}} \to S^{\xi \oplus \underline{\bbR^{k+n}}}$$
  in
  $\SPHB^G(X_1)$ which is trivial over $X_0$. By
  Lemma~\ref{lem:extending_vector_bundles_in_LO1} there is a $G$-vector bundle
  $\mu$ over $X$, a $G$-vector bundle $\xi'$ over $X_1$ and an
  isomorphism of $G$-vector bundles $v \colon \xi \oplus \xi'
  \xrightarrow{\cong} F^* \mu$.  Consider the morphism
  $u'$ in $\SPHB^G(X_1)$ which is given by the composition
\begin{multline*}
  S^{F^*\mu  \oplus \underline{\bbR^k}} \xrightarrow{S^{v^{-1} \oplus \id }}
  S^{\xi \oplus \xi'  \oplus \underline{\bbR^k}}
  \xrightarrow{\sigma_1} S^{\xi'} \wedge_{X_1} S^{\xi \oplus \underline{\bbR^k}}
  \xrightarrow{\id  \wedge_{X_1} u} S^{\xi'} \wedge_{X_1} S^{\xi \oplus \underline{\bbR^{k+n}}}
  \\
  \xrightarrow{\sigma_2} S^{\xi \oplus \xi' \oplus \underline{\bbR^{k+n}}}
  \xrightarrow{S^{v \oplus \id}} S^{F^*\mu \oplus \underline{\bbR^{k+n}}}.
\end{multline*}
By definition it represents the same element in $\pi_G^n(X_1,X_0)$ as
$u$.  Hence we can assume without loss of generality for the
representative $u$ of $a$ that the bundle $\xi$ is of the form
$F^*\mu$ for some $G$-vector bundle $\mu$ over $X$. Since
the morphism $u$ is trivial over $X_0$, we can find a morphism
$\SPHB^G(X)$
$$\overline{u} \colon S^{\mu \oplus \underline{\bbR^k}} \to S^{\mu \oplus\underline{\bbR^{k+n}}}$$
which is trivial over $X_2$ and satisfies
$F^*(\overline{u}) = u$.  Hence the morphism $\overline{u}$
defines an element in $\pi_G^n(X,X_1)$ such that
$\pi_G^n(F,f)([\overline{u}]) = [u] = a$ holds.  

It remains to prove injectivity of $\pi_G^n(F,f)$.
Consider an element $b$ in the kernel of $\pi_G^n(F,f)$.
Choose a morphism
$$u \colon S^{\xi \oplus \underline{\bbR^k}} \to S^{\xi \oplus \underline{\bbR^{k+n}}}$$
in $\SPHB^G(X)$ which is trivial over $X_2$ and represents $b$.
Then $F^*u \colon S^{F^*\xi \oplus \underline{\bbR^k}} \to S^{F^*\xi \oplus \underline{\bbR^{k+n}}}$
represents zero in $\pi_G^n(X_1,X_0)$. Hence we can find a bundle
$\mu$ over $X_1$ such that the composition
\begin{multline*} S^{\mu \oplus F^*\xi  \oplus \underline{\bbR^k}} \xrightarrow{\sigma_1}
S^{\mu} \wedge_{X_1} S^{F^*\xi \oplus \underline{\bbR^k}}  \xrightarrow{\id \wedge_{X_1} u}
S^{\mu} \wedge_{X_1}  S^{F^*\xi \oplus \underline{\bbR^{k+n}}}
\xrightarrow{\sigma_2}
S^{\mu \oplus  F^*\xi \oplus \underline{\bbR^{k+n}}}
\end{multline*}
is homotopic relative $X_0$ to the trivial morphism. As  in the proof of the surjectivity
we can arrange using Lemma~\ref{lem:extending_vector_bundles_in_LO1}
that $\mu$ is of the shape $F^*\xi'$ for some $G$-vector bundle $\xi'$ over $X$.
By replacing $u$ by $\id_{S^{\xi'}} \wedge_X u$ we can achieve that
$b = [u]$ still holds and additionally the morphism in $\SPHB^G(X_1)$
$$F^*u \colon S^{F^*(\xi \oplus  \underline{\bbR^k})} \to S^{F^*(\xi \oplus  \underline{\bbR^{k+n}})}$$
is homotopic relative $X_0$ to the trivial map. Since $u$ is trivial over
$X_2$, we can extend this homotopy trivially from $X_1$ to $X$ to show that $u$ itself is homotopic relative $X_2$
to the trivial map. But this means $b = [u] = 0$ in $\pi_G^n(X,X_2)$. This finishes the proof of
Lemma~\ref{lem:excision_for_pi_G^n}.
\end{proof}

\begin{lemma} \label{lem:two_stage_exact_sequence_of_a_triple}
Let $A \subseteq Y \subset X$ be inclusions of finite proper $G$-$CW$-complexes.
Then the sequence
$$\pi_G^n(X,Y) \xrightarrow{\pi_G^n(j)} \pi_G^n(X,A) \xrightarrow{\pi_G^n(i)}
\pi_G^n(Y,A)$$
is exact at $\pi_G^n(X,A)$ for all $n \in \bbZ$, where $i$ and $j$ denote
the obvious inclusions.
\end{lemma}
\begin{proof}
The inclusion $j \circ i \colon (Y,A) \to (X,Y)$ induces
the zero map $\pi_G^n(X,Y) \to \pi_n^G(Y,A)$ since an element
in $a \in \pi_G^n(X,Y)$ is represented by a morphism
$u \colon S^{\xi \oplus \underline{\bbR^k}} \to S^{\xi \oplus \underline{\bbR^{k+n}}}$
in $\SPHB^G(X)$ which is trivial over $Y$ and
$\pi_G^n(j \circ i)(a)$ is represented by the restriction of $u$ to $Y$.

Consider an element $a \in \pi_G^n(X,A)$ which is mapped to zero under
$\pi_G^n(i)$. Choose a morphism
$u \colon S^{\xi \oplus \underline{\bbR^k}} \to S^{\xi \oplus \underline{\bbR^{k+n}}}$
in $\SPHB^G(X)$ which is trivial over $A$ and represents $a$.
Hence we can find a $G$-vector bundle $\xi'$ over $Y$ such that
the morphism in $\SPHB^G(Y)$ given by the composition
$$S^{\xi' \oplus i^*\xi \oplus \underline{\bbR^k}}
\xrightarrow{\sigma_1}  S^{\xi'} \wedge_Y S^{i^*\xi \oplus \underline{\bbR^k}}
 \xrightarrow{\id\wedge_Y i^*u} S^{\xi'} \wedge_Y  S^{i^*\xi \oplus
\underline{\bbR^{k+n}}} \xrightarrow{\sigma_2}
S^{\xi' \oplus i^*\xi \oplus \underline{\bbR^{k+n}}}$$ is homotopic to
the trivial map relative $A$.

As  in the proof of Lemma~\ref{lem:excision_for_pi_G^n} we can arrange using
Lemma~\ref{lem:extending_vector_bundles_in_LO1} that $\xi'$ is the of the shape
$i^*\mu$ for some $G$-vector bundle $\mu$ over $X$. Hence we can
achieve by replacing $u$ by $\id_{S^{\mu}} \wedge_X u$ that $a = [u]$
still holds in $ \pi_G^n(X,A)$ and additionally the morphisms in $\SPHB^G(Y)$
$$i^*u \colon S^{j^*\xi \oplus \underline{\bbR^k}} \to S^{i^*\xi \oplus \underline{\bbR^{k+n}}}$$
is homotopic relative $A$ to the trivial map. One proves inductively
over the  number of equivariant cells in $X-Y$ and \cite[Theorem~1.2]{Lueck-Oliver(2001a)} that this
homotopy can be extended to a homotopy relative $A$ of the morphism $u$ to another morphism
$v \colon S^{\xi} \to S^{\xi \oplus \underline{\bbR^n}}$
in $\SPHB^G(X)$ which is trivial over $Y$.
But this implies that the element $[v] \in \pi_G^n(X,Y)$ represented by $v$ is mapped to $a = [u]$ under $\pi_G^n(i)$.
\end{proof}

In order to define a $G$-cohomology theory we must construct a
connecting homomorphism for pairs.
In the sequel maps denoted by $i_k$ are the obvious inclusions and
maps denoted by $\pr_k$ are the obvious projections.
Consider a pair of finite proper $G$-$CW$-complexes $(X,A)$ and $n \in \bbZ$, $n \ge 0$. We
define the homomorphism of abelian groups
\begin{eqnarray}
\delta_G^n(X,A) \colon \pi_G^n(A) & \to & \pi_G^{n+1}(X,A)
\label{delta_G^n(X,A): pi_G^n(A) to  pi_G^{n+1}(X,A)}
\end{eqnarray}
to be the composition
\begin{multline*}
\pi_G^n(A) \xrightarrow{\sigma^n_G(A)} \pi^{n+1}_G(A \times [0,1], A \times \{0,1\})
\\
\xrightarrow{\left(\pi_G^{n+1}(i_1)\right)^{-1}} \pi^{n+1}_G(X \cup_{A\times \{0\}} A \times [0,1], X \amalg  A \times \{1\})
\\
\xrightarrow{\pi_G^{n+1}(i_2)} \pi^{n+1}_G(X \cup_{A\times \{0\}} A \times [0,1], A \times \{1\})
\\
\xrightarrow{\left(\pi_G^{n+1}(\pr_1)\right)^{-1}} \pi^{n+1}_G(X,A),
\end{multline*}
where $\sigma^n_G(A)$ is the suspension isomorphism (see
Lemma~\ref{sigma^n_G(X,A)_is_a_bijection}), the map
$\pi_G^{n+1}(i_1)$ is bijective by excision
(see Lemma~\ref{lem:excision_for_pi_G^n}) and $\pi_G^{n+1}(\pr_1)$ is bijective by
homotopy invariance (see Lemma~\ref{lem:homotopy_invariance_of_pi^n_G}) since $\pr_1$ is a $G$-homotopy equivalence of pairs.

\begin{lemma} \label{lem:long_exact_sequence_of_a_pair}
  Let $(X,A)$ be a pair for proper finite $G$-$CW$-complexes.  Let $i
  \colon A \to X$ and $j \colon X \to (X,A)$ be the inclusions. Then
  the following long sequence is exact and natural in $(X,A)$:
\begin{multline*}
  \cdots \xrightarrow{\delta_G^{n-1}} \pi_G^n(X,A)
  \xrightarrow{\pi_G^n(j)} \pi_G^n(X) \xrightarrow{\pi_G^n(i)}
  \pi_G^n(A)
  \\
  \xrightarrow{\delta_G^n} \pi_G^{n+1}(X,A)
  \xrightarrow{\pi_G^{n+1}(j)} \pi_G^{n+1}(X)
  \xrightarrow{\pi_G^{n+1}(i)} \pi_G^{n+1}(A)
  \xrightarrow{\delta_G^{n-1}} \cdots.
\end{multline*}
\end{lemma}
\begin{proof}
  It is obviously natural. It remains to prove exactness.

  Exactness at $\pi^n_G(X)$ follows from Lemma~\ref{lem:two_stage_exact_sequence_of_a_triple}.

  Exactness at $\pi^n_G(X,A)$ follows from the following commutative
  diagram
  $$\begin{CD} \pi_G^{n+1}\left(X \cup_{A \times \{0\}} A \times
      [0,1], X \amalg A \times \{1\}\right) @> (\sigma_G^n(A))^{-1}
    \circ \pi_G^{n+1}(i_1) > \cong > \pi_G^n(A)
    \\
    @V \pi_G^{n+1}(i_2) VV @V \delta_G^{n+1}(X,A) VV
    \\
    \pi_G^{n+1}\left(X \cup_{A \times \{0\}} A \times [0,1], A \times
      \{1\}\right) @> \left(\pi_G^{n+1}(\pr_1)\right)^{-1} > \cong >
    \pi_G^{n+1}(X,A)
    \\
    @V \pi_G^{n+1}(i_3) VV @V \pi_G^{n+1}(i) VV
    \\
    \pi_G^{n+1}\left(X \amalg A \times \{1\}, A \times \{1\}\right) @>
    \pi_G^{n+1}(i_4) > \cong > \pi_G^{n+1}(X)
\end{CD}$$
whose left column is exact at $\pi_G^{n+1}\left(X \cup_{A \times
    \{0\}} A \times [0,1], A \times \{1\}\right) $ by
    Lemma~\ref{lem:two_stage_exact_sequence_of_a_triple}.

Exactness at $\pi^n_G(A)$ is proved analogously by applying
Lemma~\ref{lem:two_stage_exact_sequence_of_a_triple} to the inclusions
\begin{multline*}
  \left(X \times \{0,1\}, X \times \{0\} \amalg A \times \{1\}\right)
  \\
  \subseteq \left(X \times [0,1], X \times \{0\} \amalg A \times
    \{1\}\right)
  \\
  \subseteq \left(X \times [0,1], X \times \{0,1\}\right).
\end{multline*}
\end{proof}

We conclude from Lemma~\ref{lem:homotopy_invariance_of_pi^n_G},
Lemma~\ref{lem:excision_for_pi_G^n} and Lemma~\ref{lem:long_exact_sequence_of_a_pair}
that $\pi_G^*$ defines a $G$-cohomology theory
on the category of finite proper $G$-$CW$-complexes.

Consider a finite proper $G$-$CW$-complex $X$ with two subcomplexes $A,B \subseteq X$.
We want to define a multiplicative structure, i.e., a cup-product
\begin{eqnarray}
\cup \colon \pi_G^m(X,A) \times \pi_G^n(X,B) & \to & \pi_G^{m+n}(X,A \cup B).
\label{cup product on pi^*_G}
\end{eqnarray}
Given elements $a \in \pi_G^m(X,A)$ and $b \in \pi_G^n(X,B)$, choose appropriate morphisms
$u \colon S^{\xi \oplus \underline{\bbR^k}} \to S^{\xi \oplus \underline{\bbR^{m+k}}}$
and
$v \colon S^{\eta \oplus \underline{\bbR^l}} \to S^{\eta \oplus \underline{\bbR^{n+l}}}$
in $\SPHB^G(X)$ representing $a$ and $b$. Define $a \cup b$ by the class of the composition of
 morphisms in $\SPHB^G(X)$ which are trivial over $A \cup B$.
\begin{multline*}
S^{\xi \oplus  \eta \oplus \underline{\bbR^{k+l}}} \xrightarrow{\sigma_1}
S^{\xi \oplus \underline{\bbR^k}} \wedge_X S^{\eta \oplus \underline{\bbR^l}}
\\
\xrightarrow{u \wedge_X v}
S^{\xi \oplus \underline{\bbR^{k+m}}} \wedge_X S^{\eta \oplus \underline{\bbR^{l+n}}}
\xrightarrow{\sigma_2}
S^{\xi \oplus  \eta \oplus \underline{\bbR^{k+l+m+n}}}.
\end{multline*}

Next we deal with the induction structure.
Consider a group homomorphism $\alpha \colon H \to G$.
The pullback construction
for the $\alpha \colon H \to G$-equivariant map $X \to \ind_{\alpha} X = G \times_{\alpha} X, \; x \mapsto 1G \times_{\alpha} x$
defines a functor $\SPHB^G(\ind_{\alpha} X) \to \SPHB^H(X)$ which yields a homomorphism
of abelian groups
\begin{eqnarray}
\ind_{\alpha} \colon \pi^n_G\left(\ind_{\alpha} (X,A)\right) & \to & \pi_H^n(X).
\label{induction homorphism for pi^?_*}
\end{eqnarray}
Now suppose that the kernel of $H$ acts trivially on $X$.
Let $\xi \colon E \to X$ be a $H$-vector bundle.
Then $G \times_{\alpha} X$ is a proper $H$-$CW$-complex. Since $H$ acts freely on $X$, we obtain a
well-defined $G$-vector bundle
$G \times_{\alpha} \xi \colon G \times_{\alpha} E \to G \times_{\alpha} X$.
Thus we obtain a functor $\ind_{\alpha} \colon \SPHB^H(X) \to \SPHB^G(G \times_{\alpha} X)$.
It defines a homomorphism of abelian groups
\begin{eqnarray}
\ind_{\alpha} \colon \pi_H^n(X) & \to &
\pi^n_G\left(\ind_{\alpha} (X,A)\right).
\label{induction isomorphism for pi^?_*}
\end{eqnarray}
which turns out to be an inverse of the induction homomorphism~\eqref{induction homorphism for pi^?_*}.

Now we have all ingredients of an equivariant cohomology theory with a multiplicative structure.
We leave it to the reader to verify all the axioms. This finishes the proof of
Theorem~\ref{the:Equivariant_Stable_Cohomotopy_in_terms_of_equivariant_vector_bundles}.


\subsection{Comparison with the Classical Construction for Finite Groups}
\label{subsec:Comparison_with_the_Classical_Construction_for_Finite_Groups}

Next we want to show that for a finite group $G$ our construction reduces to the classical one.
We first explain why the finite group case is easier.

\begin{remark}[Advantages in the case of finite groups]
\label{rem:Advantages_in_the_case_of_finite_groups}
\em
The finite group case is easier because for finite groups the following facts are true.
The first fact is that every $G$-$CW$-complex is proper.
Hence one can view pointed $G$-$CW$-complexes, where the base point is fixed under the $G$-action and one
can carry out constructions like mapping cones without loosing the property proper. We need proper
to ensure that certain basic facts about $G$-vector bundles are true (see \cite[Section~1]{Lueck-Oliver(2001a)}).
The second fact is that every $G$-vector bundle over a finite $G$-$CW$-complex $\xi$ is a direct summand in
a trivial $G$-vector bundle, i.e., a $G$-vector bundle given by the projection $V \times X \to X$ for some $G$-representation
$V$. This makes for instance Lemma~\ref{lem:extending_vector_bundles_in_LO1} superfluous whose proof is
non-trivial in the infinite group case (see \cite[Lemma~3.7]{Lueck-Oliver(2001a)}).
\em
\end{remark}

Next we identify $\pi^n_G(X)$ defined in Subsection~\ref{subsec:The_Definition_of_Equivariant_Stable_Cohomotopy_Groups}
with the classical definitions which we have explained in
Subsection~\ref{subsec:The_Burnside_Ring_and_Equivariant_Stable_Cohomotopy} provided that $G$ is a finite group.

Consider an element  $a \in \pi_G^n(X)$  with respect to the definition given in
Subsection~\ref{subsec:The_Burnside_Ring_and_Equivariant_Stable_Cohomotopy}.
Obviously we can find a positive integer
$k \in \bbZ$ with $k +n \ge 0$ such that $a$ is represented for some complex $G$-representation $V$ by a $G$-map
$f \colon S^V \wedge S^k \wedge X_+ \to S^V \wedge S^{k+n}$. Let $\xi$ be the trivial $G$-vector bundle
$V \times X \to X$. Define a map
\begin{multline*}
\overline{f} \colon S^{\xi \oplus \underline{\bbR^k}} \xrightarrow{\sigma_1} (S^V \wedge S^k) \times X
\xrightarrow{f \times \pr_X} (S^V \wedge S^{k+n}) \times X \xrightarrow{\sigma_2}
S^{\xi \oplus \underline{\bbR^{k+n}}}.
\end{multline*}
It is a morphism in $\SPHB^G(X)$ and hence defines an element $a' \in \pi_G^n(X)$
with respect to the definition of Subsection~~\ref{subsec:The_Definition_of_Equivariant_Stable_Cohomotopy_Groups}.
Thus we get a homomorphism of abelian groups $a \mapsto a'$ from the definition of
Subsection~\ref{subsec:The_Burnside_Ring_and_Equivariant_Stable_Cohomotopy} to the one of
Subsection~\ref{subsec:The_Definition_of_Equivariant_Stable_Cohomotopy_Groups}.

Consider a morphism
$u \colon S^{\xi \oplus \underline{\bbR^{k}}} ~ \to ~ S^{\xi \oplus \underline{\bbR^{k + n}}}$
in $\SPHB^G(X)$ representing an element in $b = [u]$ in $\pi_G^n(X)$ as defined in
Subsection~\ref{subsec:The_Definition_of_Equivariant_Stable_Cohomotopy_Groups}.
Choose a $G$-vector bundle $\mu$, a complex $G$-representation $V$ and an isomorphism of (real) $G$-vector bundles
$\phi \colon \mu \oplus \xi \xrightarrow{\cong} V \times X$.
Then the morphism
\begin{multline*}
v \colon S^{V \oplus \bbR^k} \times X \xrightarrow{(\phi\oplus \id)^{-1}} S^{\mu \oplus \xi \oplus \underline{\bbR^k}}
\xrightarrow{\sigma_1} S^{\mu} \wedge_X S^{\xi \oplus \underline{\bbR^k}}
\\
\xrightarrow{\id \wedge u} S^{\mu} \wedge_X S^{\xi \oplus \underline{\bbR^{k+n}}}
\xrightarrow{\sigma_2} S^{\mu \oplus \xi \oplus \underline{\bbR^{k+n}}}
\xrightarrow{\phi \oplus \id} S^{V \oplus \bbR^{k+n}} \times X
\end{multline*}
is equivalent to $u$ and hence $b = [v]$. Since $v$ covers the identity and is fiberwise a pointed map,
its composition with the projection $S^{V \oplus \bbR^{k+n}} \times X \to S^{V \oplus \bbR^{k+n}} $ yields a map
$$\overline{v} \colon S^{V \oplus \bbR^k} \wedge X_+ \to S^{V \oplus \bbR^{k+n}}.$$
It defines an element $b' := [\overline{v}]$ in $\pi_G^n(X)$ with respect to the definition given in
Subsection~\ref{subsec:The_Burnside_Ring_and_Equivariant_Stable_Cohomotopy}. The map $b \mapsto b'$ is the inverse of the map
$a \to a'$ before.

\begin{remark} [Why consider $G$-vector bundles?]
\label{rem:Why_consider_G-vector_bundles} \em One may ask why we
consider $G$-vector bundles $\xi$ in
Section~\ref{sec:Equivariant_Stable_Cohomotopy_in_Terms_of_Real_Vector_Bundles}.
It would be much easier if we would only consider trivial
$G$-vector bundles $V \times X$ for $G$-representations $V$. Then
we would not need Lemma~\ref{lem:extending_vector_bundles_in_LO1}.
The proof that $\pi_G^*$ is a $G$-cohomology theory with a
multiplicative structure would go through and for finite groups we
would get the classical notion. The problem is that the induction
structure does not exists anymore as the following example shows.

Consider a finitely generated group such that $G_{\mrf}$ is trivial. Then any
$G$-representation $V$ is trivial (see
Subsection~\ref{subsec:Maximal_Residually_Finite_Quotients}). This implies that
a morphism $u \colon S^{\xi \oplus \underline{\bbR^k}} \to S^{\xi \oplus \underline{\bbR^{k+n}}}$
in $\SPHB^G(X)$ for $\xi = V \times X$ is the same as a (non-equivariant) map
$$S^{V \oplus \bbR^k} \wedge (G\backslash X_+) \to S^{V \oplus \bbR^{k+n}}.$$
This yields an identification of $\pi^n_G(X)$ with respect to the definition,
where all $G$-vector bundles are of the shape $X \times V$, with the (non-equivariant) stable
cohomotopy group $\pi^n_s(G\backslash X)$. If $G$ contains a non-trivial subgroup
$H \subseteq G$, then the existence of an induction structure would predict
for $X = G/H$ that $\pi^n_G(G/H)$ is isomorphic to $\pi^n_H$, which is in general
different from $\pi^n_s(G\backslash(G/H)) = \pi^n_{\{1\}}$.

So we need to consider $G$-vector bundles in order to get induction structures
and hence an equivariant cohomology theory. In particular
our definition guarantees $\pi^n_G(G/H) = \pi^n_H$ for
every group $G$ with a finite subgroup $H \subseteq G$.
\em
\end{remark}

\begin{remark}[The coefficients of equivariant stable cohomotopy]
\label{exa:The_coefficients_of_equivariant_stable_cohomotopy}
\em
It is important to have information about the values $\pi^n_G(G/H)$ for a finite subgroup $H \subseteq G$ of a group $G$.
By the induction structure and the identification above $\pi^n_G(G/H)$ agrees with
the abelian groups $\pi_H^n = \pi^H_{-n}$  defined in Subsection~\ref{subsec:The_Burnside_Ring_and_Equivariant_Stable_Cohomotopy}.
The equivariant homotopy groups $\pi^H_{-n}$ are computed in terms of the splitting due to Segal and tom Dieck
(see \cite[Theorem~7.7 in Chapter~II on page~154]{Dieck(1987)},
\cite[Proposition~2]{Segal(1971)}) by
$$\pi^n_G(G/H) ~ = ~ \pi^H_{-n} ~ = ~ \bigoplus_{(K) \in \consub(H)} \pi^s_{-n}(BW_HK).$$
The abelian group $\pi_q^s$ is finite for $q \ge 1$ by a result of
Serre~\cite{Serre(1953)} (see also \cite{Klaus-Kreck(2004)}), is $\bbZ$ for $q = 0$
and is trivial for $q \le -1$. Since $W_HK$ is finite, $\widetilde{H}_p(BW_HK;\bbZ)$ is finite for all
$p \in \bbZ$. We conclude from the Atiyah-Hirzebruch spectral sequence that
$\pi^s_{-n}(BW_HK)$ is finite for $n \le -1$. This implies
$|\pi^n_G(G/H)| <  \infty$ for $n \le -1$
and that $\pi^n_G(G/H) = 0$ for $n \ge 1$. We know already $\pi^0_H = A(H)$ from
Theorem~\ref{the:pi^0_G_cong_A(G)}. Thus we get
$$\begin{array}{lllll}
|\pi^n_G(G/H)| & < & \infty & & n \le -1;
\\
\pi^0_G(G/H) & = & A(H); & &
\\
\pi^n_G(G/H) & = & \{0\} & & n \ge 1.
\end{array}
$$\em
\end{remark}

\begin{remark}{\bf (Equivariant stable Cohomotopy for arbitrary $G$-$CW$-com\-plexes).}
\label{rem:Equivariant_Stable_Cohomotopy_for_arbitrary_G-CW-complexes}
\em In order to construct an equivariant cohomology theory or an
(equivariant homology theory) for arbitrary $G$-$CW$-complexes it
suffices to construct a contravariant (covariant) functor from the
category of small groupoids to the category of spectra (see
\cite{SauerJ(2002)}, \cite[Proposition~6.8]{Lueck-Reich(2004)}).
In a different paper  we will carry out such a construction
yielding equivariant cohomotopy and homotopy for arbitrary
equivariant $CW$-complexes and will identify the result with the
one presented here for finite proper $G$-$CW$-complexes. \em
\end{remark}


\subsection{Rational Computation of Equivariant Cohomotopy}
\label{subsec:Rational_Computation_of_Equivariant_Stable_Cohomotopy}

The cohomotopy theoretic Hurewicz homomorphism yields a transformation of cohomology theories
$$\pi^*_s(X) \to H^*(X;\bbZ)$$
from the (non-equivariant) stable cohomotopy
to singular cohomology with $\bbZ$-coefficients. It is rationally an isomorphism provided that
$X$ is a finite $CW$-complex. It is compatible with the multiplicative structures.
The analogue for equivariant cohomotopy is described next.

Let $G$ be a group and $H \subseteq G$ be a finite subgroup.
Consider a pair of finite proper $G$-$CW$-complexes $(X,A)$.
Lemma~\ref{lem:orbits_and_Weyl_groups} implies that $(X^H,A^H)$ is a pair of  finite proper $W_GH$-$CW$-complexes
and $W_GH\backslash (X^H,A^H)$ is a pair of finite $CW$-complexes.
Taking the $H$-fixed point set yields a homomorphism
$$\alpha^n_{(H)}(X,A) \colon \pi_G^n(X,A) ~ \to ~ \pi_{W_GH}^n(X^H,A^H).$$
This map is natural and compatible with long exact sequences of pairs and Mayer-Vietoris sequences.

The induction structure with respect to the homomorphism $W_GH \to \{1\}$ yields a homomorphism
$$\beta_{W_GH}^n(X^H,A^H) \colon  \pi^n_s\left(W_GH\backslash (X^H,A^H)\right) \to \pi_{W_GH}^n(X^H,A^H).$$
We claim that $\beta_{W_GH}^n(Z,B)$ is a rational isomorphism for any pair of finite proper $W_GH$-$CW$-complexes
$(Z,B)$. Since $\beta_{W_GH}^*$ is natural and compatible with the long exact sequences of pairs and
Mayer-Vietoris sequences, it suffices to prove the claim for $Z = W_GH/L$ and $B = \emptyset$ for any finite subgroup
$L \subset W_GH$. But then $\beta_{W_GH}^n$ reduces to the obvious map
$\pi^n(\pt) \to \pi^n_L$ which is a rational isomorphism by
Remark~\ref{exa:The_coefficients_of_equivariant_stable_cohomotopy}.

Let
$$h^n(X^H,A^H) \colon \pi^n_s\left(W_GH\backslash (X^H,A^H)\right) ~ \to ~ H^n\left(W_GH\backslash (X^H,A^H);\bbZ\right)$$
be the cohomotopy theoretic Hurewicz homomorphism. Let
\begin{multline*}
\gamma^n\left(W_GH\backslash(X^H,A^H)\right) \colon H^n\left(W_GH\backslash (X^H,A^H);\bbZ\right) \otimes_{\bbZ} \bbQ
\\
\to H^n\left(W_GH\backslash (X^H,A^H);\bbQ\right)
\end{multline*}
be the natural map. Define a  $\bbQ$-homomorphism by the composition
\begin{multline*}
\zeta^n_G(X,A)_{(H)} \colon  \pi_G^n(X,A) \otimes_{\bbZ} \bbQ
\xrightarrow{\alpha^n_{(H)}(X,A) \otimes_{\bbZ} \id_{\bbQ}} \pi_{W_GH}^n(X^H,A^H)\otimes_{\bbZ} \bbQ
\\
\xrightarrow{\left(\beta_{W_GH}^n(X^H,A^H) \otimes_{\bbZ} \id_{\bbQ}\right)^{-1}}
 \pi^n_s\left(W_GH\backslash (X^H,A^H)\right) \otimes_{\bbZ} \bbQ
\\
\xrightarrow{h^n(X^H,A^H) \otimes_{\bbZ} \id_{\bbQ}} H^n\left(W_GH\backslash (X^H,A^H);\bbZ\right) \otimes_{\bbZ} \bbQ
\\
\xrightarrow{\gamma^n\left(W_GH\backslash(X^H,A^H)\right)} H^n\left(W_GH\backslash (X^H,A^H);\bbQ\right).
\end{multline*}
Define
\begin{multline}
\zeta^n_G(X,A) = \prod_{(H) \in \consub(G)} \zeta^n_G(X,A)_{(H)}  \colon  \pi_G^n(X,A) \otimes_{\bbZ} \bbQ
\\
~ \to ~ \prod_{(H) \in \consub(G)}  H^n\left(W_GH\backslash (X^H,A^H);\bbQ\right).
\label{def_of_seta^g_n(X,A)}
\end{multline}

\begin{theorem} The maps
$$\zeta^n_G(X,A)\colon  \pi_G^n(X,A) \otimes_{\bbZ} \bbQ ~ \xrightarrow{\cong} ~
\prod_{(H) \in \consub(G)}  H^n\left(W_GH\backslash (X^H,A^H);\bbQ\right)$$
are bijective for all $n \in \bbZ$ and all pairs of finite proper $G$-$CW$-complexes $(X,A)$.
They are compatible with the obvious multiplicative structures.
\end{theorem}
\begin{proof}
One easily checks that $\zeta^*_G$ defines a transformation of $G$-homology theories, i.e., 
is natural in $(X,A)$ and compatible with long exact sequences of pairs and Mayer-Vietoris sequences.
Hence it suffices to show that $\zeta^n_G(G/K)$ is bijective for all $n \in \bbZ$ and finite subgroups
$K \subset G$. The source and target of $\zeta^n_G(G/K)$ are trivial for $n \not= 0$
(see Remark~\ref{exa:The_coefficients_of_equivariant_stable_cohomotopy}). The map
$\zeta^0_G(G/K)$ can be identified using
Lemma~\ref{lem:orbits_and_Weyl_groups}~\ref{lem:orbits_and_Weyl_groups:WKbackslashG/H^K}
with the rationalization of the character map
$$\character^H \colon A(H)  \to  \prod_{(H) \in \consub(G)} \bbZ$$
defined in \eqref{character_map_for_A(G)_for_finite_G}
which is bijective by Theorem~\ref{the:Burnside_ring_congruences_for_finite_groups}.
\end{proof}


\subsection{Relating Equivariant Stable Cohomotopy and Equivariant Topological $K$-Theory}
\label{subsec:Relating_Equivariant_Stable_Cohomotopy_and_Equivariant_Topological_K-Theory}

We have introduced two equivariant cohomology theories with multiplicative structure, namely equivariant cohomotopy
(see~Theorem~\ref{the:Equivariant_Stable_Cohomotopy_in_terms_of_equivariant_vector_bundles}) and equivariant topological $K$-theory
(see Subsection~\ref{subsec:Equivariant_Topological_K-Theory}).

Let $X$ be a finite proper $G$-$CW$-complex with $G$-$CW$-subcomplexes $A$ and $B$
and let $a \in \pi^m_G(X,A)$ be an element.
We want to assign to it for every $m \in \bbZ$ a homomorphism of abelian groups
\begin{eqnarray}
\phi^{m,n}_G(X,A)(a) \colon K_G^{n}(X,B) \to K_G^{m+n}(X,A \cup B).
\label{phi^{m,n}_G(X,AcupB)(a)}
\end{eqnarray}
Choose an integer $k \in \bbZ$ with $k \ge 0, k + m \ge 0$ and a  morphism
$u \colon S^{\xi \oplus \underline{\bbR^k}} \to S^{\xi \oplus \underline{\bbR^{k+m}}}$ in
$\SPHB^G(X)$ which is trivial over $A$ and represents $a$.
Let $v$ be the morphism in
$\SPHB^G(X)$ which is given by the composite
$$v \colon S^{\xi \oplus \xi \oplus \underline{\bbR^k}} \xrightarrow{\sigma} S^{\xi} \wedge_X S^{\xi \oplus \underline{\bbR^k}}
\xrightarrow{\id \wedge_X u}
S^{\xi} \wedge_X S^{\xi \oplus \underline{\bbR^{k+m}}} \xrightarrow{\sigma^{-1}} S^{\xi \oplus \xi \oplus \underline{\bbR^{k+m}}}.$$
Then $v$ is another representative of $a$. The bundle
$\xi \oplus \xi$ carries a canonical structure of a complex vector bundle and we denote this
complex vector bundle by $\xi_{\bbC}$.

Let $\sigma^k(X,A \cup B) \colon  K_G^{m+n}(X,A \cup B)
\xrightarrow{\cong} K_G^{m+n+k}\left((X,A \cup B) \times (D^k,S^{k-1})\right)$ be the suspension isomorphism.
Let $\pr_k \colon X \times D^k \to X$ be the projection and $\pr_k^*\xi_{\bbC}$ be the complex vector bundle
obtained from $\xi_{\bbC}$ by the pull back construction. Associated to it is a Thom isomorphism
\begin{multline*}
T_{\pr_k^*\xi_{\bbC}}^{m+n+k} \colon K_G^{m+n+k}\left((X,A \cup B) \times (D^k,S^{k-1})\right)
\\
\xrightarrow{\cong}
K_G^{m+n + k+ 2 \cdot \dim(\xi)}\left(S^{\pr_k^*\xi_{\bbC}},
S^{\pr_k^*\xi_{\bbC}|_{X \times S^{k-1} \cup (A \cup B) \times D^k}} \cup (X \times D^k)_{\infty}\right),
\end{multline*}
where $(X \times D^k)_{\infty}$ is the copy of $X \times D^k$ given by the various points at infinity in the fibers
$S^{\pr_k^*\xi_{\bbC}}$ and $\pr_k^*\xi_{\bbC}|_{X \times S^{k-1} \cup (A \cup B) \times D^k}$ is the restriction of
$\pr_k^*\xi_{\bbC}$ to $X \times S^{k-1} \cup (A \cup B) \times D^k$
(see \cite[Theorem~3.14]{Lueck-Oliver(2001a)}). Let
\begin{multline*}
p_k \colon \left(S^{\pr_k^*\xi_{\bbC}},S^{\pr_k^*\xi_{\bbC}|_{X \times S^{k-1} \cup (A \cup B) \times D^k}}
\cup (X \times D^k)_{\infty}\right)
\\
\left(S^{\xi \oplus \xi \oplus \underline{\bbR^k}}, S^{\xi \oplus \xi \oplus \underline{\bbR^k}|_{A \cup B}} \cup X_{\infty}\right)
\end{multline*}
be the obvious projection which induces by excision an isomorphism on $K_G^*$. Define an isomorphism
\begin{multline*}
\mu^{m+n,m+n+k+2\cdot\dim(\xi)} \colon K_G^{m+n}(X,A\cup B)
\\
 \to
K_G^{m+n+k+2 \cdot \dim(\xi)}
\left(S^{\xi \oplus \xi \oplus \underline{\bbR^k}}, S^{\xi \oplus \xi \oplus \underline{\bbR^k}|_{A \cup B}} \cup X_{\infty}\right)
\end{multline*}
by the composite $K_G^{m+n+k+2 \cdot \dim(\xi)}(p_k)^{-1} \circ T_{\pr_k^*\xi_{\bbC}}^{n+m+k} \circ  \sigma^k(X,A \cup B)$.
Define
\begin{multline*}
\mu^{n,m+n+k+2\cdot \dim(\xi)} \colon K_G^{n}(X,B)
\\
\to
K_G^{m+n+k+2 \cdot \dim(\xi)}
\left(S^{\xi \oplus \xi \oplus \underline{\bbR^{k+m}}}, S^{\xi \oplus \xi \oplus \underline{\bbR^{k+m}}|_{B}} \cup X_{\infty}\right)
\end{multline*}
analogously. Let the desired map $\phi^{m,n}_G(X,A)(a)$ be the composite
\begin{multline*}
\phi^{m,n}_G(X;A,B)(a) \colon K_G^{n}(X,B)
\\
\xrightarrow{\mu^{n,m+n+k+2\cdot \dim(\xi)}}
K_G^{m+n+k+2 \cdot \dim(\xi)}
\left(S^{\xi \oplus \xi \oplus \underline{\bbR^{k+m}}}, S^{\xi \oplus \xi \oplus \underline{\bbR^{k+m}}|_{B}} \cup X_{\infty}\right)
\\
\xrightarrow{K_G^{m+n+k+2 \cdot \dim(\xi)}(v)}
K_G^{m+n+k+2 \cdot \dim(\xi)}\left(S^{\xi \oplus \xi \oplus \underline{\bbR^k}}, S^{\xi \oplus \xi \oplus \underline{\bbR^k}|_{A \cup B}} \cup X_{\infty}\right)
\\
\xrightarrow{\left(\mu^{m+n,m+n+k+2\cdot \dim(\xi)}\right)^{-1}}
K_G^{m+n}(X,A \cup B).
\end{multline*}
We leave it to the reader to check that the definition of $\phi^{m,n}_G(X;A,B)(a)$
is independent of the choices of $k$ and $u$. The maps $\phi^{m,n}_G(X;A,B)(a)$ for the various
elements $a \in \pi^m_G(X,A)$ define pairings
\begin{eqnarray}
\phi^{m,n}_G(X;A,B) \colon \pi^m_G(X,A) \times K_G^{n}(X,B) & \to & K_G^{m+n}(X,A \cup B).
\label{phi^{m,n}_G(X,A)}
\end{eqnarray}

The verification of the next theorem is left to the reader.
\begin{theorem}{\bf (Equivariant topological $K$-theory as graded algebra over equivariant stable cohomotopy).}
\label{the:K^*_G_as_graded_module_over_pi^*_G}
\begin{enumerate}

\item \label{the:K^*_G_as_graded_module_over_pi^*_G:Naturality} Naturality\\[1mm]
The pairings $\phi^{m,n}_G(X;A.B)$ are natural in $(X;A,B)$;

\item \label{the:K^*_G_as_graded_module_over_pi^*_G:algebra_structure} Algebra structure\\[1mm]
The collection of the pairings $\phi^{m,n}_G(X;\emptyset,A)$ defines on $K^*_G(X,A)$ the structure of a graded algebra
over the graded ring $\pi^*_G(X)$ ;

\item \label{the:K^*_G_as_graded_module_over_pi^*_G:induction} Compatibility with induction\\[1mm]
Let $\phi \colon H \to G$ be a group homomorphism and $(X,A)$ be a
pair of proper finite $G$-$CW$-complexes. Then the following
diagram commutes
$$\begin{CD}
\pi^m_G\left(\ind_{\alpha}(X,A)\right) \times K_G^{n}\left(\ind_{\alpha}(X,B)\right)
@> \phi^{m,n}_G\left(\ind_{\alpha}(X;A,B)\right) >>
K_G^{m+n}\left(\ind_{\alpha}(X,A\cup B)\right)
\\
@V \ind_{\alpha} \times \ind_{\alpha} VV  @V \ind_{\alpha} VV
\\
\pi^m_H(X,A) \times K_H^{n}(X,B)
@> \phi^{m,n}_H(X;A,B) >>
K_H^{m+n}(X,A \cup B)
\end{CD}$$

\item \label{the:K^*_G_as_graded_module_over_pi^*_G:boundary_operator}
For $a \in \pi^{m-1}_G(A)$ and $b \in K^n_G(X)$ we have
$$\phi^{m,n}_G(X;\emptyset,\emptyset)(\delta(a),b) ~ = \delta\left(\phi^{m-1,n}_G(A;\emptyset,\emptyset)(a,K_G^n(j)(b))\right),$$
where $\delta \colon \pi^{m-1}(A) \to \pi^m_G(X)$ and $\delta
\colon K^{m+n-1}_G(A) \to K^{m+n}_G(X)$ are boundary operators for
the pair $(X,A)$ and $j \colon A \to X$  is the inclusion.
\end{enumerate}
\end{theorem}

For every pair  $(X,A)$ of finite proper $G$-$CW$-complexes define a homomorphism
\begin{eqnarray}
\psi^n_G(X,A) \colon \pi^n_G(X,A) \to  K_G^n(X,A), && a ~ \mapsto \phi^{n,0}_G(X,A,\emptyset)(a,1_X),
\label{psi^n_G(X,A)}
\end{eqnarray}
where $1_X \in K_G^0(X)$ is the unit element. Then Theorem~\ref{the:K^*_G_as_graded_module_over_pi^*_G}
implies

\begin{theorem}[Transformation from equivariant stable cohomotopy to equivariant topological $K$-theory]
\label{the:Transformation_from_equivariant_Cohomotopy_to_equivariant_topological_K-theory}
We obtain a natural transformation of equivariant
cohomology theories with multiplicative structure for pairs of equivariant proper finite $CW$-complexes
by the maps
$$\psi^*_? \colon \pi^*_? \to K^*_?.$$
If $H \subseteq G$ is a finite subgroup of the group $G$, then the map
$$\psi^n_G(G/H) \colon \pi^n_G(G/H) \to K^0_G(G/H)$$
is trivial for $n \ge 1$ and agrees for $n = 0$ under the identifications
$\pi^0_G(G/H) = \pi^0_H = A(H)$ and $K^0_G(G/H) = K^0_H(\pt) = R_{\bbC}(H)$ with the ring homomorphism
$$A(H) \to R_{\bbC}(H), \quad [S] \mapsto [\bbC[S]]$$
which assigns to a finite $H$-set the associated complex permutation representation.
\end{theorem}


\typeout{--------- Section 7: The Homotopy theoretic Burnside ring   -----------}

\section{The Homotopy Theoretic Burnside Ring}
\label{sec:The_Homotopy_Theoretic_Burnside_ring}

In this section we introduce another version of the Burnside ring which is of homotopy theoretic nature
and probably the most sophisticated and interesting one.


\subsection{Classifying Space for Proper $G$-Actions}
\label{subsec:Classifying_space_for_proper_G-actions}

We need the following notion due to tom Dieck~\cite{Dieck(1972)}.

\begin{definition}[Classifying space for proper $G$-actions]
\label{def:Classifying_space_for_proper_G-action} A model
for the \emph{classifying space for proper $G$-actions} is a proper $G$-$CW$-complex
$\underline{E}G$ such that $\underline{E}G^H$ is contractible for every finite subgroup
$H \subseteq G$.
\end{definition}

Recall that a $G$-$CW$-complex is proper if and only if all its isotropy groups are finite.
If $\underline{E}G$ is a model for the classifying space for proper $G$-actions, then for every
proper $G$-$CW$-complex $X$ there is up to $G$-homotopy precisely one $G$-map $X \to \underline{E}G$.
In particular two models are $G$-homotopy equivalent and the $G$-homotopy equivalence between two models
is unique up to $G$-homotopy. If $G$ is finite, a model for $\underline{E}G$ is $G/G$.
If $G$ is torsionfree, $\underline{E}G$ is the same as $EG$ which is by definition the total space
of the universal principal $G$-bundle $G \to EG \to BG$.

Here is a list of groups $G$ together with specific models
for $\underline{E}G$ with the property that the model is a finite $G$-$CW$-complex.
\begin{center}
\begin{tabular}{|p{75mm}|p{38mm}|}
\hline\hline
$G$ & $\underline{E}G$
\\
\hline\hline
word hyperbolic groups
&
Rips complex
\\
\hline
discrete cocompact subgroup $G \subseteq L$ of a Lie group $L$
with finite $\pi_0(L)$
&
$L/K$ for a maximal compact subgroup $K \subseteq L$
\\
\hline
$G$ acts by isometries properly and cocompactly on a CAT(0)-space $X$,
for instance on a tree or a simply-connected complete Riemannian manifold
with non-positive sectional curvature
&
$X$
\\
\hline
arithmetic groups
&
Borel-Serre completion
\\
\hline
mapping class groups
&
Teichm\"uller space
\\
\hline
outer automorphisms of finitely generated free groups
&
outer space
\\
\hline\hline
\end{tabular}
\end{center}

 More information and more references about $\underline{E}G$ can be found for instance in
\cite{Baum-Connes-Higson(1994)} and \cite{Lueck(2004a)}.


\subsection{The Definition of the Homotopy Theoretic Burnside Ring}
\label{subsec:The_Definition_of_the_Homotopy_Theoretic_Burnside_Ring}

We have introduced the equivariant cohomology theory with
multiplicative structure for proper finite equivariant
$CW$-complexes $\pi^*_?$ in
Section~\ref{sec:Equivariant_Stable_Cohomotopy_in_Terms_of_Real_Vector_Bundles}
and the classifying space $\underline{E}G$ for proper $G$-actions
in Subsection~\ref{def:Classifying_space_for_proper_G-action}.

\begin{definition}[Homotopy theoretic Burnside ring]
\label{def:Homotopy theoretic Burnside ring}
Let $G$ be a (discrete) group such that there exists a finite model $\underline{E}G$ for
the universal space for proper $G$-actions. Define the
\emph{homotopy theoretic Burnside ring} to be
$$\hoA(G) ~ :=  ~ \pi^0_G(\underline{E}G).$$
\end{definition}

If $G$ is finite, $\pi_G^0(\underline{E}G)$ agrees with $\pi_G^0$ which is isomorphic to the Burnside ring
$A(G)$ by Theorem~\ref{the:pi^0_G_cong_A(G)}. So the homotopy theoretic
definition $\hoA(G)$ reflects this aspect of the Burnside ring
which has not been addressed by the other definitions before.

After the program described in
Remark~\ref{rem:Equivariant_Stable_Cohomotopy_for_arbitrary_G-CW-complexes}
has been carried out, the assumption in
Definition~\ref{def:Homotopy theoretic Burnside ring} that there
exists a finite model for $\underline{E}G$ can be dropped and thus
the Homotopy Theoretic Burnside ring $\hoA(G)$ can be defined by
$\pi^0_G(\underline{E}G)$ and analyzed for all discrete groups
$G$.

If $G$ is torsionfree, $\hoA(G)$ agrees with $\pi^0_s(BG)$.

Theorem~\ref{the:Transformation_from_equivariant_Cohomotopy_to_equivariant_topological_K-theory} implies
that the map (see \eqref{psi^n_G(X,A)})
$$\psi^0_G(\underline{E}G) \colon \hoA(G) = \pi^0_G(\underline{E}G) ~ \to ~  K^0_G(\underline{E}G)$$
is a ring homomorphism. It reduces for finite $G$ to the ring
homomorphism $A(G) \to R_{\bbC}(G)$ sending the class of a finite
$G$-set to the class of the associated complex permutation
representation.


\subsection{Relation between the Homotopy Theoretic and the Inverse-Limit-Version}
\label{subsec:Relation_between_hoA(G)_and_invA(G)}

Suppose there is a finite model for $\underline{E}G$. Then there is an equivariant Atiyah-Hirzebruch
spectral sequence which converges to $\pi^n_G(\underline{E}G)$ and whose $E^2$-term is given in terms of
\emph{Bredon cohomomology}
$$E_2^{p,q} = H^p_{\bbZ\SubGF{G}{\calfin}}(\underline{E}G;\pi^q_?).$$
Here $\pi^q_?$ is the contravariant functor
$$\pi^q_? \colon \SubGF{G}{\calfin} \to \bbZ-\MODULES, \quad H \mapsto \pi^q_H$$
and naturality comes from restriction with a group homomorphism $H \to K$
representing a morphism in $\SubGF{G}{\calfin}$. Usually the Bredon cohomology is defined over the orbit category,
but in our case we can pass to the category $\SubGF{G}{\calfin}$ because of
Lemma~\ref{lem:calh_G(G/H)_and_calh_H(*)}. Details of the
construction of $H^p_{\bbZ\SubGF{G}{\calfin}}(\underline{E}G;\pi_q^?)$ can be found for instance
in \cite[Section~3]{Lueck(2004i)}. We will only need the following elementary facts. There
is a canonical identification
\begin{eqnarray}
H^0_{\bbZ\SubGF{G}{\calfin}}(\underline{E}G;\pi^q_?) & \cong  & \invlim{H \in \SubGF{G}{\calfin}}{\pi^q_H}.
\label{zero-th_Bredon_as_invlim}
\end{eqnarray}
If we combine \eqref{zero-th_Bredon_as_invlim} with Theorem~\ref{the:pi^0_G_cong_A(G)}
we get an identification
\begin{eqnarray}
H^0_{\bbZ\SubGF{G}{\calfin}}(\underline{E}G;\pi^0_?) & \cong  & \invA(G).
\label{zero-th_Bredon_and_A_inv(G)}
\end{eqnarray}
The assumption that $\underline{E}G$ is finite implies together with
Remark~\ref{exa:The_coefficients_of_equivariant_stable_cohomotopy}
\begin{eqnarray}
|H^p_{\bbZ\SubGF{G}{\calfin}}(\underline{E}G;\pi^q_?)|  < \infty  &&\mbox{if } q \le -1;
\label{info_about_E^{p,q}_2_i}
\\
H^p_{\bbZ\SubGF{G}{\calfin}}(\underline{E}G;\pi^q_?)  =  \{0\} &&
\mbox{if } p > \dim(\underline{E}G) \mbox{ or } p \le -1 \mbox{ or } q \ge 1.
\label{info_about_E^{p,q}_2_ii}
\end{eqnarray}

The equivariant Atiyah-Hirzebruch spectral sequence together with~\eqref{zero-th_Bredon_and_A_inv(G)},
\eqref{info_about_E^{p,q}_2_i} and~\eqref{info_about_E^{p,q}_2_ii}
implies
\begin{theorem}[Rationally $\hoA(G)$ and $\invA(G)$ agree]
\label{the:Rationally_hoA(G)_and_invA(G)_agree}
Suppose that there is a finite model for $\underline{E}G$. Then the edge homomorphism
$$\edge^G \colon \hoA(G) = \pi_G^0(\underline{E}G) ~ \to ~ \invA(G)$$
is a ring homomorphism whose kernel and the cokernel are finite.
\end{theorem}

The edge homomorphism appearing in Theorem~\ref{the:Rationally_hoA(G)_and_invA(G)_agree} can be made explicit.
Consider a morphism $u \colon S^{\xi \oplus \underline{\bbR^k}} \to S^{\xi \oplus \underline{\bbR^k}}$
in $\SPHB^G(\underline{E}G)$ representing the element $a \in \pi^0_G(\underline{E}G)$.
In order to specify $\edge^G(a)$ we must define for every finite subgroup $H \subseteq G$ an element
$\edge^G(a)_H \in A(H)$. Choose a point $x \in \underline{E}G^H$. Then $u$ induces
a pointed $H$-map $S^{\xi_x \oplus \bbR^k} \to S^{\xi_x \oplus \bbR^k}$. It defines an element
in $\pi^0_H$. Let $\edge^G(a)_H$ be the image of this element under the ring isomorphism
$\deg^H \colon \pi_H^0 \xrightarrow{\cong} A(H)$ appearing in Theorem~\ref{the:pi^0_G_cong_A(G)}.
One easily checks that the collection of these elements
$\edge^G_H(a)$ does define an element in the inverse limit $\invA(G)$. So essentially $\edge^G$
is the map which remembers just the system of the maps of the various fibers.

\begin{remark}{\bf (Rank of the abelian group $\hoA(G)$).}
\label{rem:rank_of_the_abelian_group_hoA(G)} \em
A kind of character map for the homotopy theoretic version would be the composition
of $\edge$ and the character map $\invcharacter^G$ of~\eqref{character^G_invA(G)_to_prod_{(H),|H|<infty}Z}.
Since we assume that $\underline{E}G$ has a finite model, there are only finitely many conjugacy
classes of finite subgroups and the Burnside ring congruences appearing
in Theorem~\ref{the:Burnside_ring_congruences_for_invA(G)} becomes easier to handle.
In particular we conclude from Example~\ref{exa:rank_of_invA(G)_for_|ccs_f(G)|_finite}
and Theorem~\ref{the:Rationally_hoA(G)_and_invA(G)_agree} that $\hoA(G)$ is a finitely generated abelian group
whose rank is the number $|\consub_f(G)|$ of conjugacy classes of finite subgroups of $G$.
\em
\end{remark}


\subsection{Some Computations of the Homotopy Theoretic Burnside Ring}
\label{subsec:Some Computations_of_the_Homotopy_Theoretic_Burnside_Ring}

\begin{example}[Groups with appropriate maximal finite subgroups] \label{exa:A_ho_conditions_M_and_NM} \em
Suppose that the group $G$ satisfies the conditions appearing in Example~\ref{exa:conditions_M_and_NM}
and admits a finite model for $\underline{E}G$.
In the sequel we use the notation introduced in Example~\ref{exa:conditions_M_and_NM}.
Then one can construct a $G$-pushout (see \cite[Section 4.11]{Lueck(2004a)})
\begin{eqnarray}
& \comsquare{\coprod_{i \in I} G \times_{M_i} EM_i}{i}{EG}
{}{}{\coprod_{i \in I} G/M_i}{}{\underline{E}G} &
\label{pushout_for_underlineEG}
\end{eqnarray}
Taking the $G$-quotient, yields a non-equivariant pushout. There are long exact Mayer-Vietoris sequence
associated to \eqref{pushout_for_underlineEG} and to the $G$-quotient.
(We ignore the problem that $G \times_EM_i$ and $EG$ may not be finite. It does not really matter
since both are free or because we  will in a different paper extend the definition of equivariant stable cohomotopy
to all proper equivariant $CW$-complexes). These are linked
by the induction maps with respect to the projections $G \to \{1\}$.
Splicing these two long exact sequences together, yields the long exact sequence
\begin{multline*}
\cdots \to \prod_{i \in I} \ker\left(\res_{M_i}^{\{1\}} \colon \pi^{-1}_{M_i} \to \pi^{-1}_s\right) \to \pi^0_s(G\backslash \underline{E}G) \to \hoA(G)
\\
\to \prod_{i \in I} \widetilde{A}(M_i) 
\to  \pi^1_s(G\backslash \underline{E}G) \to \cdots
\end{multline*}
\em
\end{example}

\begin{example}[Extensions of $\bbZ^n$ with $\bbZ/p$ as quotient] \label{exa:hoA_extensions_by_Z/p} \em
Suppose that $G$ satisfies the assumptions appearing in Example~\ref{exa:extensions_by_Z/p}.
Then $G$ admits a finite model for $\underline{E}G$. In the sequel we use the notation introduced 
in Example~\ref{exa:extensions_by_Z/p}. Then variation of the argument above yields a long exact 
sequence
\begin{multline*}
\cdots \to \prod_{H^1(\bbZ/p;A)} \ker\left(\res_{\bbZ/p}^{\{1\}} \colon \pi^{-1}_{\bbZ/p}(BA^{\bbZ/p}) \to \pi^{-1}_{s}(BA^{\bbZ/p})\right)  \to 
\pi^0_s(G\backslash \underline{E}G) 
\\
\to \hoA(G)
\to \prod_{{H^1(\bbZ/p;A)}} \ker\left(\res_{\bbZ/p}^{\{1\}} \colon \pi^0_{\bbZ/p}(BA^{\bbZ/p}) \to \pi^0_{s}(BA^{\bbZ/p})\right) 
\to  \pi^1_s(G\backslash \underline{E}G) \to \cdots
\end{multline*}
where $\bbZ/p$ acts trivially on $BA^{\bbZ/p}$. If $r$ is the rank of the finitely generated free abelian group $A^{\bbZ/p}$, then
$$\ker\left(\res_{\bbZ/p}^{\{1\}} \colon  \pi^{n}_{\bbZ/p}(BA^{\bbZ/p}) \to \pi^{n}_{s}(BA^{\bbZ/p})\right)
~ = ~  \bigoplus_{k = 0}^r \ker\left(\res_{\bbZ/p}^{\{1\}} \colon \pi^{n-k}_{\bbZ/p} \to \pi^{n-k}_{s}\right)^{\binom{r}{k}}.$$
\end{example}


\typeout{--------- Section 8: The Segal Conjectures for Infinite Groups  -----------}

\section{The Segal Conjecture for Infinite Groups}
\label{sec:The_Segal_Conjecture_for_Infinite_Groups}

We can now formulate a version of the Segal Conjecture for infinite groups.
Let $\epsilon^G \colon \hoA(G) \to \bbZ$ be the ring homomorphism
which sends an element represented by
a morphism $u \colon S^{\xi \oplus \underline{\bbR^k}} \to S^{\xi \oplus \underline{\bbR^k}}$
in $\SPHB^G(\underline{E}G)$ to the mapping degree of the map induced on the fiber
$u_x \colon S^{\xi_x \oplus \bbR^k} \to S^{\xi_x \oplus \bbR^k}$ for some $x \in \underline{E}G$.
This is the same as the composition
$$\hoA(G) \xrightarrow{\edge^G} \invA(G) \xrightarrow{\invcharacter^G}
\prod_{(H) \in \consub_f(G)} \bbZ \xrightarrow{\pr_{\{1\}}} \bbZ,$$
where $\invcharacter^G$ is the ring homomorphism defined in \eqref{character^G_invA(G)_to_prod_{(H),|H|<infty}Z}
and $\pr_{\{1\}}$ the projection onto the factor belonging to the trivial group.
We define the \emph{augmentation ideal} $\bfI_G$ of $\hoA(G)$ to be the kernel of the ring homomorphism $\epsilon^G$.
Recall that for a finite proper $G$-$CW$-complex $X$ the abelian group
$\pi^n_G(X)$ is a $\pi^0_G(X)$-module. The classifying map
$f \colon X \to \underline{E}G$ is unique up to $G$-homotopy. Suppose
that $\underline{E}G$ is finite. Then $f$ induces
a uniquely defined ring homomorphism $\pi_G^0(f) \colon \hoA(G) = \pi_G^0(\underline{E}G)
\to \pi^0_G(X)$ and we can consider $\pi^n_G(X)$ is a $\hoA(G)$-module.

\begin{conjecture}[Segal Conjecture for infinite groups]
\label{con:Segal_Conjecture_for_infinite_groups}
Let $G$ be a group such that there is a finite model for the classifying space of proper $G$-actions
$\underline{E}G$. Then for every finite proper $G$-$CW$-complex there is an isomorphism
$$\pi^n_s(EG \times_G X) \xrightarrow{\cong} \pi^n_G(X)\widehat{_ {\bfI_G}},$$
where $\pi^n_G(X)\widehat{_ {\bfI_G}}$ is the $\bfI_G$-adic completion of the
$\hoA(G)$-module $\pi^n_G(X)$.

In particular we get for all $n \in \bbZ$ an isomorphism
$$\pi^n_s(BG) \xrightarrow{\cong} \pi^n_G(\underline{E}G)\widehat{_ {\bfI_G}}$$
and especially for $n = 0$
$$\pi^0_s(BG) \xrightarrow{\cong} \hoA(G)\widehat{_ {\bfI_G}}.$$
\end{conjecture}

If $G$ is finite,
Conjecture~\ref{con:Segal_Conjecture_for_infinite_groups} reduces
to the classical Segal Conjecture (see
Theorem~\ref{the:Segal_Conjecture_for_finite_groups}).

The classical Segal Conjecture (see
Theorem~\ref{the:Segal_Conjecture_for_finite_groups}) for a finite group
does not say much if the finite group is torsionfree, i.e.\ is trivial.
The same remark holds for the Segal Conjecture~\ref{con:Segal_Conjecture_for_infinite_groups}
for infinite groups if the group under consideration is torsionfree as explained below.

Since $G$ is torsionfree, $G$ acts freely on $X$ and $\underline{E}G = EG$. Thus we obtain an identification
$$\pi^n_s(EG \times_G X) = \pi^n_s(G\backslash X).$$
The Burnside ring $\hoA(G)$ becomes $\pi^0_s(BG)$ and $\bfI_G$ corresponds to the kernel of the
ring homomorphism $\pi^0_s(BG) \to \pi^0_s(\pt)$. By assumption $BG$ is a finite $CW$-complex.
It is not hard to check that $(\bfI_G)^{\dim(BG) + 1} = \{0\}$.
Hence we obtain an identification
$$\pi^n_G(X)\widehat{_ {\bfI_G}} \cong \pi^n_G(X) \cong \pi^n_s(G\backslash X).$$
Under these two identifications the map appearing in
Conjecture~\ref{con:Segal_Conjecture_for_infinite_groups} is the identity
and hence an isomorphism.

\typeout{-------------------- References -------------------------------}

\addcontentsline{toc}{section}{References}
\bibliographystyle{abbrv}
\bibliography{dbdef,dbpub,dbpre,dbburnextra}

\end{document}